\documentclass[12pt,a4paper,final]{amsart}
\usepackage{geometry}\newgeometry{asymmetric, centering}
\usepackage{mathrsfs} % For really curly L
\usepackage[all]{xy}
\usepackage{url} 
\usepackage{amssymb}
\usepackage{graphicx}
\usepackage{psfrag}
\usepackage{color}

\newcommand{\mmtext}{}%{\color{blue}}

\newcommand{\nntext}{}%{\color{blue}} 

\newcommand{\extension}[1]{} 
\newcommand{\generalizations}[1]{}

\DeclareMathOperator{\singsupp}{singsupp}

\DeclareMathOperator{\Id}{Id}
\DeclareMathOperator{\Symop}{Sym}

\DeclareMathOperator{\Ric}{Ric}
\DeclareMathOperator{\Ein}{Ein}
\renewcommand{\div}{\hbox{div}}

\def\vec{\pmb}

\def\T{{\mathcal T}}
  
\renewcommand{\H}{{\mathbb H}} 
\newcommand{\R}{{\mathbb R}} 
\newcommand{\D}{{\cal D}} 
\newcommand{\I}{{\cal I}}

\newcommand{\N}{{\mathbb N}} 
\newcommand{\cal}{\mathcal } 
 
\renewcommand{\L}{{\mathcal L}}

\def \p {\partial}

\def\F{{\mathcal F}}

\def\hat{\widehat}
\def\tilde{\widetilde}

\def \ba {\begin {eqnarray*} }
\def \ea {\end {eqnarray*} }
\def \beq {\begin {eqnarray} }
\def \eeq {\end {eqnarray}}

\def \supp {\hbox{supp}\,}

\def \WF {\hbox{WF}\,}

\newtheorem{definition}{Definition}[section] 
\newtheorem{theorem}[definition]{Theorem} 
\newtheorem{lemma}[definition]{Lemma} 
\newtheorem{proposition}[definition]{Proposition} 
\newtheorem{corollary}[definition]{Corollary}

\theoremstyle{remark} 
\newtheorem{remark}[definition]{Remark}

\usepackage[notref, notcite]{showkeys}
\usepackage[bordercolor=white, color=white, colorinlistoftodos]{todonotes}

\makeatletter\providecommand\@dotsep{5}\def\listtodoname{List of Todos}\def\listoftodos{\hypersetup{linkcolor=black}\@starttoc{tdo}\listtodoname\hypersetup{linkcolor=blue}}\makeatother
\newcommand{\norm}[1]{\left\|#1 \right\|}

\DeclareMathOperator{\tr}{tr}
\DeclareMathOperator{\Hom}{Hom}
\DeclareMathOperator{\End}{End}

\def\tmax{\hat T}
\def\Mmax{\hat M}

\def\Sym{\mathrm{Sym}^2}

\title[Einstein-scalar field equations]
{Inverse problem for Einstein-scalar field equations
}
\date{Jan. 5, 2018}

\author[Y. Kurylev]{Yaroslav Kurylev}
\address{Department of Mathematics, University College London, 
Gower Street, London UK, WC1E 6BT.}
\email{y.kurylev@ucl.ac.uk}

\author[M. Lassas]{Matti Lassas} 
\address{Department of Mathematics and Statistics, University of Helsinki, Box 68, Helsinki, 00014, Finland}
\email{Matti.Lassas@helsinki.fi}

\author[L. Oksanen]{Lauri Oksanen}
\address{Department of Mathematics, University College London, 
Gower Street, London UK, WC1E 6BT.}
\email{l.oksanen@ucl.ac.uk}
\author[G. Uhlmann]{Gunther Uhlmann}
\address{
Department of Mathematics, University of Washington, 
Box 354350
Seattle, Washington 98195,
USA}
\address{Institute for Advanced Study, Hong Kong University of Science and Technology, Hong Kong SAR
}
\email{gunther@math.washington.edu}


\begin{document}
\begin{abstract}
The paper introduces a method to solve inverse problems for hyperbolic systems where the leading order terms are non-linear. We apply the method
to the coupled Einstein-scalar field equations and
study the question whether the structure of spacetime can be determined by making active measurements near the world line of an observer. We show that such measurements determine the topological, differential and conformal structure of the spacetime in the optimal chronological diamond type set containing the world line. 
In the case when the unknown part of the spacetime is vacuum,
we can also determine the metric itself. We exploit the non-linearity of the equation to obtain a rich set of propagating singularities, produced by a non-linear interaction of singularities that propagate initially as for linear wave equations. This non-linear effect is 
then used as a tool to solve the inverse problem for the non-linear system.
The method works even in cases where the corresponding inverse problems for linear equations remain open, and it can potentially be applied to a large class of inverse problems for non-linear hyperbolic equations encountered in practical imaging problems.
\end{abstract}

\maketitle
\tableofcontents

\section{Introduction and main results}
\label{sec_intro}

In this paper we introduce a new method to solve inverse problems for hyperbolic systems where the leading order terms 
depend non-linearly on the solution.
Earlier, we have considered a related inverse problem
for scalar wave equations where the leading order terms are given by the linear wave operator with a fixed metric \cite{KLU-august}. 
Non-linear hyperbolic systems are used in several applications, but
the models of the non-linear terms are often quite heuristic.
In order to have a canonical non-linearity, we
consider a measurement model for Einstein's equations coupled with scalar fields.

In this context several difficulties appear. In particular,
detection of singularities of waves becomes non-trivial: A well-known physical example of this is that 
only two polarisations of gravitational waves can be observed even though 
the solutions to Einstein's equations take values in the space of symmetric $4\times 4$ matrices that corresponds
to a 10-dimensional space of polarisations (see Section \ref{sec_outline} for further discussion). Moreover, the definition of measurements becomes non-trivial as any source changes the space-time itself in the future
of its support.

The method we develop utilises the
non-linearity as a tool and it enables us to solve inverse problems for non-linear equations
even in cases where the inverse problem for the corresponding linear system remains open.
Indeed, the existing uniqueness results for linear hyperbolic equations require restrictive geometric assumptions such as stationarity, whereas our result is applicable to any globally hyperbolic spacetime. We review the existing literature on inverse problems in detail in Section \ref{sec_literature}.

{\nntext
Outside the context of inverse problems,
there are many results on non-linear interaction of waves,  starting from the studies of
 Bony \cite{Bony}, Melrose and Ritter \cite{MR1,MR2} and Rauch and Reed,  \cite{R-R}.
However, these studies are different from the present paper
since they assume that the geometrical setting, in particular the locations and types of caustics, is a priori known, whereas in our inverse problem wave propagation is studied on an unknown manifold and we allow for caustics of arbitrary, including unstable, type.
This causes further difficulties in the analysis of the non-linear interaction.}

In physical terms, we study the question: Can an observer determine the structure of the surrounding spacetime by doing measurements near its world line? 
The conservation law for Einstein's equations dictates, roughly speaking, that any source in the equation must take energy from some fields in order to increase energy in other fields.
The scalar fields are included in the model that we consider to facilitate this. They correspond to spin zero particles.

Let us now give the mathematical formulation of the problem.
Let $M$ be a smooth $1+3$-dimensional manifold. 
Einstein's equations for a Lorentzian metric $g$ on $M$ are
\begin{equation*}
\Ein(g) = T,
\end{equation*} 
where $\Ein(g)$ is the Einstein tensor associated to $g$ and $T$ is a stress-energy tensor. In vacuum $T = 0$. 
The measurement model we consider is 
\begin{align}
\label{eq1}
& \Ein(g) = T, \quad T = \mathbb T(g, \phi) + \mathcal{F}^1,
\\&\label{eq2} 
\square_g \phi_l - \mathcal V_l'(\phi) = \mathcal{F}^2_l, \quad l = 1,2,...,L. \end{align}
Here $\mathcal V_l'(\phi) = \p_{\phi_l} \mathcal V(\phi)$,  $\phi = (\phi_l)$, $l=1,2,\dots,L$, are the scalar fields,
and 
$$
\F = (\F^1,\F^2) = (\F^1, \F^2_1,\dots, \F^2_L),
$$
models a source in the measurements.
The physical interpretation of $\mathcal F$ is discussed in Appendix C.
The coupling $\mathbb T = (\mathbb T_{jk})$
is given by
\begin{equation*} \mathbb T_{jk}(g, \phi) 
= 
\sum_{l=1}^L \left( \partial_j \phi_l \partial_k \phi_l - \frac{1}{2} g_{jk} g^{pq} \partial_p \phi_l \partial_q \phi_l\right) - \mathcal V(\phi) g_{jk} ,
\end{equation*}
the potential $\mathcal V$ is assumed to be a smooth function $\R^L \to \R$, and 
$\square_g$ is the wave operator associated to $g$.
For example, a typical model is 
$$
\mathcal V(\phi) = \frac 1 2 \sum_{l = 1}^L m_l^2 \phi_l^2,
$$
where $m_l \ge 0$ are constants. 
 
We say that a smooth Lorentzian manifold $(M,\hat g)$ and a function $\hat \phi$ are background spacetime and scalar fields if they satisfy (\ref{eq1})--(\ref{eq2}) with $\mathcal F = 0$ and if $(M, \hat g)$ is globally hyperbolic.
The definition of global hyperbolicity is recalled below, see Definition
\ref{def_gl_hyperb}. 
A globally hyperbolic manifold is isometric to a product manifold $\mathbb{R} \times N$ with the Lorentzian metric given by
\begin{equation*}
%\label{gl_hyp_form}
\hat g = -\beta(t,y) dt^2 + \kappa(t,y),
\end{equation*} 
where $\beta : \mathbb{R} \times N \rightarrow \mathbb{R}^+$ is smooth and $\kappa$ is a Riemannian metric on $N$ depending smoothly on $t$, see \cite{Bernal2005}.
Without loss of generality, we make the standing assumption that 
$M = \R \times N$, and consider equations (\ref{eq1})--(\ref{eq2}) with the initial conditions 
\begin{equation}
\label{init_cond}
g = \widehat{g}, \; \phi  = \widehat{\phi}, \quad \text{in $(-\infty,0) \times N$.}
\end{equation}

The source $\mathcal F$ in (\ref{eq1})--(\ref{eq2}) can not be arbitrary since the Bianchi identities 
imply that that the $\Ein(g)$ is divergence free, that is, $\div_g \Ein(g) = 0$. Hence
the stress energy tensor $T$ needs to satisfy the conservation law
\begin{equation}
\label{conservation_law}
\div_g T = 0.
\end{equation}
This again implies the compatibility condition
\begin{equation}
\label{compcond}
\div_g \mathcal{F}^1 + \sum_{l=1}^L \mathcal{F}^2_l d \phi_l = 0,
\end{equation} 
see Corollary \ref{cor_reduced} for a detailed discussion.
In order to guarantee that (\ref{eq1})--(\ref{init_cond}) admits a rich set of sources in a neighbourhood of $\mathcal F = 0$, we assume the following generic condition, that we call the {\em non-degeneracy condition} for the background scalar fields,
\begin{itemize}
\item[(ND)] $d \hat{\phi}_l(x)$, $l=1,2,3,4$, are linearly independent for all $x \in U$,
\end{itemize}
where $U \subset M$ is a set containing the support of the source $\mathcal F$ in (\ref{eq1})--(\ref{eq2}).

\begin{figure}%[tbhp]
\centering
\def\svgwidth{.4\linewidth}
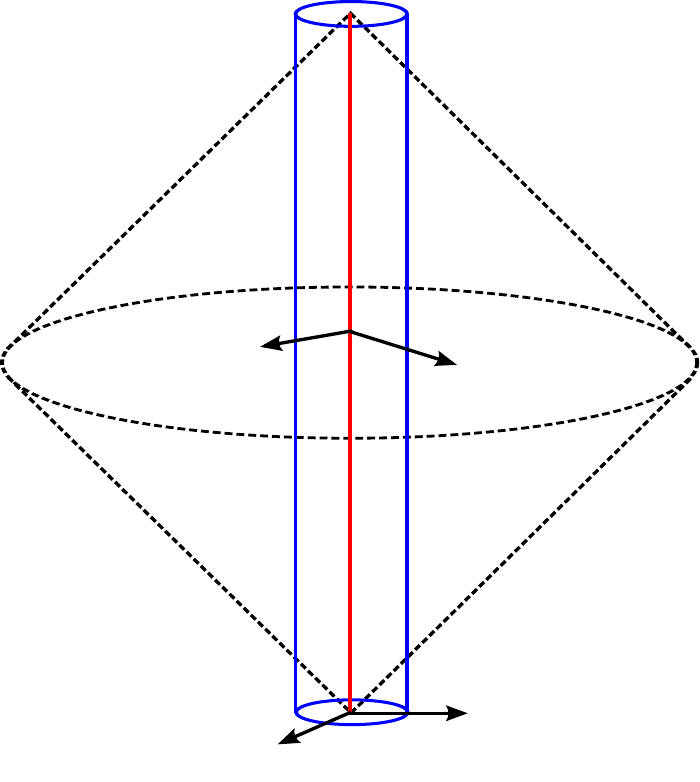
\def\svgwidth{.4\linewidth}
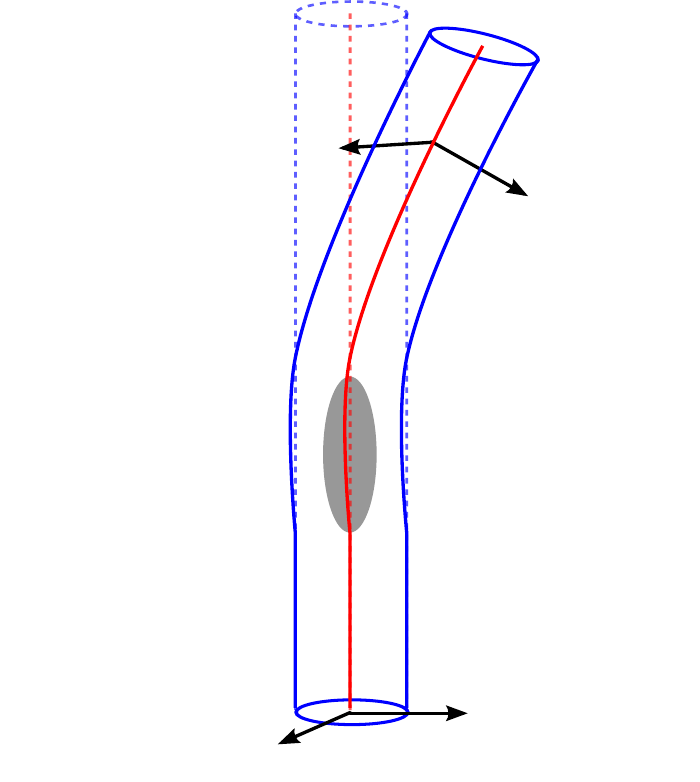
\caption{
\label{fig_schematic_active}
{\em Left.} The setting of Theorem \ref{th_main}.
The solid red line depicts the worldline $\mu$
along which the frame $Y_j = Y_j^{\hat g}$, $j=1,2,3$, defines the Fermi coordinates (the third direction is suppressed in the picture). 
The blue cylinder is the coordinate neighbourhood $\Phi(V)$
and the dashed double cone is the chronological diamond $I(\mu(0), \mu(1))$.
{\em Right.} Fermi coordinates for a perturbed metric $g$.
As the source $\mathcal F$ becomes non-zero, the Fermi coordinates start to change. The support of $\mathcal F$ is depicted as the grey area. 
}
\end{figure}

Observe that the system (\ref{eq1})--(\ref{init_cond})
is invariant with respect to isometries and that 
the metric tensor $g$ begins to change as soon as $\mathcal F$ becomes non-zero. 
In order to take this invariance into account, we model local measurements using Fermi coordinates. 
Let $x \in \{0\} \times N$ and 
let $\xi \in T_x M$ be timelike and future-pointing. Define $\mu_{\hat g}(s) = \exp^{\hat g}_{x}(s \xi)$,
where $\exp^{\hat g}$ is the exponential map on $(M,\widehat g)$. 
Let $X_j$, $j=0,1,2,3$, be a basis of $T_{x}M$, with $X_0 = \xi$,
and consider the following Fermi coordinates around the world line $\mu_{\hat g}$,
\begin{align}\label{def_Fermi}
\Phi_{\hat g}(s, y^1, y^2, y^3) = \exp^{\hat g}_{\mu_{\hat g}(s)}\left(\sum_{j=1}^3 y^j Y_j^{\hat g}\right),
% \quad \Phi_{\hat g} : V \to M,
\end{align}
where $Y_j^{\hat g}$ is the parallel transport of $X_j$ along $\mu_{\hat g}$, $j=1,2,3$. Here the parallel transport is taken with respect to $\widehat g$, and 
$$
s \in [0,1], \quad (y^1,y^2,y^3) \in \overline{B_0},
$$
where $B_0 \subset \R^3$ is a small enough neighbourhood of the origin so that the Fermi coordinates (\ref{def_Fermi}) are well-defined in $[0,1] \times \overline{B_0}$.
We write 
$$
V = (0,1) \times B_0.
$$

Let $\tmax > 0$ satisfy $\Phi_{\hat g}(\overline V) \subset (-\infty,\tmax) \times N$, and consider a Lorentzian metric $g$ on $(-\infty, \tmax) \times N$
such that the corresponding Fermi coordinates $\Phi_g : V \to M$ are well-defined. Here $\Phi_g$ is defined as $\Phi_{\hat g}$ above except that now the world line $\mu_g$, the exponential map and the parallel transports are defined with respect to $g$.
We often use the shorthand notation $\mu = \mu_{\hat g}$, $\Phi = \Phi_{\hat g}$ and
\begin{align}\label{def_M_strip}
\Mmax = (0, \tmax) \times N.
\end{align}
Then we model local measurements near the world line $\mu$ by the following data set 
\begin{align}\label{data}
\mathcal D_{r, k}(\hat g, \hat \phi) = &\{(\Phi_g^* g, \Phi_g^* \phi, \Phi_g^* \mathcal F) \in C^k(V)\;;\; 
\\\notag&\qquad
\text{$(g, \phi, \mathcal F)$ satisfies 
(\ref{eq1})--(\ref{eq2}) on $\Mmax$},\ \text{(\ref{init_cond}) holds},
\\\notag&\qquad
\supp(\mathcal F) \subset \Phi_g(V),\ \norm{\Phi_g^* \mathcal F}_{C^k(V)} < r \},
\end{align}
where $r > 0$, $k \in \N$, and $\Phi_g^*$ is the pullback under $\Phi_g : V \to M$.
This set corresponds to the graph of the response operator used in earlier studies, see (\ref{resp_op}) below.

We emphasize that in the above model, the representations of the fields $g$ and $\phi$ in the Fermi coordinates are measured only in a small neighbourhood $\Phi_g(V)$ of the world line $\mu_{g}$. The neighbourhood $\Phi_{\hat g}(V)$ is depicted as the blue cylinder Figure \ref{fig_schematic_active}.
Let us also remark that $(g,\phi,\mathcal F)$ in the defining property of 
$\mathcal D_{r, k}(\hat g, \hat \phi)$ satisfies (\ref{compcond}) since it satisfies (\ref{eq1})--(\ref{eq2}).
In particular, 
for fixed $(g,\phi)$,
the sources $\mathcal F$ in (\ref{data}) depend on $(g,\phi)$ via the constraint (\ref{compcond}).

As signals in Einstein's equations propagate with finite speed,
the data $\mathcal D_{r, k}(\hat g, \hat \phi)$ can not determine $\hat g$ on the whole manifold $M$. 
The fastest signals, for example gravitational waves, propagate at the speed of light, and 
an observer moving along $\mu$ is able to probe a region $U$ in $M$ only if signals can reach $U$ from $\mu$ and if signals from $U$ can return to $\mu$.

We will define next the set where we are able to gain information on $\hat g$. This diamond type set is close to optimal. Indeed, if the chronological future and past are replaced by the causal future and past in its definition, then 
it is not possible to obtain information on $\hat g$
in the complement of this slightly larger set, independently of the choice of the Fermi coordinates.
We give the definitions of causal future and past in Section \ref{sec_notations} below, and recall here only the concepts that are needed to formulate our main result.
For $p, q \in M$, $p \ll q$ in $(M, \hat g)$ means that $p$ and $q$ can be joined by future pointing timelike curves in $(M,\hat g)$. 
The \emph{chronological future} and \emph{past} of $p \in M$ are
\begin{equation}
\label{future_I}
I^+(p) = \{q \in M \; ; \; p \ll q \text{ in $(M,\hat g)$} \}, 
\quad
I^-(p) = \{q \in M \; ; \; q \ll p \text{ in $(M,\hat g)$} \}.
\end{equation}
We will obtain information on $\hat g$ on the chronological diamond type set 
\begin{equation}
\label{diamond}
I(p,q) = I^+(p) \cap I^-(q),
\end{equation}
containing the world line $\mu$ in the sense that 
$p = \mu(0)$ and $q = \mu(1)$.

\begin{theorem} 
\label{th_main}
\label{alternative main thm Einstein}
Suppose that smooth $\hat g$ and $\hat \phi$ satisfy (\ref{eq1})--(\ref{eq2}) 
on $M$
with $\mathcal F = 0$, and that $(M, \hat g)$ is globally hyperbolic.
Consider the Fermi coordinates $\Phi : V \to M$ around a world line $\mu$ as defined in (\ref{def_Fermi}), and 
suppose that $\hat{\phi}$ satisfies (ND) with $U = \overline{\Phi(V)}$.
Then for arbitrary small $r > 0$ and large $k \in \N$,
the data set ${\mathcal D}_{r, k}(\hat g, \hat \phi)$ determines the topology, differentiable structure and conformal class of the metric $\widehat g$ in the chronological diamond $I(\mu(0), \mu(1))$. 
\end{theorem}

In the case of a vacuum spacetime we can recover also the conformal factor. 

\begin{corollary}
\label{cor_to_main}
Let $\hat g$, $\hat \phi$, $\Phi$ and $\mu$ be as in 
Theorem \ref{alternative main thm Einstein}.
Suppose, furthermore, that $W\subset I(\mu(0), \mu(1))$ is vacuum, i.e. Ricci-flat,
and that any point in $W$ 
can be connected by a curve, lying completely in $W$, to a point in $\Phi(V)$.
Then for arbitrary small $r > 0$ and large $k \in \N$,
the data set ${\mathcal D}_{r, k}(\hat g, \hat \phi)$ determines the topology and differentiable structure of $W$ and isometry class of $\widehat g$ in $W$. 
\end{corollary}

%\HOX{Order of remarks is changed.}
\begin{remark}
The techniques considered in this paper can be used also to study
inverse problems for non-linear hyperbolic systems encountered in applications.
{\nntext In fact, one of the  motivations of this paper is to develop methods that will transfer 
tools of General Relativity to medical and seismic imaging.} For instance, in medical imaging, in the
the recently developed Ultrasound Elastography imaging technique
the elastic material parameters are reconstructed by sending 
(s-polarized) elastic waves that are imaged using (p-polarized) elastic waves,
see  
e.g.\ \cite{Hoskins2012,McLaughlin2006}.
This imaging method uses  interaction of waves and is based on the non-linearity of the system.
\end{remark}

\begin{remark}
{\nntext The question, what is the topology of our Universe, can be considered as an inverse problem for the coupled Einstein and matter field equations. This fundamental question has been studied in \cite{Luminet,Ringstrom2}. 
In addition to the topology, in this paper we consider also the determination of the conformal type of the Lorentzian metric.}
\end{remark}

\begin{remark}
Theorem \ref{alternative main thm Einstein} implies that, if we have
two non-conformal  spacetimes, a generic  measurement gives different results
on these manifolds. In particular, this implies that perfect 
spacetime cloaking, in sense of light rays, see \cite{Fridman2012,McCall2011},  is not possible with a smooth metric  in a globally hyperbolic universe.
\end{remark}

\begin{remark}
The non-degeneracy assumption (ND) could be replaced by the assumption that 
$d \hat{\phi}_l(x)$, $l=1,\dots,L$, span $T_x^* M$ for all $x \in \overline{\Phi(V)}$. We make the stronger assumption just to avoid cumbersome notation in the linearization stability proofs below. 
\end{remark}

\subsection{Notations and definitions}
\label{sec_notations}

Analogously to the chronological future and past (\ref{future_I}),
we define the causal future and past of a point $p \in M$ 
on a Lorentzian manifold $(M,g)$
by 
\begin{equation*}
%\label{future_J}
J_{(M,g)}^+(p) = \{q \in M \; ; \; p \le q \text{ in $(M,g)$} \},\quad 
J_{(M,g)}^-(p) = \{q \in M \; ; \; q \le p \text{ in $(M,g)$} \}.
\end{equation*}
where $p \le q$ means that $p$ and $q$ can be joined by a future pointing causal curve in $(M,g)$ or $p=q$.
We write also $J_{(M,g)}(p,q) = J_{(M,g)}^+(p) \cap J_{(M,g)}^-(q)$.

\begin{definition}[\cite{Bernal2007}]
\label{def_gl_hyperb}
A Lorentzian manifold $(M, g)$ is \emph{globally hyperbolic} if 
there are no closed causal paths in $M$, and 
the set $J_{(M,g)}(p,q)$ is compact for any pair of points $p,q \in M$.
\end{definition}

We denote by $\nabla_g$, $\div_g$ and $\Ric(g)$ the covariant derivative, divergence and the Ricci tensor with respect to $g$. 
Using the shorthand notation 
$$
(d\phi \otimes d\phi)_{jk} = \sum_{\ell = 1}^L \p_j \phi_\ell \p_k \phi_\ell, \quad
(\mathcal V'(\phi))_\ell = \p_{\phi_\ell} \mathcal V(\phi).
$$
where $j,k=0,1,2,3$, $\ell=1,\dots,L$,
equation 
(\ref{eq1})
can be written as 
\begin{align}
\label{ein1}
\Ric(g) - d\phi \otimes d\phi - \mathcal V(\phi) g &= F^1,
\quad F^1 = I_g \mathcal F^1,
\end{align}
where $I_g$ is the involution 
$$
I_g F^1 = F^1 - \frac{\tr_g F^1}{2} g. 
$$
Here $\tr_g F^1 = g^{jk} F^1_{jk}$ 
is the contraction with $g$.
We will use also the shorthand notation
$$
I_g \mathcal F = (I_g \mathcal F^1, \mathcal F^2).
$$

The bundle of symmetric 2-tensors over $M$
is denoted by $$
\Sym = \Symop(T^*M \otimes T^*M).
$$ For vector bundles $E$ and $F$ over $M$, we denote by $\Hom(E, F)$
the vector bundle whose fibre at $x \in M$ is the space of linear maps from the fibre $E_x$ of $E$ to the fibre $F_x$ of $F$. We write also $\End(E) = \Hom(E,E)$.

\subsection{Earlier results}
\label{sec_literature}

Inverse problems for partial differential equations
is a much studied topic, however, the present theory is largely confined to the case of linear equations. For a majority of the few existing results on non-linear equations, e.g. \cite{Isakov1993,Salo2012,Sun1996}, non-linearity is an obstruction rather than a helpful feature.
On the contrary, our method to solve the inverse problem for Einstein's equations is based on the non-linear interaction of gravitational waves, and only global hyperbolicity is assumed on the underlying spacetime.
In an earlier paper \cite{KLU-august} we considered a similar method in the case of a semi-linear wave equation, and we discuss this result in more detail below. 
We emphasize that the present study of Einstein's equations can not use the techniques in \cite{KLU-august} in a straightforward manner due to the much more complicated structure of the non-linearity. In particular, contrary to \cite{KLU-august}, the non-linearity appears also in the leading order terms. 

{\nntext As our approach uses the non-linearity in an essential way, % \HOX{This sentence is changed}
it differs from the classical methods used to study inverse problems for linear wave models, that is, }
the Boundary Control method 
\cite{Katchalov2001} and
%\cite{Belishev1987,Belishev1992,Katchalov2001,Anderson2004,Eskin2010}
the geometric optics based approach \cite{Rakesh1988a,Stefanov1989}, 
{\nntext that} have been used to solve various inverse problems for linear wave equations. {\nntext However, these} methods fail when the equations have general time-dependent coefficients. 
For the {\nntext former} method this is because the sharp unique continuation result \cite{Tataru1995} does not hold for wave equations with 
coefficients that are
smooth, but not real analytic, in the time variable \cite{Alinhac1983};
and for {\nntext the latter method} because it is not known if 
the light ray transform is injective on a globally hyperbolic Lorentzian manifold.
 {\nntext In the case of time-independent coefficients, the light ray transform reduces to the geodesic ray transform. The injectivity of the latter transform has been studied extensively in the Riemannian context
\cite{Stefanov2,Ilmavirta,Paternain1,Paternain3,Paternain2,Stefanov1,Uhlmann-Vasy}.}

The inverse problem for linearized Einstein's equations is open, and the same is true for other similar linear wave equations. We solve the inverse problem for the physical, non-linear Einstein's equations.

In the case of linear equations, it is now classical to model active measurements by using the Dirichlet-to-Neumann map \cite{Sylvester1987}.
In the case of Einstein's equations, modelling active measurements 
is more subtle since the spacetime can not be observed from outside, and (\ref{data}) appears to be the first model for active measurements in this case.
Let us now recall from \cite{KLU-august} the model of active measurements for the non-linear wave equation
\begin{align}%\label{}
\label{eq: wave-eq general}
&\square_{\hat g}u(x)+a(x)\,u(x)^2=f(x)
\quad\hbox{on }M,\hspace{-.5cm}\\ \nonumber
&u(x)=0,\quad \hbox{for }x\in M_0\setminus J^+_{(M,\hat g)}(\supp(f)),
\end{align}
where $a$ is non-vanishing and smooth, and $(M,\hat g)$ is globally hyperbolic. 
Letting $\mu$ be as in Theorem \ref{th_main}, 
and letting $V \subset M$ be a neighbourhood of $\mu([0,1])$,
the active measurements can be modelled by the set
\begin{align}\label{data_wave}
\mathcal D_{r, k}^{\text{wave}}(\hat g) = \{(u|_V, f) \;;\; 
&\text{$(u,f)$ satisfies 
(\ref{eq: wave-eq general}) on $\Mmax$},\ 
\\\notag&
\supp(f) \subset V,\ \norm{f}_{C^k(V)} < r \}.
\end{align}
The set (\ref{data_wave}) is the analogue of (\ref{data}) but simpler, since the complications related to the coordinate invariance of the Einstein's equations do not arise in this case. 
Let us remark that in \cite{KLU-august} the results were actually formulated by using the response operator 
    \begin{equation}\label{resp_op}
L_V f = u|_V
    \end{equation}
For small enough $r > 0$ and large enough $k \in \N$, the set (\ref{data_wave}) is the graph of $L_V$. 
It was shown in \cite{KLU-august} that the set (\ref{data_wave})
determines the conformal class of the metric $\hat g$ in $I(\mu(0),\mu(1))$.
This is the analogue of Theorem \ref{th_main} for (\ref{eq: wave-eq general}).

\begin{remark}
The result in \cite{KLU-august} deals with an inverse problem for ``near field'' measurements, modelled by  (\ref{data_wave}) or the corresponding response operator.
In the case of inverse problems for linear equations, near field
measurements are typically modelled by the Dirichlet-to-Neumann map, and they are
equivalent to scattering measurements or ``far field'' information \cite{Berezanski1958}. Analogous equivalence
for non-linear equations has not yet been studied but it is plausible that the theory of near field measurements is useful when studying inverse scattering problems also in this case. On related inverse scattering
problems, {\nntext see \cite{Graham2003,JSb,Melrose2014,SB}.

Finally, even though in this paper we consider near-filed measurements, the techniques we use to
study non-linear interaction of waves takes its inspiration from microlocal scattering theory
\cite{BSb,VBW,Guillarmou2,HMV,H-V2,H-V1,MW,V1}. In the case of the near field measurements,
the inverse problems for Einstein manifolds have been studied in the Riemannian context in \cite{SB}. 
}
\end{remark}

\section{Outline of the proof}
\label{sec_outline}

\begin{figure}
\centering
\def\svgwidth{.4\linewidth}
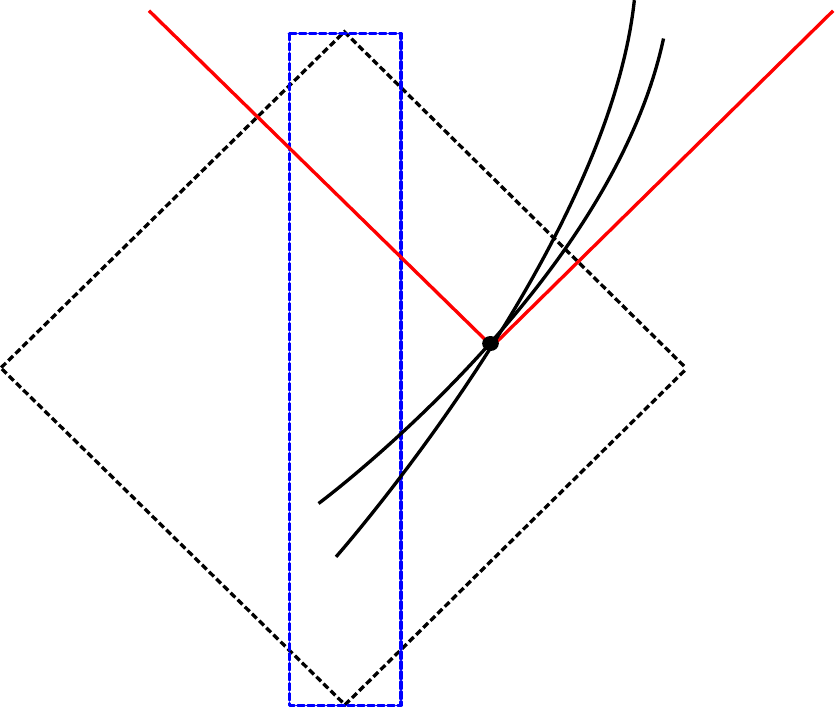
\quad
\includegraphics[width=.4\linewidth]{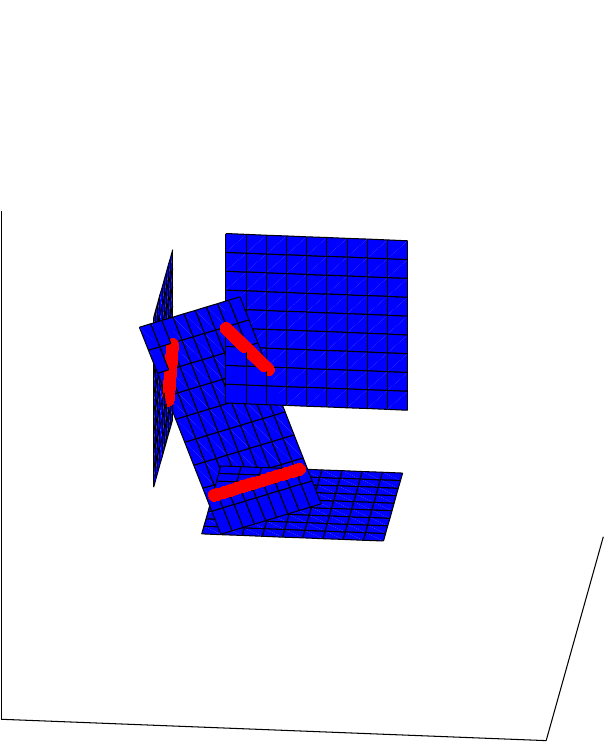}
\\
\includegraphics[width=.4\linewidth]{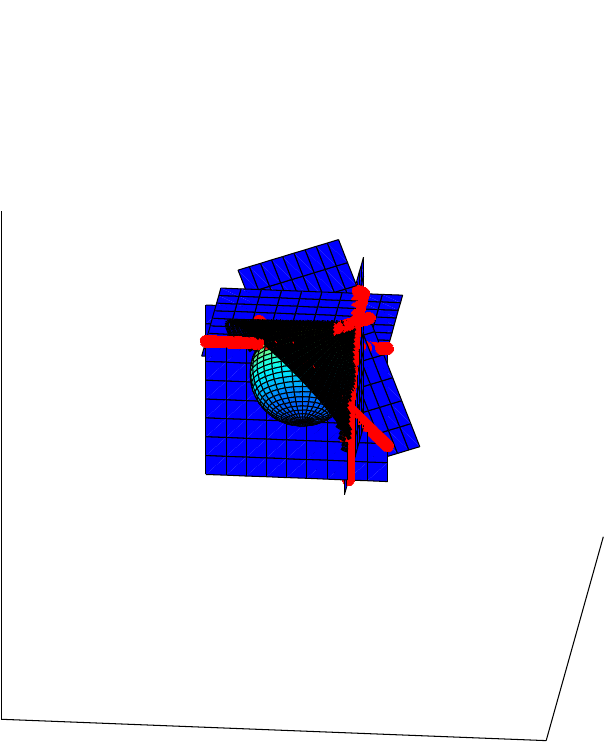}
\quad
\includegraphics[width=.4\linewidth]{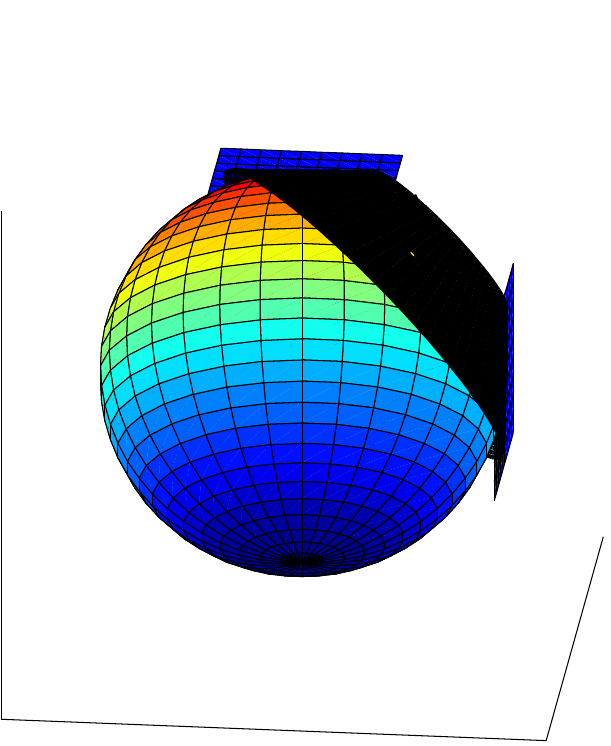}
\caption{
\label{fig_frames}
{\em Top left.} A schematic picture of outline of the proof of Theorem \ref{th_main} as discussed in Section \ref{sec_outline}.
Waves propagating near geodesics (depicted as black curves), sent from $\Phi(V)$,
interact at the point $q \in I(\mu(0), \mu(1))$. 
Via a non-linear interaction, they can produce propagating singularities (in red) analogous to those generated by a point source at $q$. 
{\em Rest of the panels.} Interactions of pieces of the four plane waves $u_j$ as time passes.
{\em Top right.} The pieces interact pairwise on the red line segments, no new propagating singularities appear. 
{\em Bottom left.} Interactions of triples produce conic singularities and the interaction of all the four pieces of waves produces a spherical singularity. Only one of the four conic singularities (in black) is shown in the picture. 
{\em Bottom right.} The spherical singularity expands as time passes. It can be detected far away from the point where the interaction happened. 
%Animation of the full interaction is included as two movies in the {\em supporting information}. In the first movie S1, only one of the four conic singularities is shown,
%and in the second one S2, all the four conic singularities are shown.
}
\end{figure}

In order to give a rough idea of the proof, let us consider here only the simple model case (\ref{eq: wave-eq general}) with $\hat g = -dt^2+dx^2+dy^2+dz^2$, the Minkowski metric, and $a=1$ identically. 
The idea is to use the source $f$ to generate propagating singularities that interact via the non-linearity. The propagating singularities will be constructed so that their interaction produces a new singularity that behaves like a point source. Recording the singularities produced by each such ``artificial point source'' in $I(\mu(0), \mu(1))$ allows for a reduction to the purely geometric problem solved in \cite[Th. 1.2]{KLU-august}.
Let us also mention that the geometric step 
has been generalized for manifolds with boundary in
{\nntext \cite{H-U}. Another generalization involving broken geodesics is given in \cite{Lars}.}

To explain this in more detail, but still roughly, let $f = \epsilon h$, $\epsilon > 0$,
and write the asymptotic expansion of the solution of (\ref{eq: wave-eq general}),
\begin{equation*}
u = \epsilon w_1 + \epsilon^2 w_2 + \epsilon^3 w_3 + 
\epsilon^4 w_ 4 + \mathcal O(\epsilon^5),
\end{equation*} 
where, using the notation $Q = \square_{\hat g}^{-1}$ for the causal inverse,
\begin{align}
\label{asymp_expansion}
&w_1 = Q h, \quad
w_2 = -Q(w_1^2), \quad
w_3 = - 2 Q(w_1 w_2), \quad
\\\notag
&w_4 = -Q(w_2^2) - 2 Q(w_1 w_3).
\end{align}
We say, for example, that $w_3$ results from the interaction of $w_1$ and $w_2$.

Let us make a linear change of coordinates in the Minkowski space, such that in the new coordinates $x^j, j=1,2,3,4$, the hyperplanes $K_j = \{x^j = 0\}$ are lightlike, that is, $T_x K_j$ contains a lightlike vector for all $x \in \R^4$. The plane waves
$u_j(x) = (x^j)_+^m$,
where $m > 0$ and
\begin{align}
\label{linear_wave}
 (x^j)_+^m = \begin{cases}(x^j)^m,  & x^j \geq 0, \\ \ \ \ 0, &\mbox{otherwise,}\end{cases}
\end{align}
are solutions to the wave equation $\square_{\hat g} u = 0$.
They are singular on the hyperplanes $K_j$. Moreover, the singularity on $K_j$ has a specific direction: the wave front set of $u_j$ is the conormal bundle of $K_j$, that is, $\WF(u_j) = N^* K_j \setminus 0$, see Section \ref{sec_a_of_sing} for this notation.

%\begin{figure}%[tbhp]
%\centering
%\includegraphics[width=.4\linewidth]{eplanes001}
%\quad
%\includegraphics[width=.4\linewidth]{eplanes050}
%\\\vspace{-0.2cm}
%\includegraphics[width=.4\linewidth]{eplanes120}
%\quad
%\includegraphics[width=.4\linewidth]{eplanes175}
%\caption{
%\label{fig_frames}
%Interactions of pieces of the four plane waves $u_j$ as time passes.
%{\em Top left.} The pieces approach each other but do not yet interact.
%{\em Top right.} The pieces interact pairwise on the red line segments, no new propagating singularities appear. 
%{\em Bottom left.} Interactions of triples produce conic singularities and the interaction of all the four pieces of waves produces a spherical singularity. Only one of the four conic singularities (in black) is shown in the picture. 
%{\em Bottom right.} The spherical singularity expands as time passes. It can be detected far away from the point where the interaction happened. 
%Animation of the full interaction is included as two movies in the {\em supporting information}. In the first movie S1, only one of the four conic singularities is shown,
%and in the second one S2, all the four conic singularities are shown.
%}
%\end{figure}

Interactions of the waves $u_j$ produce new singularities,
for instance, the product $u_1 u_2$ is singular on $U \cup X$
where $U = N^* K_1 \cup N^* K_2$ and $X = N^* (K_1 \cap K_2)$.
The set $X$ contains covectors that are not in $U$,
however, in this case of the pairwise interaction,
all the lightlike covectors in $X$ are also in $U$.
As the solution operator $Q$ propagates only lightlike singularities,
no new propagating singularities are produced in pairwise interactions.

When three waves interact, say $u_1, u_2, u_3$, new propagating singularities are produced since $N^*(K_1 \cap K_2 \cap K_3)$ contains lightlike covectors that are not in any $N^* K_j$, $j=1,2,3$.
%In the $1+2$-dimensional case, 
The interaction of three waves has been analysed in \cite{Rauch1982,Melrose1985,Melrose1996}.
A conic singularity comes out of the line of intersection, 
and the geometry of the singularity is similar
to shock discontinuities that happen at the interface between subsonic and supersonic speeds, see Figure \ref{fig_frames}.

The intersection of the four hyperplanes $K_j$ is the origin, and its conormal bundle is the fibre of the cotangent space at the origin.
In particular, the conormal bundle contains the light-cone 
$L_0^+ \R^4$, and hence
the interaction of the four waves $u_j$ produces a new singularity that corresponds to a spherical wave, in the sense that the singularity coincides with that generated by a point source at the origin. Summarizing, the four wave interactions can produce artificial point sources.

The proof of Theorem \ref{th_main} is based on an analysis of the non-linear interaction of four solutions to the linearized Einstein-scalar field equations. Of course, 
as the background metric $\hat g$ can be very general, and as the sources $\mathcal F$ in (\ref{data}) need to satisfy support constraints, we can not use explicit solutions like (\ref{linear_wave}). These are replaced by more general conormal distributions in the actual proof below. 

We begin by showing that the Einstein-scalar field equations admit a rich set of solutions with propagating singularities. That is, in Section \ref{sec_construction_singsol} we show that (\ref{eq1})--(\ref{eq2}) can be linearized along families of solutions so that the thus obtained solutions to the corresponding linearized equations have propagating singularities, and that the geometric structure of the singularities is constrained only by a linearized, microlocal version of the conservation law (\ref{conservation_law}), see Definition \ref{def_SY_comp}.

In Section \ref{sec_interaction} we study the analogue of the interaction terms (\ref{asymp_expansion}) for the Einstein-scalar field equations. Although on the level of the geometry of propagation of singularities the analysis of the interaction terms is very similar with the above model case, it is more complicated to show that interaction terms do not vanish, say, due to them cancelling each other. 
We show this by considering the principal symbols of the interaction terms, viewed as conormal distributions. 
Let us also point out that the analysis of the principal symbols is not analogous to that in \cite{KLU-august}, in particular, the asymptotic approach in \cite[Section 3.4]{KLU-august} does not seem to work in our setting. 

Finally in Section \ref{sec_reduction_to_los} we deal with the global geometric aspects of the proof and finish the reduction to \cite[Th. 1.2]{KLU-august}. Also 
we treat there the complications related to the coordinate invariance of the Einstein's equations, that are manifest both when sending and receiving propagating singularities.

To get a flavour of the complications when trying to detect singularities in $g$, say gravitational waves, 
observe that the Fermi coordinates $\Phi_g$ in the definition of the data set (\ref{data}) may be non-smooth when $g$ is non-smooth. 
The use of non-smooth coordinates $\Phi_g$ can hide singularities in $g$. However, 
if $g$ is a conormal distribution on a submanifold of codimension one, certain singularities are generically visible in the coordinates $\Phi_g$.
More precisely,  the principal symbol of $\Phi_g^* g$ is non-vanishing generically on a rank 2 subbundle of the rank 10 bundle of symmetric 2-tensors.
Indeed, as the Fermi coordinates are used, $\Phi_g^* g$ has 6 degrees of freedom to begin with, and on top of this, the conservation law (\ref{conservation_law}) imposes 4 additional constraints. 
Definition \ref{def_SY_comp} below gives the analogue (\ref{conservation_law}) for principal symbols.
In physical terms, if $g$ is a gravitational wave,
the rank 2 subbundle corresponds to the two polarizations of that can be observed, {\nntext see
 \cite {Maggiore}. 
 More detailed physical explanations on the used measurement model can be found from the earlier extended preprint version of this paper  \cite{preprint}.}

%\begin{center}
%$  $\hspace{-10cm}
%\includegraphics[height=4.0cm]{eplanes001.eps}\quad \includegraphics[height=4.0cm]{eplanes050.eps}\quad\includegraphics[height=4.0cm]{eplanes120.eps}\quad\includegraphics[height=4.0cm]{eplanes175.eps}\hspace{-10cm}$  $
%\end{center}
%{\it FIGURE 2. Four plane waves propagate in space. 
%When the planes intersect, the non-linearity of the hyperbolic system  produces
%new waves.
%The four figures show
%the waves before the interaction of the waves start, when 2-wave interactions have  started, when all  4 waves have just interacted, and later
%after the interaction. {\bf Left:} Plane waves before  interacting.
% {\bf Middle left:}
% The 2-wave interactions (red line segments) appear  but do not cause
%new propagating singularities. {\bf Middle right and Right:}  All plane waves have intersected and new 
% waves have appeared. The 3-wave interactions cause 
%new conic waves (black surface). Only one such wave is shown in the figure. The 4-wave interaction causes
%a point source in space-time that sends a spherical wave in all future light-like 
%directions. This spherical wave is essential in our considerations. For an animation
%on these interactions, see the supplementary video \cite{Video}.}
%\medskip

\section{Construction of solutions with propagating singularities}
\label{sec_construction_singsol}

In this section we will show that 
Einstein's equations coupled with scalar fields 
(\ref{eq1})--(\ref{eq2}) 
admit a rich set of solutions with propagating singularities. 
These solutions are constructed for sources $\mathcal F$ having small support. 
The locally generated propagating singularities 
are then used in later sections to probe the background spacetime $(M,\hat g)$.

\subsection{Einstein's equations as a system of quasilinear wave equations}

In order to treat (\ref{eq1})--(\ref{eq2}) as 
a hyperbolic system, we need to fix a suitable gauge. 
We follow the presentation in \cite[pp. 152-154]{Ringstrom}.
Let $\hat g$ be a smooth Lorentzian metric tensor on $M$.
Let $\hat \Gamma_{pq}^r$ be the Christoffel symbols of $\hat g$
and define $\hat H_k = g_{kr} g^{pq} \hat \Gamma_{pq}^r$.
Define also 
$$
\Gamma_{pkq} = \frac 1 2 ( \p_p g_{kq} + \p_q g_{pk} - \p_k g_{pq} ),
\quad \Gamma_k = g^{pq} \Gamma_{pkq}.
$$
As the difference of two connections is a tensor, we see that the local expression $\hat H_k - \Gamma_k$ defines a tensor, that we denote by $H_{\hat g}(g)$. That is, 
\begin{align}\label{def_harm_tensor}
H_{\hat g}(g)_k = \hat H_k - \Gamma_k, \quad k=0,1,2,3.
\end{align}
We define the reduced Ricci tensor
$$
\Ric_{\hat g}(g) = \Ric(g) + \nabla_g H_{\hat g}(g),
$$
%where $\nabla$ is the covariant derivative with respect to $g$,
and consider the reduced Einstein's equations coupled with scalar fields
\begin{align}
\label{einred_full}
\Ric_{\hat g}(g) - d\phi \otimes d\phi - \mathcal V(\phi) g = F^1&,
\quad \text{in $\hat M$,}
\\\notag
\Box_g \phi - \mathcal V'(\phi) = F^2&,
\quad \text{in $\hat M$,}
\\\notag
g = \widehat{g}, \ \phi  = \widehat{\phi}&, \quad \text{in $(-\infty,0) \times N$.}
\end{align}
This is a quasilinear wave equation. Indeed, locally \cite[eq. (14.3)]{Ringstrom}, 
$$
\Ric_{\hat g}(g)_{jk} = - \frac 1 2 g^{pq} \p_p \p_q g_{jk}
+ \nabla_j \hat H_k + g^{pq} g^{rs} (\Gamma_{prj} \Gamma_{qsk}
+ \Gamma_{prj} \Gamma_{qks} + \Gamma_{prk} \Gamma_{qjs}),
$$
where only the first term contains second order derivatives of $g$.
Here $\nabla_j$ is the covariant derivative $\nabla_g$ in the coordinate direction $\p_j$.

It follows from a typical fixed point argument
that (\ref{einred_full}) has a unique solution $(g,\phi)$ when $(F^1,F^2)$
is in a small neighbourhood of the origin in $C_0^k(U; \Sym \oplus \R^L)$, with $U \subset \Mmax$ bounded and $k \in \N$ large, see \cite{HKM} or \cite{Kato1975} for such an argument and 
Appendix B for further discussion.

Let $(g,\phi)$ solve (\ref{einred_full}) and set 
$\mathcal F = I_g F$ where $F =(F^1,F^2)$.
Then $(g,\phi, \mathcal F)$ solves (\ref{eq1})--(\ref{eq2}) 
if and only if the following gauge condition holds 
\begin{align}\label{gauge_cond}
H_{\hat g}(g) = 0.
\end{align}
This is equivalent with the compatibility condition (\ref{compcond}), see Corollary \ref{cor_reduced} in Appendix A. 
In the next section we will give a construction of sources $\F$ that 
satisfy (\ref{compcond}),
and that generate a rich set of propagating singularities,
as demonstrated later in Section \ref{sec_muloc_linstab}.

It will be occasionally useful that the subprincipal part of the second equation in (\ref{einred_full}) simplifies under the condition (\ref{gauge_cond}).
Indeed, in local coordinates (\ref{gauge_cond})
is equivalent with $g^{pq} \Gamma^j_{pq} = g^{pq} \hat\Gamma^j_{pq}$, $j=0,1,2,3$, and then 
\begin{align}\label{Box_phi_subprin}
\Box_g \phi = g^{-1}(D,D) \phi - g^{pq} \hat \Gamma^j_{pq} \p_{x^j} \phi,
\end{align}
where $g^{-1}(D,D) \phi = g^{pq} \p_p \p_q \phi$.
In particular, the subprincipal part of the second equation in (\ref{einred_full}) does not contain derivatives of $g$.

\subsection{Linearization stability}

We will next formulate a linearization stability result for the system (\ref{eq1})--(\ref{init_cond}).
Linearization stability has been studied
extensively in the context of the initial value problem for (\ref{eq1})--(\ref{eq2}), i.e. when $\mathcal{F} = 0$, see e.g. 
\cite{Brill1987,Brill1973,Choquet-Bruhat1973,Fischer1973,Girbau2010}.
The distinctive feature of our result is that the source $\mathcal{F}$
can be non-zero but its support is localized. 
The localization is essential when considering inverse problems, in particular, it allows us to satisfy the constraint $\supp(\mathcal F) \subset \Phi_g(V)$ in (\ref{data}).

\begin{proposition}
\label{prop_lin_stab_basic}
Let $U \subset \hat M$ be open, and suppose that 
$L \ge 4$
and that $\hat \phi$ satisfies the non-degeneracy condition (ND).
If $(f^1,f^2) \in C^k(\hat M; \Sym \oplus \R^L)$ is compactly supported in $U$, $k$ is large enough,
and if 
\begin{align}
\label{compcond_lin}
\div_{\hat g} f^1 + f^2 \cdot d \hat \phi = 0,
\end{align}
then there is a neighbourhood $\mathcal E \subset \R$ of the origin and a family of solutions  
$(g(\epsilon), \phi(\epsilon), \F(\epsilon))$, $\epsilon \in \mathcal E$, to (\ref{eq1})--(\ref{init_cond}) that satisfies
\begin{align}
\label{supp_cond_F}
\p_\epsilon \F^j(\epsilon)|_{\epsilon = 0} = f^j,\ j=1,2,\quad 
\quad \bigcup_{j=1,2} \supp(\F^j(\epsilon)) \subset \bigcup_{j=1,2} \supp(f^j).
\end{align}
Here $f^2 \cdot d \hat \phi$ is the inner product on $\R^L$, that is,
$f^2 \cdot d \hat \phi = \sum_{l = 1}^L f_l^2 d \hat \phi_l$.
\end{proposition}

In the context of the initial value problem for (\ref{eq1})--(\ref{eq2}), there are counterexamples to linearization stability, that is, there are 
solutions to the linearization of (\ref{eq1})--(\ref{eq2})
that do not arise as the derivative of a one-parameter family of solutions to the original non-linear problem. 
In fact, linearization stability for the initial value problem for Einstein's equations in vacuum is characterized by the absence of Killing fields \cite{Arms1979}. The generic condition (ND)
can be viewed as an analogue of this property in the present setting. 

We emphasize that by (\ref{supp_cond_F}), the support of the source functions $\mathcal F(\epsilon)$ can be chosen to be small if the same is true for the supports of $f^1$ and $f^2$. This is crucial from the point of view of the inverse problem that we consider, since it allows us to choose sources supported near the world line $\mu$. 
If the sources $\mathcal F(\epsilon)$ are interpreted as being produced by a measurement device, see Appendix C, $\mathcal F(\epsilon)$ having a small support means that the device can be located in a small neighbourhood of $\mu$.

In what follows, we will use only a microlocal version of Proposition \ref{prop_lin_stab_basic}, but the above version might be of independent interest and we include its proof. The proof is based on the following lemma.

\begin{lemma}
\label{lem_adaptive}
Let $U \subset \hat M$ be open and suppose that $\hat \phi$ satisfies the non-degeneracy condition (ND).
Then there are vector bundle homomorphisms
$$
A_\phi^1 \in C^{\kappa}(U;\Hom(T^* U, \R^L)),
\quad
A_\phi^2 \in C^{\kappa}(U;\Hom(\R^L, \R^L)),
$$
%covering the identity and 
that satisfy the following:
\def\aTensor{S}
\begin{itemize}
\item[(i)] $A_\phi^1(v) \cdot d \phi = v$
for all $v \in C^\infty(U; T^* U)$.
\item[(ii)] $A_\phi^2(w) \cdot d \phi = 0$
for all $w \in C^\infty(U; \R^L)$.
\item[(iii)] 
$A_\phi^1(w \cdot d \phi ) + A_\phi^2(w) = w$
for all $w \in C^\infty(U; \R^L)$.
\item[(iv)] The map $\phi \mapsto (A^1_\phi, A^2_\phi)$
is $C^\infty$-smooth in a neighbourhood of $\hat \phi$,
and depends only on $d \phi$.
\end{itemize}
\end{lemma}

We postpone the proof of the lemma in Appendix A.

\begin{proof}[Proof of Proposition \ref{prop_lin_stab_basic}]
Using Lemma \ref{lem_adaptive}, define 
\begin{align}
\label{adaptive_F}
\F^1(\epsilon) = \epsilon f^1, \quad \F^2(\epsilon) = \epsilon A_\phi^1(-\div_g f^1) + \epsilon A_\phi^2(f^2).
\end{align}
By the discussion in Appendix B, the system 
(\ref{einred_full}),
with $F$ being replaced by $I_g \F(\epsilon)$,
has a unique solution $(g(\epsilon),\phi(\epsilon))$ on $\Mmax$ when $\epsilon \in \mathcal E$ and $\mathcal E \subset \R$ is small enough neighbourhood of the origin.
The compatibility condition (\ref{compcond}) holds by Lemma \ref{lem_adaptive}. Indeed, by the properties (i) and (ii), we have
$
\F^2 \cdot d \phi = - \div_g \F^1.
$
%Corollary \ref{cor_reduced} guarantees that $v(\epsilon)$ is a solution of (\ref{eq_Ein_nonhomog}).
Moreover, as $f^1$ and $f^2$ satisfy the linearized compatibility condition (\ref{compcond_lin}), 
it follows from property (iii) of Lemma \ref{lem_adaptive} that
\begin{align*}
\p_\epsilon \F^2(\epsilon)|_{\epsilon = 0} = A_{\hat \phi}^1(f^2 \cdot d \hat \phi) + A_{\hat \phi}^2(f^2)
= f^2.
\end{align*} 
\end{proof}

\subsection{Notations for analysis of singularities}
\label{sec_a_of_sing}

We refer to \cite{H4} for detailed exposition of the concepts that we use for analysis of singularities.
We use the notation 
$\WF(u)$ 
for the wave front set of a distribution $u \in \D'(M)$.
The wave front is a subset of $T^* M \setminus 0$,
the cotangent bundle with the zero section removed,
and its projection on the base space $M$
is called the singular support, $\singsupp(u)$.
Fibres $T_x^* M$ are occasionally identified with the sets $\{(x,\xi) \in T^*M;\ \xi \in T_x^* M\}$, $x \in M$. 
%The singleton $\{0\}$, both for the origin in $T_x^* M$, $x \in M$, and for the zero section in $T^*M$, is denoted by $0$.

For a conic Lagrangian submanifold 
$\Lambda \subset T^*M \setminus 0$ and a vector bundle $E$ over $M$, 
we denote by $I^{p}(\Lambda;E)$ the space of Lagrangian distributions of order $p \in \R$ associated to $\Lambda$,
and taking values in $E$.
We write also 
$I(\Lambda;E) = \bigcup_{p \in \R} I^p(\Lambda;E)$, and sometimes omit writing $E$ if it is clear from the context.
Recall that 
$$
\WF(u) \subset \Lambda, \quad u \in I(\Lambda;E).
$$
For a submanifold $K \subset M$, we define 
the space of conormal distributions on $K$ by
$I(K;E) = I(N^* K \setminus 0;E)$, where $N^* K$ is the conormal bundle of $K$. 
The singleton $\{x\}$, $x \in M$, is considered to be a zero dimensional manifold with $N^* \{x\} = T_x M$.
We follow the convention that conormal bundles contain the zero section, but the wave front sets do not.

The principal symbol of $u \in I(\Lambda;E)$
is denoted by $\sigma[u](x,\xi)$, $(x,\xi) \in \Lambda$, 
see \cite[Th. 25.1.9]{H4} for the definition,
and we omit writing out the half density and Maslov factors. 
Moreover, we denote the lift of the bundle $E$ to $\Lambda$ still by $E$, see Definition \ref{def_bundle_lift} for a more detailed discussion in the particular case that is important for us. 
The spaces of symbols are denoted by $S^p(\Lambda; E)$, and we write also $$
S(\Lambda; E) = \bigcup_{p \in \R} S^p(\Lambda;E).
$$
When $\Lambda \subset T^*(M \times M) \setminus 0$ the distributions in $I(\Lambda; \End(E))$ correspond to kernels of Fourier integral operators acting on sections of $E$. We denote by $\Delta$ the twisted diagonal of $T^*(M \times M) \setminus 0$, that is, 
$$
\Delta = \{(x,\xi,x,-\xi);\ (x,\xi) \in T^* M \setminus 0\}.
$$
Then $I(\Delta; \End(E))$ are the pseudodifferential operators, and in this case we consider the symbol to be defined on $T^* M \setminus 0$.

For two conic Lagrangian submanifolds $\Lambda_0, \Lambda_1 \subset T^* M \setminus 0$ intersecting cleanly, that is,
\begin{equation*}
T_\lambda(\Lambda_0) \cap T_\lambda(\Lambda_1) = T_\lambda(\Lambda_0 \cap \Lambda_1), \quad \lambda \in \Lambda_0 \cap \Lambda_1,
\end{equation*} the space of paired Lagrangian distributions associated to $\Lambda_0, \Lambda_1$ is denoted by $I^{p,l}(\Lambda_0, \Lambda_1;E)$,
$p,l \in \R$.
Microlocally away from the intersection $\Lambda_0 \cap \Lambda_1$, it holds that $u \in I^{p+l}(\Lambda_0 \setminus \Lambda_1;E)$ and $u \in I^p(\Lambda_1 \setminus \Lambda_0;E)$, see \cite{MU1}. 
We write also $I(\Lambda_0, \Lambda_1;E) 
= \bigcup_{p,l \in \R} I^{p,l}(\Lambda_0, \Lambda_1;E)$.
Let us point out that \cite{MU1} discusses only the case where $E$ is a trivial line bundle, but the generalization for arbitrary vector bundles is straightforward.
We refer also to \cite{Guillemin1981} in relation to theory of paired Lagrangian distributions.
The first use of paired Lagrangian distributions in the context of cosmology was \cite{Guillemin1989}.  

For us, the primary example of a paired Lagrangian distribution is the causal inverse of the wave operator. This is discussed in detail in next section.

\subsection{Propagation of singularities for linear wave equations}
\label{sec_prop_sing}

Let $(M,\hat g)$ be a background spacetime. We denote by $\Box = \Box_{\hat g}$ the wave operator associated to $\hat g$, and by $\Sigma$ and $L_x M \subset T_x^* M$, $x \in M$, the characteristic variety of $\Box$ and light cones of $(M,\hat g)$.
We follow the convention that the characteristic variety does not contain the zero section but each light cone contain the corresponding vertex, that is, 
\begin{align}\label{def_Sigma}
\Sigma = \bigcup_{x \in M} L_x M \setminus 0.
\end{align}
For a set $S \subset T^* (M \times M)$ we define its twist by 
$$
S' = \{(x,\xi,y,-\eta);\ (x,\xi,y,\eta) \in S\}.
$$
Observe that $(S')' = S$.

Let $V$ be a first order differential operator acting on sections of a vector bundle $E$ over $M$. 
As $(M, \hat g)$ is globally hyperbolic, the causal inverse $$Q = (\Box+V)^{-1}$$
is well-defined and its 
Schwartz kernel is in $I(\Delta, \Lambda; \End(E))$
where $\Lambda$ is the twist of the canonical relation $\Lambda'$ of $\Box$ defined by
\begin{align}
\label{def_can_rel}
\Lambda' = \{ (x,\xi,y,\eta) \in \Sigma^2;\ 
(y,\eta) = \beta(s; x,\xi),\ s \in \R \}.
\end{align}
Here the scalar wave operator $\Box$ acts diagonally and 
$\beta(\cdot; x,\xi)$
is the bicharacteristic of $\Box$ 
satisfying $\beta(0; x,\xi) = (x,\xi)$. That is,
$$
\beta(s; x,\xi) = (\gamma_{x,\xi}(s), \dot\gamma_{x,\xi}^*(s)), 
$$
where $\gamma_{x,\xi}$ is the geodesic of $(M,\hat g)$
with the initial data $(x,\xi_*)$, and we use the notation
$\xi_*^j = \hat g^{jk} \xi_k$, $\xi \in T_x^* M$, and $v^*_j = \hat g_{jk} v^k$, $v \in T_x M$, for isomorphism between the tangent and cotangent vectors. 

Recall that $Q$ propagates singularities along $\Lambda' \setminus \Delta$ in the following sense:
if $(x,\xi,y,\eta) \in \Lambda' \setminus \Delta$ then there is linear isomorphism $\tilde q(x,\xi,y,\eta)$ between the fibres $E_x$ and $E_y$ of $E$ such that 
\begin{align}\label{tildeq}
\sigma[Qf](y,\eta) = \tilde q(x,\xi,y,\eta)\sigma[Qf](x,\xi),
\end{align}
whenever $\beta(\cdot; x,\xi)$ does not intersect $\WF(f)$ between $(y,\eta)$ and $(x,\xi)$. 
To be more precise, 
suppose that $(y,\eta) = \beta(s; x, \xi)$ with $s > 0$.
Then (\ref{tildeq}) holds if $\beta(s'; x,\xi) \notin \WF(f)$
for $s' \in [0,s]$.
The case $s < 0$ is analogous. 

The linear map $\tilde q(x,\xi,y,\eta)$ is obtained by solving the transport equation
\begin{align}\label{eq_transport}
-i\mathscr{L}_H \tilde q + c \tilde q = 0
\end{align}
along the bicharacteristic $\beta(\cdot; x,\xi)$, with 
the initial condition given by the identity map, see \cite[Th. 5.3.1]{FIO2}.
Here $H$ is the Hamilton vector field corresponding to the principal symbol $-\hat g(\xi,\xi)$ of $\Box$, $c$ is the subprincipal symbol of $\Box+V$, and $\mathscr{L}_H$ is the Lie derivative along $H$. Solving the transport equation in the opposite direction gives the inverse of $\tilde q(x,\xi,y,\eta)$.

We define for a set $S \subset T^* M$,
\begin{align}\label{def_flowout_S}
\Lambda(S) = \{(y,\eta) \in \Sigma;\ 
(x,\xi,y,\eta) \in \Lambda',\ (x,\xi) \in S \},
\end{align}
and say that $S$
$\Lambda$-invariant if 
\begin{align}\label{def_L_invariant}
S \subset \Sigma
\quad \text{and} \quad
\Lambda(S) = S.
\end{align}
Moreover, we say that a smooth submanifold $K \subset M$ is 
$\Lambda$-invariant if 
$N^* K \setminus 0$ is.

Let $Y$ be a submanifold of $M$.
We assume that each lightlike geodesic intersects $Y$ at most once, and that the intersection is transversal.
If $f \in I(Y;E)$, then $u = Qf$ is in $I(N^* Y \setminus 0, \Lambda(N^* Y);E)$.
A more precise statement can be formulated 
by using the conic Lagrangian submanifold with boundary,
the flowout $N^* Y$ in the future direction,
$$
\Lambda^+(N^* Y) = \{(y,-\eta) \in \Sigma;\ 
(y,\eta) = \beta(s; x,\xi),\ s \in \R,\ (x,\xi) \in \Sigma \cap N^* Y,\ x \le y \}.
$$
Taking into account the causality, it holds that 
$u \in I(N^* Y \setminus 0, \Lambda^+(N^* Y);E)$, see \cite{MU1}.
There the authors give also 
a relation between $\sigma[u]$, as a symbol in $S(\Lambda^+(N^* Y);E)$, and $\sigma[f]$, as a symbol in $S(Y;E)$,
at points $(x,\xi) \in \p \Lambda^+(N^* Y)$.
The equation (6.7) in \cite{MU1} defines a map
\begin{align}\label{transport_init}
\sigma[f](x,\xi) \mapsto \sigma[u](x,\xi), \quad (x,\xi) \in \p \Lambda^+(N^* Y),
\end{align}
that is a linear automorphism on the fibre $E_x$.
We denote the composition of $\tilde q$ and (\ref{transport_init}) by $q$. Then $q = q(x,\xi,y,\eta)$ gives a linear isomorphism
\begin{align}\label{def_prop_q}
q : E_x \to E_y, \quad q : \sigma[f](x,\xi) \mapsto \sigma[u](y,\eta), \quad (y,\eta) = \beta(s;x,\xi),\ s \in \R.
\end{align}

\subsection{Microlocal linearization stability}
\label{sec_muloc_linstab}
\def\SSym{\mathbb S}

In this section we give a microlocal version of Proposition \ref{prop_lin_stab_basic}. We begin by considering the principal part of the compatibility condition (\ref{compcond_lin}).

\begin{definition}
\label{def_bundle_lift}
For a submanifold $Y \subset M$,
$\Sym_Y$ is the lift of 
$$\Sym = \Symop(T^*M \otimes T^*M)$$ on
$N^* Y \setminus 0$. That is,
$\Sym_Y$ is a vector bundle over $N^* Y \setminus 0$
and at a point $(x,\xi) \in N^* Y \setminus 0$
its fibre is 
$$
(\Sym_Y)_{(x,\xi)} = \Symop(T^*_x M \otimes T^*_x M).
$$
When considering the space of principal symbols of conormal disributions on $Y$, taking values on $\Sym_Y$, we write simply $S(Y;\Sym)$.
\end{definition}

\begin{definition}
\label{def_SY_comp}
Let $(M,\hat g)$ be a background spacetime,
and let $Y \subset M$ be a submanifold.
We define a subbundle $\SSym_Y$ of $\Sym_Y$ over
$N^* Y \setminus 0$
by the following equation for $h \in \Sym_Y$,
\begin{align}
\label{def_S_Y}
h(\xi_*, v) - \frac{\tr_{\hat g} h}{2} \hat g(\xi_*, v) = 0,
\quad \text{for all } (x,\xi) \in N^* Y \setminus 0,\ v \in T_x M,
\end{align}
where
$\xi_*$ is the 
vector $\xi_*^j = \hat g^{jk} \xi_k$.
We denote by $S(Y; \SSym)$
the subspace of symbols $S(Y; \Sym)$ taking values on $\SSym_Y$, that is,
\begin{align}\label{def_SSY}
S(Y; \SSym) = \{ s \in S(Y; \Sym);\ s(x,\xi) \in \SSym_Y,\ (x,\xi) \in N^* Y \setminus 0 \}.
\end{align}
\end{definition}

\begin{lemma}
\label{lem_subbundle_rank}
The bundle $\SSym_Y$ has rank 6.
\end{lemma}
\begin{proof}
Equation (\ref{def_S_Y}) is equivalent with $\tilde h = I_{\hat g} h$ and 
\begin{align}
\label{def_S_Y_tilde}
\tilde h(\xi_*, v) = 0,
\quad \text{for all } (x,\xi) \in N^* Y \setminus 0,\ v \in T_x M.
\end{align}
Write $e_1 = (1, 0, 0, 0)$, $(\tilde h \xi_*)_j = \tilde h_{jk} \xi_*^k$,
and $\hat \xi = \xi_* / |\xi_*|$ where $|\cdot|$ is the Euclidean norm in the coordinates. Choose a 
neighbourhood $W$ of $(x,\xi)$ and 
a smooth map 
$R : W \to O(4)$ such that 
$R(\hat \xi) \hat \xi = e_1$ for $(x,\xi) \in W$.
Define for $a = (a_1, \dots, a_6)$,
\begin{align}\label{param_SSym}
\mathcal H_{ker} : \R^6 \to \Symop(\R^4 \otimes \R^4), \quad 
\mathcal H_{ker}(a) = 
\begin{pmatrix}
0 & 0 & 0 & 0
\\
0 & a_1 & a_2 & a_3
\\
0 & a_2 & a_4 & a_5
\\
0 & a_3 & a_5 & a_6
\end{pmatrix}.
\end{align}
Now 
$
\tilde h \xi_* = |\xi_*| \tilde h R(\hat \xi)^T e_1,
$
whence, in local coordinates, (\ref{def_S_Y_tilde}) is equivalent with
\begin{align}\label{def_S_Y_local}
R(\hat \xi) \tilde h R(\hat \xi)^T \in \mathcal H_{ker}(\R^6),
\quad \text{for all } (x,\xi) \in N^* Y \setminus 0.
\end{align}
Now for each $(x,\xi) \in W$, the map $a \mapsto I_{\hat g} R(\hat \xi)^T \mathcal H_{ker}(a) R(\hat \xi)$ gives a parametrization of the fibre $(\SSym_Y)_{(x,\xi)}$, and thus $\SSym_Y$  has rank $6$.
\end{proof}

\begin{proposition}
\label{prop_muloc_linstab}
Let $(M,\hat g)$ be a background spacetime, 
let $Y \subset M$ be a submanifold, 
let $U \subset M$ be a neighbourhood of $\overline Y$,
and suppose that 
$\hat \phi \in C^{\infty}(M; \R^L)$ satisfies the non-degeneracy condition (ND).
Write $E = \Sym \oplus \R^L$ and 
\beq\label{notations numbered 2}
I_Y^p = \{f \in  I^p(Y; E);\ \sigma[f] \in S^{p+1}(Y; \SSym \oplus \R^L) \}.
\eeq
Then there are 
$p_0 \in \R$  
and $k_0 \in \N$
such that for $p \le p_0$ and $k \ge k_0$
there are maps
\begin{align*}%\label{}
&A_{(g,\phi)} : I_Y^p \to I^p(Y; E),
\quad (g,\phi) \in C^k(M; E),
\\
&\mathbb F : S^{p+1}(Y; \SSym \oplus \R^L) \to I_Y^p,
\end{align*}
with the following properties for each $s \in S^{p+1}(Y; \SSym \oplus \R^L)$, 
\begin{itemize}
\item[(A1)] 
$f = \mathbb F(s)$ satisfies the support condition $\supp(f) \subset U$, $A_{(g,\phi)}$ is a differential operator, and
$\mathcal F = A_{(g,\phi)}(f)$ satisfies the compatibility condition (\ref{compcond}).
\item[(A2)] 
$F = I_{\hat g} A_{(\hat g, \hat \phi)}(f)$
satisfies $\sigma[F] = s$.
\item[(A3)] The map $(g,\phi) \mapsto A_{(g,\phi)}(f)$
is $C^\infty$-smooth in a neighbourhood of $(\hat g, \hat \phi)$, and depends on the derivatives of $g$ and $\phi$ only up to the first order.
\end{itemize} 
\end{proposition}
\begin{proof}
Observe that for each $(x,\xi) \in N^* Y \setminus 0$ the map 
\begin{align*}
%\label{subbundle_surj}
\iota_\xi h = I_{\hat g} h(\xi_*, \cdot)
\end{align*}
is surjective from $(\Sym_Y)_{(x,\xi)}$ to $T^*_x M$.
Indeed, the kernel of $\iota_\xi$ is given by the 
equation (\ref{def_S_Y_tilde}), and as seen in the proof of Lemma \ref{lem_subbundle_rank}, the kernel is of dimension 6.
As the rank of $\Sym$ is 10, the map $\iota_\xi$ has 4 dimensional range.
Let us now choose a right inverse $r_\xi$ of $\iota_\xi$
such that $r_\xi \in C^\infty(N^* Y \setminus 0; \Hom(T^*M,\Sym))$ and that $r_\xi$ maps $S^{p+1}(Y; T^* M)$ to $S^{p}(Y; \Sym)$.

Let $s = (s^1, s^2) \in S^{p+1}(Y; \SSym \oplus \R^L)$.
As the principal symbol map is an isomorphism 
$$
I^p(Y; E) / I^{p-1}(Y; E) \to S^{p+1}(Y; E) / S^{p}(Y; E),
$$
see e.g. \cite[Th. 25.1.9]{H4}, we can choose 
\begin{align*}%\label{}
&\mathbb F^1_{pri} : S^{p+1}(Y; \SSym) \to I^p(Y; \Sym),
\quad
\mathbb F^2 : S^{p+1}(Y; \R^L) \to I^p(Y; \R^L),
\\&
\mathbb F^1_{sub} : S^{p}(Y; \Sym) \to I^{p-1}(Y; \Sym),
\end{align*}
so that 
$\sigma[\mathbb F^1_{pri}(s^1)] = s^1$, 
$\sigma[\mathbb F^2(s^2)] = s^2$, 
$\sigma[\mathbb F^1_{sub}(s^3)] = s^3$,
and that $\mathbb F^1_{pri}(s^1)$,
$\mathbb F^2(s^2)$
and 
$\mathbb F^1_{sub}(s^3)$ are all supported in $U$.
As $I_{\hat g}$ is an involution, by (\ref{def_S_Y})
it holds for $s^1 \in S^{p+1}(Y; \SSym)$ that the principal symbol of 
$\div_{\hat g} I_{\hat g} \mathbb F^1_{pri}(s^1)$ vanishes. In other words,
$$
\sigma[\div_{\hat g} I_{\hat g} \mathbb F^1_{pri}(s^1)]
= \iota_\xi \sigma[\mathbb F^1_{pri}(s^1)]
\in S^{p+1}(Y; T^* M). 
$$ 
Thus we can define the map 
\begin{align*}%\label{}
&z : S^{p+1}(Y; \SSym) \times S^{p+1}(Y; \R^L)
\to S^{p}(Y; \Sym),
\\&
z(s^1, s^2) = r_\xi \left( - \iota_\xi \sigma[\mathbb F^1_{pri}(s^1)] -s^2 \cdot d\hat \phi \right),
\end{align*}
and also the maps
\begin{align}\label{def_E_and_A}
\mathbb F(s^1, s^2) &= \left(\mathbb F^1_{pri}(s^1) + \mathbb F^1_{sub}(z(s^1,s^2)),\, \mathbb F^2(s^2) \right),
\\\notag
A_{(g,\phi)}(f^1, f^2) &= (I_g f^1, A_\phi^1(-\div_g I_g f^1) + A_\phi^2(f^2)).
\end{align}
Observe that the definition of $A_{(g,\phi)}$ is closely related to the choice in (\ref{adaptive_F}).

Now the the compatibility condition (\ref{compcond}) follows from (i) and (ii) in Lemma \ref{lem_adaptive}.
Setting $(f^1, f^2) = \mathbb F(s^1, s^2)$
and $(F^1, F^2) = I_{\hat g} A_{(\hat g, \hat \phi)}(f^1, f^2)$, we have 
$$
\sigma[F^1] = \sigma[I_{\hat g}^2 f^1] = \sigma[f^1] = \sigma[\mathbb F^1_{pri}(s^1)] = s^1.
$$
Moreover,  
\begin{align*}%\label{}
&\sigma[\div_{\hat g} I_{\hat g}  f^1]) 
= \iota_\xi \sigma[\mathbb F^1_{pri}(s^1) + \mathbb F^1_{sub}(z(s^1,s^2))] 
\\&\quad= 
\iota_\xi \sigma[\mathbb F^1_{pri}(s^1)] 
+ \iota_\xi r_\xi\left(- \iota_\xi \sigma[\mathbb F^1_{pri}(s^1)] -s^2 \cdot d\hat \phi\right)  
= -s^2 \cdot d\hat \phi,
\end{align*} 
and therefore using (iii) in Lemma \ref{lem_adaptive}
it holds that
\begin{align*}%\label{}
\sigma[F^2] &= A_{\hat \phi}^1(-\sigma[\div_{\hat g} I_{\hat g}  f^1]) + A_{\hat \phi}^2(\sigma[f^2])
= A_{\hat \phi}^1(s^2 \cdot d \hat \phi) + A_{\hat \phi}^2(s^2) = s^2.
\end{align*}
The smoothness follows from (iv) in Lemma \ref{lem_adaptive}.
\end{proof}

Consider $s \in S^p(Y; \SSym \oplus \R^L)$
with negative enough $p \in \R$ and define the family $F(\epsilon) = \epsilon I_g A_{(g,\phi)}(E(s))$, $\epsilon \in \R$.
Now the system (\ref{einred_full})
has a unique solution $(g,\phi)$ on $\Mmax$ when $\epsilon \in \mathcal E$ and $\mathcal E \subset \R$ is small enough neighbourhood of the origin.
Moreover, setting $\mathcal F = I_g F$,
we have that the family $(g,\phi,\mathcal F)$, $\epsilon \in \mathcal E$, satisfies (\ref{eq1})--(\ref{eq2}).
As $\sigma[\p_\epsilon F^j|_{\epsilon = 0}] = s$,
we say that the equation (\ref{eq1})--(\ref{eq2})
has microlocal linearization stability 
for principal symbols in $S^p(Y; \SSym \oplus \R^L)$.

\subsection{Propagation of singularities for linearized Einstein's equations}
\label{sec_prop_sing_ein}

Let $\mathcal E \subset \R$ be a neighbourhood of the origin, and let us consider a family 
\begin{align}\label{lin_family_1d}
(g(\epsilon), \phi(\epsilon), F(\epsilon))
\in C^1(\mathcal E; C^2(\Mmax;(\Sym \oplus \R^L)^2))
\end{align}
of solutions to (\ref{einred_full}) near a background spacetime and scalar fields $(M,\hat g)$ and $\hat \phi$ in the sense that $g(0) = \hat g$, $\phi(0) = \hat \phi$
and $F(0) = 0$. 
We write 
%$\dot u = \p_\epsilon u |_{\epsilon = 0}$, $u = g,\phi,F,\dots$, and 
%\quad f = (f^1, f^2) = \dot F.
\begin{align}\label{def_v}
v(\epsilon) = (g(\epsilon) - \hat g, \phi(\epsilon) - \hat \phi).
\end{align}
Then $\dot v = \p_{\epsilon} v|_{\epsilon=0}$ satisfies a linear differential equation that we denote by 
\begin{align}\label{ein12_linearized}
\Box_{\hat g, \hat \phi} \dot v = \dot F \quad \text{on $\hat M$},
\end{align}
where $\dot F = \p_{\epsilon} F|_{\epsilon=0}$.
We stress that (\ref{ein12_linearized}) is coordinate invariant since (\ref{einred_full}) is, and express next (\ref{ein12_linearized}) in local coordinates.

To simplify the notation, we rescale the first equation in (\ref{einred_full}) by $-2$. Then (\ref{einred_full}) can be written locally in the form 
\begin{align}
\label{ein_abs}
g^{-1}(D,D) v + A(x, v, Dv) = F,
\end{align}
where, for a local trivialization $U \subset \R^4 \times \R^{10+L}$ of the bundle $\Sym \oplus \R^L$, 
$$
(x, y ,z) \mapsto A(x, y ,z) : U \times \R^{4(10+L)} \to \R^{10+L}
$$ 
is a smooth map.
Differentiating (\ref{ein_abs}) with respect to $\epsilon$ at $\epsilon = 0$, we obtain
\begin{align}
\label{ein12_lin_loc}
\hat g^{-1}(D,D) \dot v + D_z A(x,0,0) D \dot v + D_y A(x,0,0) \dot v = \dot F.
\end{align}
Thus $\Box_{\hat g, \hat \phi}$ is of the form $\Box + V$
where $\Box = \Box_{\hat g}$ and $V$ is a first order differential operator on the vector bundle $\Sym \oplus \R^L$.
%Note that the sum in (\ref{ein12_lin_loc})
%identifies only the principal symbol $\hat g^{-1}(\xi,\xi)$ 
%of $\Box_{\hat g, \hat \phi}$; it is not coordinate invariant in the sense that 
%the local representation of $V \dot v$ is not given by
%the two latter terms on the left-hand side of (\ref{ein12_lin_loc}). 

As the linearized equation (\ref{ein12_linearized}) is of the form discussed in Section \ref{sec_prop_sing},
we know how it propagates singularities.
However, we still need to understand how the 
constraint (\ref{def_S_Y}) is propagated.
Put differently, we analyse next what is the image of the rank 6 subbundle $\SSym_Y$ of $\Sym_Y$
under the isomorphism $q$ defined by (\ref{def_prop_q}).

Assume now that the family (\ref{lin_family_1d})
is constructed as in the end of Section \ref{sec_muloc_linstab}.
Then $\sigma[\dot F^1]$
takes values in $\SSym_Y$, and $g$
satisfies the gauge condition (\ref{gauge_cond}).
Differentiating (\ref{gauge_cond}) with respect to $\epsilon$
at $\epsilon = 0$, and writing the result in local coordinates, gives
$$
\hat g^{pq} ( \p_p \dot g_{kq} + \p_q \dot g_{pk} - \p_k \dot g_{pq} )
= r_k \dot g, \quad k=0,1,2,3,
$$
where $\dot g = \p_{\epsilon} g|_{\epsilon=0}$ and $r_k$ is a zeroth order operator. 
Taking the principal symbol of this expression yields
\begin{align*}
\hat g^{pq} (\xi_p \sigma[\dot g_{kq}] + \xi_q \sigma[\dot g_{pk}] - \xi_k \sigma[\dot g_{pq}]) = 0, \quad k=0,1,2,3,
\end{align*} 
or equivalently, writing $h = \sigma[\dot g]$,
\begin{align}\label{harmonicity_cond}
\xi_*^p h_{pk} - \frac{\tr_{\hat g} h} 2 \xi_k = 0, \quad k=0,1,2,3.
\end{align}
Supposing that $N^* K \setminus 0 \subset \Lambda^+(N^* Y)$
for a submanifold $K \subset M$,
these four conditions for $h$ 
coincide with the definition of 
the rank 6 (equivalently, codimension 4) subbundle $\SSym_K$ on $N^* K$, see Definition \ref{def_SY_comp} with $Y$ being replaced by $K$.
As $q$ is a linear isomorphism,
we see that the map 
\begin{align}\label{prop_subbundles}
\sigma[\dot F](x,\xi) \mapsto \sigma[\dot v](y,\eta), \quad (y,\eta) = \beta(s;x,\xi),\ s \in \R,
\end{align}
gives an isomorphism 
between the fibres $(\SSym_Y \oplus \R^L)_{(x,\xi)}$
and $(\SSym_K \oplus \R^L)_{(y,\eta)}$.
Here $(x,\xi)$ is in $\p \Lambda^+(N^* Y)$ and $\beta(\cdot;x,\xi)$ is assumed to intersect $\p \Lambda^+(N^* Y)$ only at $(x,\xi)$.

Recall that the gauge condition (\ref{gauge_cond})
implies that (\ref{Box_phi_subprin}) holds. 
This again implies that, writing $\dot v = (\dot g, \dot \phi)$, the second term in (\ref{ein12_lin_loc}) is of the form
$$
D_z A(x,0,0) D \dot v
= \begin{pmatrix}
a_{11}(x) & a_{12}(x)
\\ 
0 & a_{22}(x)
\end{pmatrix}
\begin{pmatrix}
D \dot g
\\
D \dot \phi
\end{pmatrix}
,
\quad a_{22}(x) D \dot \phi
= -\hat g^{pq} \hat \Gamma^j_{pq} \p_{x^j} \dot \phi.
$$
In particular, the subprincipal symbol of $\Box_{\hat g, \hat \phi}$, which appears as $c$ in (\ref{eq_transport}),
is of the above upper triangular form. 
This observation leads to the following remark.

\begin{remark}
\label{rem_restriction_S}
The restriction of the map (\ref{prop_subbundles}) on the subspace 
$(\SSym_Y)_{(x,\xi)}$ of the fibre $(\SSym_Y \oplus \R^L)_{(x,\xi)}$
is an isomorphism between the two subspaces of dimension $6$,
\begin{align}\label{isomorp_rest}
(\SSym_Y)_{(x,\xi)} \to (\SSym_K)_{(y,\eta)}.
\end{align}
\end{remark}

\section{Interaction of singular solutions}
\label{sec_interaction}

In this section we show how a non-linear interaction of singular solutions to (\ref{ein12_linearized}) can behave, on the microlocal level, as a point source. In order to have a product calculus for principal symbols, we consider the case where the singular solutions are conormal distributions.

As the singularities of a solution to (\ref{ein12_linearized}) are invariant under the map $\Lambda$, see (\ref{def_flowout_S}), we are restricted to the case where the conormal distributions are associated to submanifolds of codimension one in $M$. 
This again implies that we need four such submanifolds to get a single point as their intersection. For this reason, we begin by studying the linearization of (\ref{einred_full}) with respect to a four dimensional family of solutions.
 {\nntext We note that the interaction of two waves having very strong, impulse type singularities
has been studied in  \cite{L-Rodnianski}.} 

\subsection{Fourfold linearization}

Let $\epsilon \in \mathcal E$, a neighbourhood of the origin in $\R^4$, and consider a family 
\begin{align}\label{lin_family}
(g(\epsilon), \phi(\epsilon), F(\epsilon))
\in C^4(\mathcal E; C^2(\Mmax;(\Sym \oplus \R^L)^2))
\end{align}
of solutions to (\ref{einred_full}) near a background spacetime and scalar fields $(M,\hat g)$ and $\hat \phi$. Define $v$ as in (\ref{def_v}).
Analogously to (\ref{ein12_linearized}), 
we see that 
\begin{align}
\label{def_vj_deriv}
v_j = \p_{\epsilon_j} v|_{\epsilon = 0}, \quad j=1,2,3,4,
\end{align}
satisfies the wave equation
\begin{align}\label{eq_vj}
\Box_{\hat g, \hat \phi} v_j = f_j,
\end{align}
where $f_j = \p_{\epsilon_j} F|_{\epsilon = 0}$.
Let us now consider interactions of linearized solutions $v_j$.
To simplify the computations we assume that 
\begin{align}\label{no_source_interactions}
\p_{\epsilon_j} \p_{\epsilon_k} F|_{\epsilon = 0}
=
\p_{\epsilon_j} \p_{\epsilon_k} \p_{\epsilon_l} F|_{\epsilon = 0}
=
\p_{\epsilon_1} \p_{\epsilon_2} \p_{\epsilon_3} \p_{\epsilon_4} F|_{\epsilon = 0} = 0
\end{align}
for all distinct $j$, $k$ and $l$.

Recall that $A(x,y,z)$ is defined by (\ref{ein_abs}).
We use the shorthand notation
$$
A^k_j(x,y,z) = D_y^j D_z^k A(x,y,z)
$$
for the derivatives of $A$, that we consider as multilinear forms, and omit writing $(x,y,z)$ as a parameter for $A$ below.
Note that $A$ is a quadratic polynomial in $z$,
and therefore $A_j^k = 0$ when $k > 2$.
We write also
\def\P{\mathcal P}
\begin{align}\label{def_vderivs}
v_{j_1\dots j_n} = \p_{\epsilon_{j_1}} \dots \p_{\epsilon_{j_n}} v|_{\epsilon = 0},
\end{align}
denote by $S_n$ the permutations on $\{1,2,\dots,n\}$, and define
$$
P(u_1,u_2,\dots,u_n) 
= \hat g^{-1} g_1 \hat g^{-1} g_2
\dots \hat g^{-1} g_n \hat g^{-1}(D,D),
$$
where $u_j = (g_j, \phi_j) \in C^2(\Mmax;\Sym \oplus \R^L)$.

Differentiating (\ref{ein_abs}) with respect to $\epsilon_j$ and $\epsilon_k$, $j \ne k$, 
we see that $v_{jk}$ satisfies the wave equation 
\def\I{\mathscr I}
\def\RR{\mathscr R}
\def\FF{\mathscr F}
\def\Q{\mathscr Q}
\begin{align}
\label{eq_vjk}
\Box_{\hat g, \hat \phi} v_{jk} = 
\I_2(v_j, v_k) + \RR_2(v_j, v_k),
\end{align}
where
\begin{align*}
%\label{def_hjk}
&\I_2(u_1, u_2) = \sum_{\tau \in S_2} \tilde \I_2(u_{\tau(1)}, u_{\tau(2)}), \quad \tilde \I_2(u_1, u_2) = P(u_1) u_2 - \frac 1 2 A^2(D u_1, D u_2),
\\
&\RR_2(u_1, u_2) = \sum_{\tau \in S_2} \tilde \RR_2(u_{\tau(1)}, u_{\tau(2)}), \quad \tilde \RR_2(u_1, u_2) =
- A^1_1(Du_1, u_2)
- \frac 1 2 A_2(u_1, u_2). 
\end{align*}
Observe that $\I_2$ and $\RR_2$ are bilinear differential operators of order $2$ and $1$, respectively.
Here the order is the total number of derivatives. 

Differentiating (\ref{ein_abs}) with respect to $\epsilon_j$, $\epsilon_k$ and $\epsilon_l$, with $j$, $k$ and $l$ distinct, 
we see that $v_{jkl}$ satisfies the wave equation
\begin{align}
\label{eq_vjkl_pre}
\Box_{\hat g, \hat \phi} v_{jkl} = 
\sum_{\tau \in S(j,k,l)} 
& \Bigl(
\I_3^1(v_{\tau(1)}, v_{\tau(2)}, v_{\tau(3)})
+\I_3^2(v_{\tau(1)\tau(2)}, v_{\tau(3)})
\\\notag&
+\RR_3^1(v_{\tau(1)}, v_{\tau(2)}, v_{\tau(3)})
+\RR_3^2(v_{\tau(1)\tau(2)}, v_{\tau(3)})
\Bigr),
\end{align}
where $
S(j,k,l) = \{\tau : \{1,2,3\} \to \{j,k,l\};\ \text{$\tau$ bijection}\},
$
and $\I_3^1$ and $\I_3^2$
are the following trilinear and bilinear second order differential operators
\begin{align*}
\I_3^1(u_1, u_2, u_3)
&= 
-P(u_1, u_2) u_3
-\frac 1 2 A^2_1(D u_1,D u_2, u_3),
\\
\I_3^2(u_1, u_2)
&= \frac 1 2 \sum_{\tau \in S_2} \left(
P(u_1) u_2
- \frac 1 2 A^2(D u_1, D u_2) \right).
\end{align*}
Here $\RR_3^1$ and $\RR_3^2$
are trilinear and bilinear first order differential operators
of similar structure.
%expand $v_{\tau(1)\tau(2)}$ using (\ref{eq_vjk}),
%$$
%v_{\tau(1)\tau(2)} = 
%\Box_{\hat g, \hat \phi}^{-1} \I_2(v_{\tau(1)}, v_{\tau(2)})
%+ \Box_{\hat g, \hat \phi}^{-1} \RR_2(v_{\tau(1)}, v_{\tau(2)}),
%$$
We define 
\begin{align*}%\label{}
\I_3(u_1, u_2, u_3) &= 
\sum_{\tau \in S_3}
\Bigl( \I_3^1(u_{\tau(1)}, u_{\tau(2)}, u_{\tau(3)})
+ \I_3^2(\Box_{\hat g, \hat \phi}^{-1} \I_2(u_{\tau(1)}, u_{\tau(2)}), u_{\tau(3)}) \Bigr),
%\\
%\RR_3(v_j, v_k, v_l)
%&= \sideset{}{'}\sum_{\tau \in S(j,k,l)}
%\Bigl( 
%\RR_3^1(v_{\tau(1)}, v_{\tau(2)}, v_{\tau(3)})
%+ \I_3^2(\Box_{\hat g, \hat \phi}^{-1} \RR_2(v_{\tau(1)},v_{\tau(2)}), v_{\tau(3)})
%\\&\qquad\qquad
%+\RR_3^2(\Box_{\hat g, \hat \phi}^{-1}(\I_2(v_{\tau(1)}, v_{\tau(2)}) + \RR_2(v_{\tau(1)}, v_{\tau(2)})), v_{\tau(3)})
%\Bigr).
\end{align*}
and use an analogous definition for $\RR_3$.
Then (\ref{eq_vjkl_pre}) can be rewritten in the form
\begin{align}
\label{eq_vjkl}
\Box_{\hat g, \hat \phi} v_{jkl} 
&= 
\I_3(v_j, v_k, v_l) + \RR_3(v_j, v_k, v_l).
\end{align}
Recall that $\Sigma$ is the characteristic variety of the wave operator $\Box_{\hat g, \hat \phi}$.
In the symbol computations in Sections \ref{sec_abs_interaction} and \ref{sec_prin4} below, we will consider the causal inverse $\Box_{\hat g, \hat \phi}^{-1}$ on the 
microlocal region $T^* M \setminus \Sigma$ where it is a pseudodifferential operator of order $-2$.
In this region $\I_3$ and $\RR_3$ are trilinear pseudodifferential operators of orders $2$ and $1$, respectively. 

Differentiating (\ref{ein_abs}) again, we see that 
$v_{1234}$ satisfies the wave equation 
\begin{align}
\label{eq_v1234}
\Box_{\hat g, \hat \phi} v_{1234} = 
\I_4(v_1, v_2, v_3, v_4) + \RR_4(v_1, v_2, v_3, v_4),
\end{align}
where
\begin{align}\label{fourfold_expansion}
&\I_4(v_1, v_2, v_3, v_4) 
= \sum_{\tau \in S_4} \tilde \I_4 (v_{\tau(1)}, v_{\tau(2)}, v_{\tau(3)}, u_4),
\end{align}
and writing $\Q_j = \Box_{\hat g, \hat \phi}^{-1} \I_j$, $j=2,3$, 
\begin{align*}
&\tilde \I_4 (u_1, \dots, u_4) = 
P(u_1, u_2, u_3) u_4
-\frac 1 2 P(u_1, 
\Q_2(u_2, u_3)
) u_4
-\frac 1 2 P(
\Q_2(u_1, u_2)
, u_3) u_4
\\\notag&\quad
+\frac 1 4 P(
\Q_2(u_1, u_2) 
) 
\Q_2(u_3, u_4) 
- \frac 1 2 P(u_1, u_2) 
\Q_2(u_3, u_4)
+\frac 1 6 P(
\Q_3(u_1, u_2, u_3)
) u_4
\\\notag&\quad
+\frac 1 6 P(u_1) 
\Q_3(u_2, u_3, u_4)
-\frac 1 4 A^2(D 
\Q_2(u_1, u_2)
, D 
\Q_2(u_3, u_4)
)
\\\notag&\quad
-\frac 1 6 A^2(D 
\Q_3(u_1, u_2, u_3)
, D u_4)
%\\\notag&\quad
-\frac 1 2 A^2_1(D 
\Q_2(u_1, u_2)
, D u_3, u_4)
\\\notag&\quad
-\frac 1 4 A^2_1(D u_1, D u_2, 
\Q_2(u_3, u_4))
%\\\notag&\quad
-\frac 1 4 A^2_2(Du_1,Du_2,u_3,u_4).
\end{align*}
Also, in the microlocal region $T^* M \setminus \Sigma$, 
$\RR_4$ is a 4-linear pseudodifferential operator of order $1$
of similar structure. 

We will next consider the case that $v_j$, ${j=1,2,3,4}$, are conormal distributions associated to a $\Lambda$-invariant submanifolds of codimension one. 
We derive first an upper bound, in the sense of the inclusion relation, for the outgoing singular support,
$\singsupp(v_{1234})$, assuming that not all the four incoming singular supports, $\singsupp(v_j)$, $j=1,2,3,4$, intersect, see Lemma \ref{lem_F4_WF_bound}. 
Then we proceed to study the principal symbol of the
fourfold interaction terms $\I_4 + \RR_4$. 
We consider the terms where $\Q_2$ appears at most once in Lemma \ref{lem_type1} and the terms where  $\Q_2$ appears twice in Lemma \ref{lem_type2}.

\subsection{Upper bounds for wave front sets and singular supports}

We omit the proof of the following elementary lemma.

\begin{lemma}
\label{lem_transversality}
Let $M$ be a smooth manifold and let $K_1, K_2 \subset M$ be transversal submanifolds. Then for all $x \in K_1 \cap K_2$ it holds that 
$
N_x^* (K_1 \cap K_2) = N_x^* K_1 \oplus N_x^* K_2.
$
%In particular, 
%\begin{align}
%\label{prod_okay}
%(x, \xi) \in N^* K_1 \setminus 0 \quad \text{implies} \quad %(x,-\xi) \notin N^* K_2.
%\end{align}
\end{lemma}
%\begin{proof}
%\HOX{The proof is elementary. Remove from the final version.}
%Let $x \in K_1 \cap K_2$.
%By transversality 
%$$
%T_x M = T_x K_1 + T_x K_2, \quad
%T_x (K_1 \cap K_2) = T_x K_1 \cap T_x K_2.
%$$
%In particular, for $j=1,2$,
%\begin{align*}
%N_x^* K_j &= \{\xi \in T_x^* M;\ (\xi, v) = 0\ \text{for $v \in T_x K_j$}\}
%\\&\subset 
%\{\xi \in T_x^* M;\ (\xi, v) = 0\ \text{for $v \in T_x K_1 \cap T_x K_2$}\}
%= N_x^*(K_1 \cap K_2).
%\end{align*}
%Let $v \in T_x M$ and let $v_j \in T_x K_j$ satisfy $v = v_1 + v_2$.
%Then for $\xi \in N_x^* K_1 \cap N_x^* K_2$ it holds that 
%$$
%(\xi,v) = (\xi, v_1) + (\xi, v_2) = 0,
%$$
%and therefore the sum $N_x^* K_1 + N_x^* K_2 \subset N_x^*(K_1 \cap K_2)$
%is direct.

%For any submanifold $K \subset M$ it holds that 
%the dimension of the fibre 
%$N_x^* K$ as a vector space coincides with the codimension of $K$ as a submanifold. By transversality
%\begin{align*}
%\dim(N_x^* K_1 \oplus N_x^* K_2)
%&= \codim(K_1) + \codim(K_1) = \codim(K_1 \cap K_2) 
%\\&= \dim(N_x^*(K_1 \cap K_2)),
%\end{align*}
%and we have shown that 
%$
%N_x^* (K_1 \cap K_2) = N_x^* K_1 \oplus N_x^* K_2.
%$

%For (\ref{prod_okay}), observe that $\xi \in N_x^* K_1$
%and $-\xi \in N_x^* K_2$ imply that 
%$$
%0 = \xi -\xi \in N_x^* K_1 \oplus N_x^* K_2.
%$$
%As the sum is direct, it holds that $\xi=0$.
%\end{proof}

%Below we allow $u$ to take values on a vector bundle over $M$, but suppress this in the notation.
Recall that $\Lambda$-invariance is defined by (\ref{def_L_invariant}).

\begin{lemma}[Twofold interaction]
\label{lem_2_waves}
Let $u_1, u_2 \in \mathcal D'(M)$
and suppose that $\WF(u_j)$ is $\Lambda$-invariant
for $j=1,2$.
Let $Q \in I(\Delta, \Lambda)$. The following holds
\begin{itemize}
\item[(i)] If $(x, \xi) \in \WF(u_1)$ implies $(x,-\xi) \notin \WF(u_2)$, then
$$
\singsupp(Q(u_1 u_2)) \subset \singsupp(u_1) \cup \singsupp(u_2).
$$
\item[(ii)]
If $u_j \in I(K_j)$, $j=1,2$,
for transversal submanifolds $K_1$ and $K_2$ in $M$, then
$$
\WF(Q(u_1 u_2)) \subset N^* K_1 \cup N^* K_2 \cup N^*(K_1 \cap K_2).
$$
\end{itemize}
\end{lemma} 
\begin{proof}
By the assumption in (i) and \cite[Th. 8.2.10]{H1}, 
\begin{align*}
\WF(u_1 u_2) \subset \{(x, \xi+\eta) \in T^* M \setminus 0;\ 
&\text{$(x, \xi) \in \WF(u_1)$ or $\xi=0$, and}
\\
&\text{$(x, \eta) \in \WF(u_2)$ or $\eta=0$}
\}.
\end{align*}
As the sum of two lightlike vectors is lightlike if and only if it is a scalar multiple of one of the summands, see e.g. \cite[Lem. 27, p. 141]{ONeill1983},
claim (i) follows from \cite[Prop. 2.1]{GU1} 
due to the $\Lambda$-invariance of $\WF(u_j)$, $j=1,2$.

Let us now consider (ii). Due to the transversality,
we can apply again \cite[Th. 8.2.10]{H1}, 
and analogously with the above, 
\begin{align}
\label{no_new_propag}
L_x M \cap N_x^*(K_1 \cap K_2) = 
N_x^* K_1 \cup N_x^* K_2.
\end{align}
Claim (ii) follows again from \cite[Prop. 2.1]{GU1}.
\end{proof}

{\nntext

\begin{remark} The non-linear interaction of waves may cause extraordinary singularities. For
example, M. Beals  showed for the wave equation 
$$\square u(t,y)+b(t,y)u(t,y)^3=0$$ in $\R^{1+3}$ that there are solutions for which the singular support of the Cauchy data $(u|_{t=0},\p_tu|_{t=0})$  is the point $\{0\}$,
but the singular support of $u$  contains the entire solid cone $\{(t,y)\in \R^4;\ |y|<t\}$,
see  \cite[Thm.\ 2.10]{Beals} and  \cite{Beals-spread}. 
% 
% This happens e.g.\ when both the plane wave $u_4$  and the conic wave $w_{321}=
%\square_g^{-1}(au_3\, \square_g^{-1}(au_2u_1))$, produced by the interaction
%of three waves $u_1,u_2$ and $u_3$ (see Fig.\ 3), propagate along the same geodesic $\gamma_{x_4,\xi_4}\subset K_4$.
%In this case it may be that the wave front set of $u_4$  contains a point $(x,\zeta)\in N^*K_4$ and 
%the wave front set of $w_{321}$  contains the point $(x,-\zeta)\in N^*K_4$ with the opposite direction. In this case
%it is difficult to analyze the product $u_4w_{321}$. 
This is related to the case when the condition in the claim (i) of Lemma \ref{lem_2_waves}
is violated.
%This difficulty, as well as the possible caustics 
%of $w_{321}$ (see Fig.\ 6),  are the reasons why in the claim (ii) below we restrict ourselves to a geometrically nice case.
\end{remark}
}

\begin{definition}
We say that submanifolds 
$
K_1, \dots, K_n \subset M
$ 
are transversal if 
for any disjoint subsets $I,J \subset \{1,\dots,n\}$
the manifolds $\bigcap_{i \in I} K_i$ and 
$\bigcap_{j \in J} K_j$ are transversal.
\end{definition}
Observe that if $K_1,\dots,K_4 \subset M$ are transversal submanifolds of codimension 1, then Lemma \ref{lem_transversality} implies for any $x \in \bigcap_{j=1}^4 K_j$ that
$$
T_x^* M = N_x^* K_1 \oplus \dots \oplus N_x^* K_4.
$$

\def\X{\mathcal X}
\def\U{\mathcal U}
\def\L{\mathcal L}
\begin{lemma}[Threefold interaction]
\label{lem_3_waves}
Let $K_1,K_2,K_3 \subset M$ be transversal $\Lambda$-invariant submanifolds of codimension 1. Let $u_j \in I(K_j)$, $j=1,2,3$, and $Q_1, Q_2 \in I(\Delta, \Lambda)$.
Then 
$$
\WF(Q_2(u_3 Q_1(u_2 u_1))) \subset \X \cup \L \cup \U
$$ 
where
\begin{align}\label{def_sets_XLU}
&\X = N^*(K_1 \cap K_2 \cap K_3), \quad \L = \Lambda(\X), \quad
\\\notag&\U = \bigcup_{\tau \in S_3} 
\left( N^* K_{\tau(1)} \cup N^*(K_{\tau(2)} \cap K_{\tau(3)}) \right).
\end{align}
\end{lemma}
\begin{proof}
By Lemma \ref{lem_2_waves}, and using \cite[Th. 8.2.10]{H1} as in its proof, 
$$
\WF(u_3 Q_1(u_2 u_1)) \subset \X \cup \U.
$$
By (\ref{no_new_propag}) it holds that $\Lambda(\U \cap \Sigma) \subset \U$,
and the claim follows as in Lemma \ref{lem_2_waves}.
\end{proof}

We denote by $\pi : T^* M \to M$ the canonical projection to the base manifold.
The equation (\ref{no_new_propag}) implies that 
the sets in the previous lemma satisfy
\begin{align}
\label{new_propag_3way}
\pi(\X \cup \L \cup \U) = \pi(\L) \cup K_1 \cup K_2 \cup K_3.
\end{align}

\begin{lemma}
\label{lem_F4_WF_bound}
Let $K_j \subset M$, $j=1,2,3,4$, be transversal $\Lambda$-invariant submanifolds of codimension 1, and
define for $\tau \in S_4$
the sets  
\begin{align}
\label{bad_set_K}
\mathcal K = 
\bigcup_{\tau \in S_4} \pi(\L_\tau) \cup \bigcup_{j=1}^4 K_{j},
\quad
\L_\tau = \Lambda(\X_\tau),
\quad 
\X_\tau = N^*(\bigcap_{j=1}^3 K_{\tau(j)}).
\end{align}
Suppose that $\bigcap_{j=1}^4 K_j = \emptyset$
and that $v_{j} \in I(K_{j}; \Sym \oplus \R^L)$.
Then for a solution $v_{1234}$ of (\ref{eq_v1234})
it holds that 
$$
\singsupp(v_{1234}) \subset \mathcal K.
$$
\end{lemma}
\begin{proof}
We consider only the singularities caused by the intersection
$$
X = K_1 \cap K_2 \cap K_3,
$$ 
the other cases are similar.
We use the notation from Lemma \ref{lem_3_waves}.
As the four manifolds $K_j$ do not all intersect, 
the distribution $v_4$ is smooth near $X$. 
Then, near $X$, all the singularities of 
$$
v_{1234} = 
\Box_{\hat g, \hat \phi}^{-1} \I_4(v_1, v_2, v_3, v_4) + \Box_{\hat g, \hat \phi}^{-1} \RR_4(v_1, v_2, v_3, v_4)
$$
are due to terms of the form 
$Q_2(u_3 Q_1(u_2 u_1))$ where $Q_1$ and $Q_2$ are in $I(\Delta,\Lambda)$, and $u_j \in I(K_{\tau'(j)})$, $j=1,2,3$, for some $\tau' \in S_3$.
Therefore
Lemma \ref{lem_3_waves} and (\ref{new_propag_3way}) imply that  $\singsupp(v_{1234}) \cap X$ is a subset of 
$\mathcal K$.

There may be interactions also away from $X$ since $\pi(\mathcal L)$ may intersect $K_4$.
These interactions represented by a sum of terms like 
 $$
w_1 = P(v_{123}) v_{4},
\quad 
w_2 = P(v_{4})v_{123}, 
\quad
w_3 = -A^2(D v_{123}, D v_{4})
$$
in (\ref{fourfold_expansion}).
Lemma \ref{lem_3_waves} implies that $\singsupp(v_{123}) \subset \pi(\mathcal L) \cup K_1 \cup K_2 \cup K_3$,
and hence by 
Lemma \ref{lem_2_waves} 
$$
\singsupp(w_j) \subset \singsupp(v_{123}) \cup \singsupp(v_{4}) 
\subset 
\pi(\mathcal L) \cup K_1 \cup \dots \cup K_4 \subset 
\mathcal K,
$$
as long as $(x',\xi') \in \WF(v_{123})$ implies $(x',-\xi') \notin \WF(v_{4})$. 
Note that $(x',\pm\xi')$ lie on the image $\beta$ of the same bicharacteristic. If $(x',\xi') \in \mathcal L$ then $\pi(\beta)$ intersects $X$, and if $(x',-\xi') \in N^* K_4$ then $\beta \subset N^* K_4$ by the $\Lambda$-invariance. 
But this implies that $K_4$ intersects $X$ which is a contradiction with the assumption  $\bigcap_{j=1}^4 K_j = \emptyset$.
\end{proof}

\subsection{Principal symbols of fourfold interaction terms}
\label{sec_abs_interaction}

We refer to \cite[Lem. 3.3]{Lassas2016} for a proof of the following lemma.
We ignore the $2\pi$ factors that do not play any role in the analysis below. 

\begin{lemma}
\label{lem_product}
Let $M$ be a smooth manifold and let $K_1, K_2 \subset M$ be transversal submanifolds. 
Let $u_j \in I(K_j)$, $j=1,2$,
and let $\chi$ be a pseudodifferential operator of order zero that 
vanishes in a conical neighbourhood of $N^* K_1 \cup N^* K_2$.
Then 
$\chi(u_1 u_2) \in I(K_1 \cap K_2)$
and for all $(x,\xi) \in N^*(K_1 \cap K_2)$ it holds that
\begin{align*}
\sigma[\chi(u_1 u_2)](x,\xi) &= 
\sigma[\chi](x,\xi)\, \sigma[u_1](x,\xi_1)\, \sigma[u_2](x,\xi_2),
\end{align*}
where $\xi = \xi_1 + \xi_2$, $\xi_j \in N_x^* K_j$.
\end{lemma}
%\begin{proof}
%Follows from the proof of Lemma 1.1 in \cite{GU}.
%\end{proof}

\begin{lemma}[Type $1+3$ fourfold interaction]
\label{lem_type1}
Let $K_1,\dots,K_4 \subset M$  
be transversal $\Lambda$-invariant submanifolds of codimension 1. Let $u_j \in I(K_j)$, $Q_1, Q_2 \in I(\Delta, \Lambda)$ and let 
\begin{align}
\label{x_cond}
x \in \bigcap_{j=1}^4 K_j \setminus \pi (\L \setminus \X),
\end{align}
where $\X$ and $\L$ are defined by (\ref{def_sets_XLU}).
Let 
$\xi \in L_x M \subset N_x^* K_1 \oplus \dots \oplus N_x^* K_4$ be non-zero,
and define 
$\eta_\ell = \sum_{j=1}^\ell \xi_j$, $\ell = 1,2,3,4$,
where $\xi_j \in N_x^* K_j$ give the direct sum decomposition
$\xi = \sum_{j=1}^ 4 \xi_j$.
Suppose that 
\begin{align}
\label{xi_cond}
\xi_j \ne 0,\ j=1,2,3,4.
\end{align}
Write $s_j = \sigma[u_j]$, $j=1,2,3,4$, and
$q_j(y,\eta) = \sigma[Q_j](y, \eta)$, $j=1,2$,
for $y \in M$ and $\eta \in T_y^* M \setminus L_y M.
$
Then there are $h \in I(\{x\})$ and 
$r \in \D'(M)$ such that 
\begin{align}
\label{symb_type1}
&u_4 Q_2( u_{3} Q_1(u_2 u_1)) = h + r, \quad 
(x, \xi) \notin \WF(r), 
\\\notag
&\sigma[h](x,\xi) = 
s_4(x,\xi_4) q_2(x,\eta_{3}) 
s_3(x,\xi_3) 
q_1 (x,\eta_{2}) s_2(x,\xi_2) s_1(x,\xi_1).
\end{align}
\end{lemma}
Above it is assumed that the product 
in (\ref{symb_type1}) is a well defined distribution.
This is the case if the conormal distributions $u_j$ have negative enough orders, or if $N^* K_4$ and $\L$ are transversal.
Note that $Q_j$ is a pseudodifferential operator on $\Delta \setminus \Lambda$ and thus we can consider $\sigma[Q_j]$ to be defined on $T^*M \setminus 0$.
Observe also that the condition (\ref{x_cond}) guarantees
that propagating wave fronts from earlier threefold interactions do not interfere with the fourfold interaction at the point $x$.

\begin{proof}
Observe that $\eta_3 \notin L_x M$.
Indeed, 
$\xi = \eta_3 + \xi_4 \in L_x M$
and 
$\xi_4 \in L_x M$ by the $\Lambda$-invariance of $K_4$.
If also $\eta_3 \in L_x M$ then by the transversality of $K_1, \dots, K_4$ it holds that $\xi_4 = 0$ or $\eta_3 = 0$,
both of which lead to a contradiction with (\ref{xi_cond}). 
Also $\eta_2 \notin L_x M$ as the sum of $\xi_j \in L_x M$, $j=1,2$.

We recall that $\Sigma$ is the characteristic variety defined by (\ref{def_Sigma}).
For $\ell=2,3,4$, let $\chi_\ell$ be a microlocalization near $(x,\eta_\ell)$, that is, $\chi_\ell$ is a pseudodifferential operator of order zero such that $1-\chi_\ell$ is smoothing in a conical neighbourhood of $(x,\eta_\ell)$.
The microlocalizations can be chosen so that $\chi_\ell$ 
is smoothing outside a conical neighbourhood $U_\ell$ of $(x,\eta_\ell)$ satisfying
\begin{itemize}
\item[(i)] $U_\ell \cap \Sigma = \emptyset$, $\ell = 2,3$,
\item[(ii)] $U_\ell \cap \left( N^*K_\ell 
\cup N^*(\bigcap_{j=1}^{\ell-1} K_j) \right) = \emptyset$, $\ell = 2,3,4$.
\end{itemize}
Indeed, (i) can be fulfilled since $\eta_\ell \notin L_x M$, $\ell = 2,3$, and (ii) since 
$$
\eta_\ell \notin N_x^*K_\ell \cup N_x^*(\bigcap_{j=1}^{\ell-1} K_j), \quad \ell=2,3,4.
$$
The latter holds since $\eta_\ell \in N_x^* K_\ell$
leads to the contradiction $\xi_j = 0$, $j = 1,\dots,\ell-1$,
and since $\eta_\ell \in N_x^*(\bigcap_{j=1}^{\ell-1} K_j)$
leads to the contradiction $\xi_\ell = 0$.

Moreover, the microlocalization $\chi_4$ can be chosen 
to be of the form $$\chi_4(x,D) = \chi_{42}(x,D) \chi_{41}(x),$$
where $\chi_{41} \in C_0^\infty(M)$ satisfies
\begin{itemize}
\item[(iii)] if $y \in \supp(\chi_{41})$ and $(y,\eta) \in N^* K_4$ then $(y,-\eta) \notin \L$.
\end{itemize}
Indeed, 
for a small enough neighbourhood $U \subset M$ of $x$,
the assumption (\ref{x_cond}) implies that it is enough to consider the part $\L_U$ of $\L$ that flows out from $U$, that is, 
the part
$$
\L_U = \Lambda(N^*(U \cap K_1 \cap K_2 \cap K_3)).
$$
It follows from the transversality of $K_1,\dots,K_4$ that, for small enough $U$, 
$$
(y,\eta) \in N^* K_4 \quad \text{implies} \quad (y,-\eta) \notin \L_U.
$$
Taking $\chi_{41} \in C_0^\infty(U)$, we have (iii).

We omit writing the point $x$ as a parameter below. 
Define 
$$
h = \chi_4 (u_4 Q_2 \chi_{3}( u_{3} Q_1 \chi_{2} (u_2 u_1))).
$$
As $U_2 \cap (N^*K_1 \cup N^* K_2) = \emptyset$, 
Lemma \ref{lem_product} implies that 
$$
\chi_2(u_2 u_1) \in I(K_1 \cap K_2), \quad
\sigma[\chi_2(u_2 u_1)](\eta_2) = s_2(\xi_2) s_1(\xi_1).
$$
As $U_2 \cap \Sigma = \emptyset$, 
$Q_1$ is a pseudodifferential operator in $U_2$. Therefore 
$$
Q_1 \chi_2(u_2 u_1) \in I(K_1 \cap K_2), \quad
\sigma[Q_1\chi_2(u_1 u_2)](\eta_2) = q_1(\eta_2) \sigma[\chi_2(u_1 u_2)](\eta_2).
$$
Analogously, as $U_3 \cap (N^*K_3 \cup N^* (K_1 \cap K_2) \cup \Sigma) = \emptyset$,
\begin{align*}
&Q_2 \chi_3 (u_3 Q_1 \chi_2(u_2 u_1)) \in I(K_1 \cap K_2 \cap K_3), 
\\
&\sigma[Q_2 \chi_3 (u_3 Q_1 \chi_2(u_2 u_1))](\eta_3) = q_2(\eta_3) s_3(\xi_3) \sigma[Q_1\chi_2(u_1 u_2)](\eta_2),
\end{align*}
and as $U_4 \cap (N^*K_4 \cup N^* (K_1 \cap K_2 \cap K_4)) = \emptyset$,
$$
h \in I(K_1 \cap K_2 \cap K_3 \cap K_4), \quad
\sigma[h](\xi) = s_4(\xi_4) \sigma[Q_2 \chi_3 (u_3 Q_1 \chi_2(u_2 u_1))](\eta_3).
$$

Observe that $u_4 Q_2( u_{3} Q_1(u_2 u_1)) = h + r_1 + r_2 + r_3$,
where
\begin{align*}
&r_1 = (1-\chi_4)(u_4 Q_2( u_{3} Q_1(u_2 u_1))),
\quad 
r_2 = \chi_4 (u_4 Q_2 (1-\chi_3)(u_3  Q_1(u_2 u_1))),
\\
&r_3 = \chi_4 (u_4 Q_2 \chi_3(u_3 (1-\chi_2) Q_1(u_2 u_1))).
\end{align*}
It remains to show that $(x,\xi) \notin \WF(r_j)$, $j=1,2,3$.
This is immediate for $r_1$.

Let us consider $r_2$. We use the shorthand notation 
$w = Q_2 (1-\chi_3)(u_3  Q_1(u_2 u_1))$
and $\WF_x(w) = \WF(w) \cap T_x M$.
Observe that $\U \cap T_x M \subset \X$ where $\U$ is given by (\ref{def_sets_XLU}).
This together with Lemma \ref{lem_3_waves} and $x \notin \pi(\L \setminus \X)$
imply that $\WF_x(w) \subset \X$.
Using the transversality,
we see that 
\begin{align*}
(y, \eta) \in N^* K_4 \setminus 0 \quad \text{implies} \quad (y,-\eta) \notin \X \cup \U.
\end{align*}
The above, together with (iii) and \cite[Th. 8.2.10]{H1}, imply that 
any $\zeta \in \WF_x(u_4 w)$ has a direct sum decomposition of the form $\zeta = \zeta_4 + \omega_3$ where
\begin{align}
\label{r2_split}
\zeta_4 \in \WF_x(u_4) \cup 0 \subset N_x^* K_4,\quad \omega_3 \in \WF_x(w) \cup 0 \subset N_x^* K_1 \oplus N_x^* K_2 \oplus N_x^* K_3.
\end{align}
As $\eta_3 \notin L_x M$, it holds that $\eta_3 \notin \WF_x(Q_2(1-\chi_3)u)$
for all $u \in \D'(M)$.
In particular, $\eta_3 \notin \WF_x(w)$. 
Thus 
$$
\xi = \xi_4 + \eta_3 \notin \WF_x(u_4 w) \supset \WF_x(r_2).
$$

Let us consider $r_3$. Define 
$w_1 = (1-\chi_2) Q_1(u_2 u_1)$ and $w_2 = Q_2 \chi_3(u_3 w_1)$.
Analogously to (\ref{r2_split}), any 
$\zeta \in \WF_x(u_4 w_2)$ 
decomposes as $\zeta = \zeta_4 + \omega_3$ where
$\zeta_4 \in N_x^* K_4$ and $\omega_3 \in \WF_x(w_2) \cup 0$.
Suppose now that $\omega_3 \notin L_x M$.
Then $\omega_3 \in \WF(u_3 w_1)$, and 
Lemma \ref{lem_2_waves} and \cite[Th. 8.2.10]{H1} imply that 
$\omega_3 = \zeta_3 + \omega_2$ where
$\zeta_3 \in N_x^* K_3$ and $\omega_2 \in \WF_x(w_1) \cup 0$.
As $\eta_3 \notin L_x M$, $\eta_2 \notin \WF_x(w_1)$ and $\eta_2 \ne 0$, it holds that 
$$
\xi = \xi_4 + \eta_3 = \xi_4 + \xi_3 + \eta_2 \notin \WF_x(u_4 w_2)
\supset \WF_x(r_3).
$$
\end{proof}

An analogous proof gives the following.

\begin{lemma}[Type $2+2$ fourfold interaction]
\label{lem_type2}
Let $K_1,\dots,K_4 \subset M$  
be transversal $\Lambda$-invariant submanifolds of codimension 1. Let $u_j \in I(K_j)$
and $Q_1, Q_2 \in I(\Delta, \Lambda)$.
Let $(x,\xi) \in T^*M$,
and use the notation from Lemma \ref{lem_type1}.
Suppose that (\ref{xi_cond}) holds.
Then there are $h \in I(\{x\})$ and 
$r \in \D'(M)$ such that 
\begin{align}
\label{symb_type2}
&Q_2(u_4 u_3) Q_1(u_2 u_1) = h + r, \quad 
(x, \xi) \notin \WF(r), 
\\\notag
&\sigma[h](x,\xi) = 
q_2(x,\xi_3 + \xi_4) s_4(x,\xi_4) s_3(x,\xi_3) q_1(x,\xi_1 + \xi_2) s_2(x,\xi_2)s_1(x,\xi_1).
\end{align}
\end{lemma}

\subsection{Principal symbol of the fourfold interaction}
\label{sec_prin4}

Let $(x,\xi) \in T^*M \setminus 0$, $f \in \mathcal D'(M)$,
and suppose that there are $h \in I(\{x\})$ and $r \in \mathcal D'(M)$
such that $f = h + r$ and $(x,\xi) \notin \WF(r)$. Then we define $\sigma[f](x,\xi) = \sigma[h](x,\xi)$.

Let $K_1,\dots,K_4 \subset M$ be
transversal $\Lambda$-invariant submanifolds of codimension 1,
$x \in \bigcap_{j=1}^4 K_j$, $\xi \in L_x M \setminus 0$
and 
\begin{align}\label{xi_decomp}
\xi = \xi_1 + \xi_2 + \xi_3 + \xi_4 \quad \text{where} \quad
\xi_j \in N_x^* K_j.
\end{align}
Suppose that the covectors $\xi_j$
 satisfy (\ref{xi_cond}) and
\begin{align}\label{x_cond_sym}
x \in \bigcap_{j=1}^4 K_j \setminus 
\bigcup_{\tau \in S_4}
\pi (\L_\tau \setminus \X_\tau),
\end{align}
where $\X_\tau$ and $\L_\tau$
are as in (\ref{bad_set_K}).
Note that (\ref{x_cond_sym}) is a version of the global geometric condition (\ref{x_cond})
that is symmetric with respect to the indices $1,2,3,4$.

We denote by $L_x^+ M$ the future-pointing covectors in $L_x M$.
If $b_1, \dots, b_4 \in L_x^+ M$ are such that 
\begin{align}\label{def_bj}
N_x^* K_j = \R b_j, 
\quad j=1,2,3,4,
\end{align}
then for fixed $(x,\xi)$, the covectors $\xi_j$ 
defined by (\ref{xi_decomp}) 
can be viewed as functions of $\vec b = (b_1, \dots, b_4)$, and we write
\begin{align}\label{xi_as_fun_of_b}
\xi_j = \xi_j(\vec b).
\end{align}
Furthermore, in view of the 
constraint (\ref{def_S_Y}) for the symbols of solutions to (\ref{ein12_linearized}), we write 
$$
\SSym_j(\vec b) = (\SSym_{K_j})_{(x,\xi_j(\vec b))}.
$$

Supposing that $v_j \in I(K_j; \Sym \oplus \R^L)$ have negative enough orders, 
the interaction terms $\I_4$ 
and $\RR_4$
in (\ref{fourfold_expansion}) are well-defined. 
By applying Lemmas \ref{lem_type1} and \ref{lem_type2}
componentwisely, we see that there are $r_1, r_2 \in \mathcal D'(M)$ such that $(x,\xi) \notin \WF(r_j)$, %$j=1,2$, 
and 
$
\I_4(v_1, v_2, v_3, v_4) - r_1$ and $\RR_4(v_1, v_2, v_3, v_4) - r_2$ are in $I(\{x\}; \Sym \oplus \R^L)$.
Moreover, $\RR_4(v_1, v_2, v_3, v_4) - r_2$ is one degree smoother as a conormal distribution than $\I_4(v_1, v_2, v_3, v_4) - r_1$. 
Thus 
\begin{align}\label{interaction_prinsymb}
\sigma[\Box_{\hat g, \hat \phi} v_{1234}](x,\xi) 
= \sigma[\I_4(v_1, v_2, v_3, v_4)](x,\xi).
\end{align}
The principal symbol at $(x,\xi)$ can be viewed as the map 
\def\sigmaJ{\mathscr S}
\begin{align}
\label{map_nonvanish_symb}
\sigmaJ(\vec b, \vec s) = \sigma[\I_4(v_1, v_2, v_3, v_4)](x,\xi), \quad (\vec b, \vec s) \in X_1,
\end{align}
where $\vec s = (s_1, s_2, s_3, s_4)$, $s_j = \sigma[v_j](x,\xi_j)$, 
$v_j$ is as in (\ref{def_vj_deriv}),
$\xi_j$ is as in (\ref{xi_as_fun_of_b}), and 
\begin{align*}%\label{}
X_1 &= \{ (\vec b, s_1, s_2, s_3, s_4);\ s_j \in  \SSym_j(\vec b),\ \vec b \in X_0 \},
\\
X_0 &= \{ \vec b \in (L_x^+ M)^4;\ \text{$\vec b$ is a basis of $\R^4$ and}
\\&\qquad\qquad
\text{$\xi_j = \xi_j(\vec b)$, $j=1,2,3,4$, satisfy (\ref{xi_cond})} \}.
\end{align*}
Note that transversality of $K_1, \dots, K_4$
implies that $\vec b$ is a basis. Conversely, if 
$\vec b \in (L_x^+ M)^4$ is a basis then there are 
transversal $\Lambda$-invariant submanifolds of codimension 1, such that (\ref{def_bj}) holds. 
An explicit construction is given in the proof of Theorem \ref{th_detection} below.

\begin{proposition}
\label{prop_nonvanish_symb}
Let $(x, \xi) \in \Sigma$.
There is an open and dense subset $\mathcal U$ of $X_1^6$
such that for all $(\vec b^n, \vec s^n)_{n=1}^6 \in \mathcal U$
the span of $\sigmaJ(\vec b^n, \vec s^n)$, $n=1,\dots,6$,
is a $6$-dimensional subspace of $
\Sym_x \subset
(\Sym \oplus \R^L)_x$. 
\end{proposition}

The proof of Proposition \ref{prop_nonvanish_symb} is divided in two steps. We will first show that 
there is a point $(\vec b^n, \vec s^n)_{n=1}^6 \in X_1^6$
such that the span of $\sigmaJ(\vec b^n, \vec s^n)$, $n=1,\dots,6$,
is $6$-dimensional.
Then we use an argument based on the fact 
that $\sigmaJ$ can be viewed as a meromorphic function.

\subsubsection{A choice of incoming directions and principal symbols}

We use normal coordinates at $x$ with respect to $\hat g$
in the computations below. In these coordinates, $\hat g$ is the Minkowski metric and 
$\hat \Gamma_{pq}^r = 0$ at $x$, see e.g. \cite[Prop. 33, p. 73]{ONeill1983}.
Moreover, we choose the coordinates so that $\xi = (1,1,0,0)$.

We begin by choosing the symbols $\vec s^n = (s_1,\dots,s_4)$ so that they do not depend on $n$, and write 
$$
s_j = (g_{(j)}, \phi_{(j)}), \quad j = 1,2,3,4,
$$
where 
$s_j = \sigma[v_j](x,\xi_j)$, 
$v_j$ is as in (\ref{def_vj_deriv}), and 
the parenthesis are to emphasize that the indices are not tensor indices. 
We will use an analogous notation for the higher order derivatives in $\epsilon$, that is, 
$$
\sigma[v_{jk}](x, \eta_{jk}) =  (g_{(jk)}, \phi_{(jk)}),
\quad
\sigma[v_{jkl}](x,\eta_{jkl}) = 
 (g_{(jkl)}, \phi_{(jkl)}),
$$
where $\eta_{jk} = \xi_j + \xi_k$ and $\eta_{jkl} = \xi_j + \xi_k + \xi_l$.
Here $v_{jk}$ and $v_{jkl}$ are defined by (\ref{def_vderivs}), and we emphasize again that the indices are not tensorial.

We choose 
$g_{(j)} = 0$ for all $j=1,2,3,4$, and 
$$
\phi_{(1)} = \phi_{(2)}
= (1,0,\dots,0) \in \R^L,
\quad 
\phi_{(3)} = \phi_{(4)} 
= (0,1,0,\dots,0) \in \R^L.
$$
Observe that then $s_j \in \SSym_j(\vec b)$ for any $\vec b$.
Let us use the shorthand notations
$$
q(\eta) = 1/\hat g(\eta,\eta), \quad 
\xi \hat \otimes \eta = \xi \otimes \eta + \eta \otimes \xi.
$$
Differentiating (\ref{einred_full}) with respect to 
$\epsilon_j$ and $\epsilon_k$, $j \ne k$, and evaluating the principal symbol at $\epsilon = 0$, we see that when $x \notin \supp(F)$,
$$
g_{(jk)} = -2 q(\eta_{jk})\, \xi_j \hat \otimes \xi_k,
\quad jk = 12,34.
$$ 
Moreover, 
$g_{(jk)} = 0$ when $jk \ne 12,34$, and 
$\phi_{(jk)} = 0$ for all pairs $jk$.

We write 
$$
\Gamma_{\alpha \gamma \beta}(\eta)g =  \frac 1 2 ( \eta_\alpha g_{\gamma \beta} + \eta_\beta g_{\alpha \gamma} - \eta_\gamma g_{\alpha \beta} ), \quad
\Gamma_{\alpha \gamma \beta}\, g_{(jk)} = 
\Gamma_{\alpha \gamma \beta}(\eta_{jk})\, g_{(jk)},
$$
and define also 
$$
p(g_{(jk)}) \phi_l = \hat g g_{(jk)} \hat g(\eta_l, \eta_l) \phi_l, 
\quad 
p(g_{(jk)}) g_{(lm)} = \hat g g_{(jk)} \hat g(\eta_{lm}, \eta_{lm}) g_{(lm)}. 
$$
%A straighforward computation shows that 
%\begin{align}
%\label{gamma_contraction_vanishes}
%\hat g^{\alpha \beta} \Gamma_{\alpha \gamma \beta}\, g_{(jk)} = 0, \quad jk=12,34,\ \gamma=0,1,2,3.
%\end{align}
Recalling (\ref{Box_phi_subprin}), we differentiate (\ref{einred_full}) with respect to 
$\epsilon_j$, $\epsilon_k$, $\epsilon_l$, where $j$, $k$ and $l$ are distinct, and evaluate the principal symbol at $\epsilon = 0$. This yields for $x \notin \supp(F)$,
$$
\phi_{(12j)} = q(\eta_{12j}) p(g_{(12)}) \phi_{(j)}, 
\quad
\phi_{(34k)} =  q(\eta_{34k}) p(g_{(34)}) \phi_{(k)},
$$
where $j=3,4$, $k=1,2$.
Moreover, $\phi_{(jkl)} = 0$ for the other triplets $jkl$,
and $g_{(jkl)} = 0$ for all triplets.

Analogously, we see that the scalar field components of $\sigma[\I_4(v_1, v_2, v_3, v_4)](x,\xi)$ vanish 
and the metric tensor part is
\begin{align}
\label{prinsymb_full}
&p(g_{(12)}) g_{(34)}
+ p(g_{(34)}) g_{(12)}
-2 \eta_{123} \hat\otimes \xi_4
-2 \eta_{124} \hat\otimes \xi_3
-2 \eta_{341} \hat\otimes \xi_2
-2 \eta_{342} \hat\otimes \xi_1
\\\notag&\quad+ 
2 \hat g^{pq} \hat g^{rs} (\Gamma_{prj}g_{(12)}\, \Gamma_{qsk}g_{(34)}
+ \Gamma_{prj}g_{(12)}\, \Gamma_{qks}g_{(34)} + \Gamma_{prk}g_{(12)}\, \Gamma_{qjs}g_{(34)})
\\\notag&\quad+ 
2 \hat g^{pq} \hat g^{rs} (\Gamma_{prj}g_{(34)}\, \Gamma_{qsk}g_{(12)}
+ \Gamma_{prj}g_{(34)}\, \Gamma_{qks}g_{(12)} + \Gamma_{prk}g_{(34)}\, \Gamma_{qjs}g_{(12)}).
\end{align}

Let us now choose the incoming directions $\vec b^n = (b_1^n,\dots,b_4^n)$. We will choose each $b_j^n$ to be in the following set of six Pythagorean quadruples
$$
B = \{ (p^2 + m^2 + n^2, 2 m p, 2 n p, 
 p^2 - (m^2 + n^2));\ m=1,2,3,\ n=1, 2,\ p=1 \}.
$$
Then we consider the following choices
$$
\mathbb B = \{\vec b = (b_1, b_2, b_3, b_4);\ b_j \in B,\ \text{$\vec b \in X_0$} \}.
$$
The set $\mathbb B$ contains 9 elements, up to reindexing. Indeed, out of the 15 subsets of $B$ with 4 elements, 6 fail to be bases or satisfy (\ref{xi_cond}). 
A straightforward computation shows that the principal symbols (\ref{prinsymb_full}) corresponding to the incoming directions in $\mathbb B$
span a 6-dimensional subspace of the fibre $\Sym_x$.
We have verified this using a computer algebra software.
%see the supplement TODO. 
%\HOX{Add a ref to a supplement. Should we put this e.g. in github?}
All the computations were performed using rational arithmetic.
Thus there is a point $(\vec b^n, \vec s^n)_{n=1}^6 \in X_1^6$
such that the span of $\sigmaJ(\vec b^n, \vec s^n)$, $n=1,\dots,6$,
is $6$-dimensional and contained in $\Sym_x$.
This completes the first step of the proof of Proposition \ref{prop_nonvanish_symb}.

\subsubsection{Changing incoming directions}

We choose local coordinates and identify $T_x^* M$ with $\R^4$.
Consider $\vec b \in X_0$ and write
\begin{align}
\label{def_xij_rj}
\xi = \xi_1 + \xi_2 + \xi_3 + \xi_4, \quad \xi_j = r_j(\vec b) b_j.
\end{align}
As the coefficients $r_j(\vec b)$ can be obtained by solving a system of linear equations (for example, using Cramer's rule), we have that 
$$
r_j(\vec b) = \frac{p_j(\vec b)}{q(\vec b)},
\quad j = 1,2,3,4,
$$
where $p_j$ and $q$ are polynomials in $\vec b \in (\R^4)^4$
that do not vanish at any point on $X_0$.
We can also consider $p_j$, $j=1,2,3,4$, and $q$ as functions on $\tilde X_0 = (L_+ M)^4$.
Here $\tilde X_0$ is a connected, real-analytic manifold, and
$p_j$ and $q$ are real-analytic functions on $\tilde X_0$.
Observe also that the equation (\ref{def_S_Y}) is linear in $h$ and 
of full rank when $\xi \in L_x^+ M$,
and therefore the manifold 
$$
\tilde X_1 = \{ (\vec b, s_1, s_2, s_3, s_4);\ s_j \in  \SSym_j(\vec b),\ \vec b \in \tilde X_0 \}
$$
is real-analytic, see e.g. \cite[Prop. 1.9.2]{Krantz1992}. 

The expression (\ref{fourfold_expansion})
is of the form $P_0(\vec \xi, \vec s)/Q_0(\vec \xi)$ where $P_0$ and $Q_0$ are polynomials in $\vec \xi = (\xi_1, \xi_2, \xi_3, \xi_4) \in (\R^4)^4$. 
Here $P_0$ takes values on $(\Sym \oplus \R^L)_x = \R^{10+L}$
and $Q_0$ on $\R$.
Moreover, 
for $\vec \xi$ as in (\ref{def_xij_rj}), with $\vec b \in X_0$,
it holds that $Q_0$ is non-zero. It follows that 
$$
\sigmaJ(\vec b, \vec s) = \frac{P_1(\vec b, \vec s)} {Q_1(\vec b)},
$$
where $P_1$ and $Q_1$ are real-analytic on $\tilde X_1$ and $\tilde X_0$, respectively, and $Q_1$ is nowhere vanishing on $X_0 \subset \tilde X_0$.

In the previous section, we showed that there is 
$(\vec b^n, \vec s^n)_{n=1}^6 \in X_1^6$
such that the span of $\sigmaJ(\vec b^n, \vec s^n)$, $n=1,\dots,6$,
is $6$-dimensional subspace
of the fibre $(\Sym \oplus \R^L)_x$. 
We set $e_n = \sigmaJ(\vec b^n, \vec s^n)$, $n=1,\dots,6$,
and choose  
$$
e_n \in (\Sym \oplus \R^L)_x, \quad n=7, \dots, 10+L,
$$
such that
the vectors $e_n$, $n=1, \dots, 10+L$ 
form a basis of $(\Sym \oplus \R^L)_x$.
Let $[e_1,\dots,e_{L+10}]$ be the matrix with
the columns $e_n$,
and consider the map 
\begin{align*}
%\label{basis_det}
(\vec b^n, \vec s^n)_{n=1}^6
\mapsto 
\det([\sigmaJ(\vec b^1, \vec s^1), \dots,
\sigmaJ(\vec b^6, \vec s^6), e_{7}, \dots, e_{L+10}]).
\end{align*}
This can be written as 
\begin{align*}
\frac{\det([P_1(\vec b^1, \vec s^1), \dots,
P_1(\vec b^6, \vec s^6), e_{7}, \dots, e_{L+10}])}
{\prod_{n=1}^6 Q_1(\vec b^n)}
= \frac{P_2((\vec b^n, \vec s^n)_{n=1}^6)}{Q_2((\vec b^n, \vec s^n)_{n=1}^6)}
\end{align*}
where $P_2$ and $Q_2$ are real-analytic functions on the connected, real-analytic manifold $\tilde X_1^6$ that do not vanish at the point $(\vec b^n, \vec s^n)_{n=1}^6$ considered in the previous section.
Thus the unique continuation principle for real-analytic functions on real-analytic manifolds, see e.g. \cite[Lem. VI.4.3]{Helgason1962},
implies that there is an open and dense set $\mathcal U \subset \tilde X_1^6$
such that both $Q_2$ and $P_2$ are nowhere vanishing on $\mathcal U$.
This finishes to proof of Proposition \ref{prop_nonvanish_symb}.
\hfill$\square$

\section{Construction of light observation sets}
\label{sec_reduction_to_los}

In this section we prove the main theorem formulated in Section \ref{sec_intro}.
We begin by relating the singular solutions constructed in Section \ref{sec_construction_singsol} with the data set (\ref{data}), and by showing that singularities due to non-linear interactions of singular solutions, as studied in Section \ref{sec_interaction}, can be detected from (\ref{data}).
Then we relate the detected singularities with geometric objects on $(M, \hat g)$, the earliest light observation sets. 
This allows us to reduce the reconstruction of the Lorentzian structure to a purely geometric problem solved in \cite{KLU-august}.

\subsection{Sending singularities}

Let us show how to use the data set (\ref{data}) to find families of sources 
parametrized by $\epsilon \in \mathcal E$ as in (\ref{lin_family}).
Recall that (\ref{data})
gives sources $\mathcal F$ in the Fermi coordinates $\Phi_g$
that depend on the source itself. For this reason,
we need to construct first common coordinates in order to parametrize the sources smoothly.

We will use the wave maps.
For a metric tensor $g$ close to $\hat g$,
the wave map operator is defined in local coordinates by
$$
\square_{g,\hat g} \Psi^p
= 
g^{jk}(\frac \p{\p x^j}\frac \p{\p x^k} \Psi^p -\Gamma^{n}_{jk}(x)\frac \p{\p x^n}\Psi^p+
\hat \Gamma^p_{BC}(\Psi(x))
\frac \p{\p x^j}\Psi^B\,\frac \p{\p x^k}\Psi^C),
$$
where $\hat \Gamma^p_{BC}$ and $\Gamma^{j}_{kl}$ are the Christoffel symbols of the metrics $\hat g$
and $g$, respectively.
This operator invariant under diffeomorphisms in the following sense: if 
$\Box_{g,\hat g} \Psi = 0$ and if 
$\psi$ and $\mu$ are diffeomorphisms, then setting 
$\tilde \Psi = \mu \circ \Psi \circ \psi$,
it holds that $\Box_{\psi^* g, \mu_* \hat g} \tilde \Psi = 0$, where $\mu_* = (\mu^{-1})^*$, the pushforward by $\mu$, see e.g. \cite{Choquet-Bruhat2000}.

%Recall that $\Phi : V \to M$ stands for (\ref{def_Fermi}), the Fermi coordinates with respect to the background metric $\hat g$, and that the corresponding coordinates 
%with respect to a metric $g$ are denoted by $\Phi_g$.

%Proposition \ref{prop_nonvanish_symb}, Lemma \ref{lem_wave_coord}

\begin{lemma}
\label{lem_wave_coord}
Let $(\hat g, \hat \phi)$ be as in Theorem \ref{th_main}
and write $\hat g' = \Phi_{\hat g}^* \hat g$. 
Let $p \in V$, let $B \subset V$ be a 
neighbourhood of $p$ satisfying $\overline B \subset V$,
and let $q \in V$ satisfy $J_{(V,\hat g')}(p,q) \subset B$.
If $k \in \N$ is large and $r > 0$ is small, then 
for all 
$(g',\phi',\mathcal F')$ in $\mathcal D_{r,k}(\hat g, \hat \phi)$
satisfying $\supp(\mathcal F') \subset J_{(V,g')}^+(p)$,
the equation
\begin{align}
\label{wave_map_V}
&\square_{g',\hat g'} \Psi=0\quad \text{in $J_{(V,g')}(p,q)$},
\\\notag
&\Psi= \Phi_{\hat g}^{-1} \circ \Phi_g,\quad \text{in $B \setminus J_{(V,g')}^+(p)$},
\end{align}
has a unique solution that coincides with $\Phi_{\hat g}^{-1} \circ \Phi_g$ on $J_{(V,\hat g')}(p,q)$.
Here $g$ is a metric that satisfies 
\begin{itemize}
\item[(G1)] $(g',\phi',\mathcal F') = \Phi_{g}^* (g, \phi, {\mathcal F})$,
\item[(G2)] $(g, \phi, {\mathcal F})$
satisfies (\ref{eq1})--(\ref{eq2}) on $\hat M$,
and (\ref{init_cond}) holds,
\item[(G3)] $\Ric_{\hat g}(g) = \Ric(g)$.
\end{itemize}
Moreover, for small enough $B$, the equation (\ref{wave_map_V}) can be solved 
given the data set $\mathcal D_{r,k}(\hat g, \hat \phi)$.
\end{lemma}
\begin{proof}
As $(g',\phi',\mathcal F') \in \mathcal D_{r,k}(\hat g, \hat \phi)$ there is $(\tilde g, \tilde \phi, \tilde{\mathcal F})$
on $\hat M$
such that (G1) and (G2) hold for this triple.
Let $\tilde T > \hat T$
and write $\tilde M = (0,{\tilde T}) \times N$.
For small $r$ and large $k$, $(\tilde g,\tilde \phi,\tilde{\mathcal F})$ can extended so that it solves the equation (\ref{eq1})--(\ref{eq2}) on $\tilde M$ with $\tilde{\mathcal F}$ vanishing on $(\hat T, {\tilde T}) \times N$.
Then, for small $r$ and large $k$, 
the solution $\tilde \Psi$ of 
\begin{align}
\label{wave_map_global}
&\square_{\tilde g,\hat g} \tilde \Psi=0\quad \text{in $\tilde M$},
\\\notag
&\tilde \Psi(x)=x,\quad \text{in $(-\infty, 0) \times N$},
\end{align}
is a diffeomorphism and $\hat M \subset \tilde \Psi(\tilde M)$.
Moreover, the reduced Ricci tensor of $g = \tilde \Psi_* \tilde g$ coincides with the actual Ricci tensor
in the sense of (G3), see e.g.
\cite[Sec. VI.7.2 and Th. 4.2 in Appendix\ III]{Choquet-Bruhat2009}.
Setting $(g,\phi,\mathcal F) = \tilde \Psi_* (\tilde g, \tilde \phi, \tilde{\mathcal F})$, we have that 
$(g,\phi,\mathcal F)$
satisfies (G1)-(G3).
Note that (G1) holds since the Fermi coordinates coincide for isometric manifolds satisfying the initial condition (\ref{init_cond}). In other words, the Fermi coordinates have the naturality property described the following commuting diagram:
\begin{align*}%\label{}
\xymatrix{
& (V,g') \ar[dl]_{\Phi_{\tilde g}} \ar[dr]^{\Phi_{g}} &      
\\
(\tilde M, \tilde g) \ar[rr]_{\tilde \Psi} && (\tilde \Psi(\tilde M), g) 
}
\end{align*}

Properties (G2) and (G3)
imply that 
\begin{align}\label{wave_map_id}
\square_{g,\hat g} \Id =0, \quad \text{in $\hat M$},
\end{align}
where $\Id(x) = x$.
Indeed, in local coordinates, 
$$
\square_{g,\hat g} \Id^p 
= g^{jk}(-\Gamma^{n}_{jk} \delta_n^p +
\hat \Gamma^p_{BC} \delta_j^B \delta_k^C)
= - g^{pq}(\hat H_q - \Gamma_q),
$$
where $\hat H_q - \Gamma_q$ are the components of the tensor $H_{\hat g}(g)$, see (\ref{def_harm_tensor}).
Moreover, the equality $\Ric_{\hat g}(g) = \Ric(g)$
implies that 
$$\div_g \Ein_{\hat g}(g) = \div_g \Ein(g) = 0,$$
where the reduced Einstein tensor is defined by $\Ein_{\hat g}(g) = I_g \Ric_{\hat g}(g)$.
As $(g, \phi, {\mathcal F})$ satisfies the initial condition (\ref{init_cond}) we may apply Lemma \ref{lem_from_div_to_wgauge}, when $k \in \N$ is large and $r > 0$ is small, to conclude that $H_{\hat g}(g) = 0$.

Now (\ref{wave_map_id}), together with the diffeomorphism invariance of the wave map operator, imply that 
$
\square_{g',\hat g'} (\Phi_{\hat g}^{-1} \circ \Phi_g) = 0.
$
Hence $\Phi_{\hat g}^{-1} \circ \Phi_g$ is the solution of (\ref{wave_map_V}).

Let us now show how to solve (\ref{wave_map_V})
given $\mathcal D_{r,k}(\hat g, \hat \phi)$. We begin by showing that 
$\hat g'$ is determined by $\mathcal D_{r,k}(\hat g, \hat \phi)$, 
which again follows from the fact that $\mathcal F' = 0$ implies $g' = \hat g'$.
Indeed, in this case also $\mathcal F = 0$,
and (G2) and (G3) imply that $(g,\phi,0)$
satisfies (\ref{einred_full}).
But this quasilinear wave equation has a unique solution for small $r$ and large $k$, and therefore $g = \hat g$.
Thus we find $\hat g'$ as the first component the unique triple in 
$\mathcal D_{r,k}(\hat g, \hat \phi)$ with the third component vanishing. 

As $\hat g'$ is known, we can evaluate 
the wave map operator $\square_{g',\hat g'} \Psi$ for any $g'$ on $J_{(V,g')}(p,q)$ and any $\Psi$ close to the identity.
It remains to show that the initial condition in (\ref{wave_map_V})
is determined by $\mathcal D_{r,k}(\hat g, \hat \phi)$.
This follows from the identity 
\begin{align}
\label{wc_init_id}
\exp_{p}^{\hat g'} \circ (\exp_{p}^{g'})^{-1}
= \Phi_{\hat g}^{-1} \circ \Phi_g, \quad \text{in $B \setminus J_{(V,g')}^+(p)$},
\end{align}
that we will show next. Here $\exp^{g'}$
denotes the exponential map with respect to the metric $g'$.

The assumption  $\supp(\mathcal F') \subset J_{(V,g')}^+(p)$
implies that $g = \hat g$
on $\hat M \setminus J_{(M,\hat g)}^+(\Phi(p))$,
using again the the fact that (\ref{einred_full}) has a unique solution.
Here we used the notation $\Phi = \Phi_{\hat g}$.
In particular, it holds that  $\Phi(p) = \Phi_g(p)$, $d\Phi = d\Phi_g$ at $p$, and 
$$
\Id = \exp_{\Phi(p)}^{\hat g} \circ (\exp_{\Phi(p)}^{g})^{-1}
\quad \text{in $U_1 \setminus J_{(M,\hat g)}^+(\Phi(p))$},
$$
for a suitable neighbourhood $U_1 \subset \hat M$ of $\Phi(p)$.
Equation (\ref{wc_init_id}) is obtained by
combining the previous identity with the naturality of the exponential map, that is, with the fact that the following two diagrams commute:
\begin{align*}%\label{}
\xymatrix{
U_2 \ar[r]^{d \Phi} \ar[d]_{\exp_{p}^{\hat g'}}
& 
U_3 \ar[d]^{\exp_{\Phi(p)}^{\hat g}} 
\\
(V, \hat g')  \ar[r]^{\Phi_{\hat g}} 
&
(\hat M, \hat g)
}
\qquad
\xymatrix{
U_2 \ar[r]^{d \Phi} \ar[d]_{\exp_{p}^{g'}}
& 
U_3 \ar[d]^{\exp_{\Phi(p)}^{g}} 
\\
(V, g')  \ar[r]^{\Phi_g} 
&
(\hat M, g)
}
\end{align*}
where $U_2 \subset T_p V$ and $U_3 \subset T_{\Phi(p)} \hat M$
are suitable neighbourhoods of the respective origins.
Indeed, in $B \setminus J_{(V,g')}^+(p)$
$$
\exp_{p}^{\hat g'} \circ (\exp_{p}^{g'})^{-1}
= 
\Phi_{\hat g}^{-1} \circ \exp_{\Phi(p)}^{\hat g} \circ d\Phi
\circ d\Phi^{-1} \circ (\exp_{\Phi(p)}^{g})^{-1} \circ \Phi_g
= \Phi_{\hat g}^{-1} \circ\Phi_g.
$$
\end{proof}

\subsection{Receiving singularities}

The only difficulty in receiving of singularities is that the Fermi coordinates $\Phi_g$
%, defined in Section \ref{sec_intro}, 
can be non-smooth when $g$ has singularities. 
This problem is solved by using the following two lemmas:
the first one shows that at least some singularities 
of the fourfold linearization of $g$
are visible when the principal symbol of the fourfold linearization is non-vanishing in a large enough subbundle
and suitable coordinates are used, and 
the second one shows that changing from the Fermi coordinates to normal coordinates gives suitable coordinates.

In the first lemma, we denote by $\Sym(K)$ the rank 6 subbundle of $\Sym$ 
over a submanifold $K \subset M$ of codimension one, that consists of the symmetric 2-tensors on $K$. 
Later, we will use this lemma when $K$ is a $\Lambda$-invariant submanifold, in fact, in the proof of Theorem \ref{th_detection} below, the submanifold $K$ will be contained 
in the image of $L_x^+ M$ under the exponential map.
Here $x$ is as in Theorem \ref{th_detection}, and there are no cut points on the geodesic segment joining $x$ and $y$,
with $y$ as in the lemma.

\newcommand{\pair}[1]{\left\langle #1 \right\rangle}
\begin{lemma}
\label{lem_detect}
Consider a family $g \in C^4(\mathcal E; C^2(M;\Sym))$ where $\mathcal E$ is a neighbourhood of the origin in $\R^4$, and let $K \subset M$ be a submanifold of codimension one.
Let $y \in K$ and let $U \subset M$ be a small neighbourhood of $y$.
Suppose that for $\alpha = (1,1,1,1)$ it holds that
\begin{align}\label{detect_smoothness_alpha}
\p_\epsilon^\alpha g|_{\epsilon =0} &\in I(K;\Sym),
\\\label{detect_smoothness_beta}
\p_\epsilon^\beta g|_{\epsilon =0} &\in 
C^\infty(U;\Sym), \quad  \text{for all $\beta < \alpha$}.
\end{align}
Here the $\beta < \alpha$ means the partial order in the sense of multi-indices. 
Suppose, moreover, that the restriction of  
$\sigma[\p_\epsilon^\alpha g|_{\epsilon =0}](y,\cdot)$ on $\Sym(K)$ does not vanish.
Let $V \subset \R^4$ be open and let
$\Psi \in C^4(\mathcal E; C^2(V; U))$.
Suppose that $\Psi(\epsilon)$
gives a system of coordinates near $y$ and that 
\begin{align}\label{detect_smoothness_W}
%\p_\epsilon^\alpha \Psi^* g|_{\epsilon =0} \in C^2(W;\Sym),
%\quad
\p_\epsilon^\beta \Psi^* g|_{\epsilon =0} &\in 
C^\infty(V;\Sym), \quad  \text{for all $\beta < \alpha$}.
\end{align}
Then $\p_\epsilon^\alpha \Psi^* g|_{\epsilon =0}$ is not smooth at $\Psi^{-1}(y)$.
Here $(\Psi^* g)(\epsilon,\cdot) = \Psi(\epsilon, \cdot)^* g(\epsilon, \cdot)$.
\end{lemma}
\begin{proof}
For small enough $U$, we can choose such smooth coordinates $(x^0,x^1,x^2,x^3)$ that $K \cap U = \{x \in U; x^0 = 0\}$.
Then the curvature $(0,4)$-tensor can be expressed in the coordinates as follows 
\begin{align*}
R_{ik\ell m} &=\frac{1}{2}\left(
\frac{\partial^2g_{im}}{\partial x^k \partial x^\ell} 
+ \frac{\partial^2g_{k\ell}}{\partial x^i \partial x^m}
- \frac{\partial^2g_{i\ell}}{\partial x^k \partial x^m}
- \frac{\partial^2g_{km}}{\partial x^i \partial x^\ell} \right)
\\&\quad +g_{np} \left(
\Gamma^n{}_{k\ell} \Gamma^p{}_{im} - 
\Gamma^n{}_{km} \Gamma^p{}_{i\ell} \right),
\end{align*}
where $\Gamma^m{}_{ij}$ are the Christoffel symbols of $g$.
We write $g^{(\alpha)} = \p_\epsilon^\alpha g|_{\epsilon =0}$. As $g^{(\alpha)} \in I^p(K;\Sym)$ and as its principal symbol does not vanish, 
there is the maximal smoothness index $s=-p-1$
such that $g^{(\alpha)}$ is in the Besov type space ${^\infty}H_{(s)}(U;\Sym)$, see \cite[Def. 18.2.6]{H3}. Moreover, 
it holds for $i,m = 1,2,3$ that
$$
\frac{\partial^2 g^{(\alpha)}_{k\ell}}{\partial x^i \partial x^m} \in {^\infty}H_{(s)}(U),
\quad 
\frac{\partial^2 g^{(\alpha)}_{i\ell}}{\partial x^k \partial x^m} \in {^\infty}H_{(s-1)}(U), 
\quad 
\frac{\partial^2 g^{(\alpha)}_{km}}{\partial x^i \partial x^\ell} \in {^\infty}H_{(s-1)}(U). 
$$
Thus, writing $R^{(\alpha)} = \p_\epsilon^\alpha R|_{\epsilon = 0}$, 
$$
R^{(\alpha)}_{i 00 m} = \frac{1}{2}
\frac{\partial^2 g^{(\alpha)}_{im}}{\partial x^0 \partial x^0} 
+ r, \quad i,m = 1,2,3,
$$
where $r \in {^\infty}H_{(s-1)}(U)$.
We used here the fact that all the terms in $R^{(\alpha)}$,
where the four derivatives with respect to $\epsilon$
do not act on the same instance of $g$, are smooth. This follows from (\ref{detect_smoothness_beta}).

Consider a family of curves $\gamma \in C^4(\mathcal E; C^1(-1,1;U))$ and suppose that $\gamma(\epsilon)$ is a geodesic with respect to $g(\epsilon)$. 
Write $\hat \gamma = \gamma|_{\epsilon = 0}$
and suppose, moreover, that $\hat \gamma$ intersects $K$ transversally at $y = \hat \gamma(0)$.
We can consider $g^{(\alpha)}_{im} \circ \hat\gamma$
as a conormal distribution associated to $\{0\} \subset (-1,1)$.
As $\sigma[g^{(\alpha)}](\hat \gamma(0),\cdot)$ does not vanish on $\Sym(K)$, there are $i,m = 1,2,3$ such that 
\begin{align}\label{detect_ps}
\sigma[R^{(\alpha)}_{i 00 m} \circ \hat\gamma
](0, \tau)
= - \frac 1 2 \xi_0^2\, \sigma[g^{(\alpha)}_{im} \circ \hat\gamma
](0, \tau), 
\end{align}
does not vanish for some $\tau \in \R \setminus 0$.

By differentiating the geodesic equation $\nabla_{\dot \gamma} \dot \gamma = 0$ with respect to $\epsilon$
we see that the path $\gamma_{(\alpha)} = \p_\epsilon^\alpha \gamma|_{\epsilon = 0}$ satisfies
\begin{align}\label{detect_eq_gamma}
\p_t \dot \gamma_{(\alpha)}^k + (\hat \Gamma_{ij}^k \circ \hat \gamma) \dot{\hat\gamma}^j 
\dot \gamma_{(\alpha)}^i 
+
(\hat \Gamma_{ij}^k \circ \hat \gamma) \dot \gamma_{(\alpha)}^j
\dot{\hat\gamma}^i
+ 
\gamma_{(\alpha)}^\ell (\p_{x^\ell}\hat \Gamma_{ij}^k \circ \hat \gamma) \dot{\hat\gamma}^j
\dot{\hat\gamma}^i
= h^k + r^k,
\end{align}
where $k = 0,1,2,4$, 
$h^k = -(\p_\epsilon^\alpha \Gamma_{ij}^k|_{\epsilon = 0} \circ \hat \gamma) \dot{\hat\gamma}^j 
\dot{\hat\gamma}^i$
and $r^k$ is smooth due to (\ref{detect_smoothness_beta}).
This is a linear system of differential equations for $\gamma_{(\alpha)}$ with smooth coefficients.
There is a maximal smoothness index $s' \in \R$
such that 
$$
g^{(\alpha)}_{jk} \circ \hat \gamma \in {^\infty}H_{(s')}(-1,1),
\quad
h^k \in {^\infty}H_{(s'-1)}(-1,1).
$$
Thus (\ref{detect_eq_gamma}) implies that 
$\gamma_{(\alpha)}^k, \dot \gamma_{(\alpha)}^k \in {^\infty}H_{(s'-1)}(-1,1)$
for $k=0,1,2,3$.

As the intersection of $\hat \gamma$ and $K$ is transversal, 
there are $C^4$-smooth maps $t(\epsilon)$ and $y(\epsilon)$,
mapping a neighbourhood of the origin in $\mathcal E$ to small neighbourhoods of the origin in $\R$ and $y$ in $M$, respectively, 
 such that $\gamma(\epsilon, t(\epsilon)) = y(\epsilon) \in K$.
Let $Z(\epsilon,t)$, $V(\epsilon,t)$ and $W(\epsilon,t)$ be the parallel transports 
of $\p_{x^0}, \p_{x^i}, \p_{x^m} \in T_{y(\epsilon)} M$ along $\gamma(\epsilon)$ with respect to $g(\epsilon)$.
The vector field $Z_{(\alpha)} = \p_\epsilon^\alpha Z|_{\epsilon = 0}$
satisfies an equation similar to (\ref{detect_eq_gamma}),
and analogously with the above, we see that its components satisfy
$Z_{(\alpha)}^k \in {^\infty}H_{(s'-1)}(-1,1)$.
The same is true for $V$ and $W$.
It follows, using again (\ref{detect_smoothness_beta}), that the principal symbol of $\p_\epsilon^\alpha \pair{R(Z, V) Z, W}|_{\epsilon = 0}$ is (\ref{detect_ps}).
As the principal symbol does not vanish, 
we see that the function 
\begin{align}\label{R_nonsmooth}
\p_\epsilon^\alpha \pair{R(Z, V)Z,W}|_{\epsilon = 0}
\end{align}
is not in $C^{\infty}(-t_0,t_0)$ for any $t_0 > 0$.

To get a contradiction suppose that 
there is a neighbourhood $V'$ of $\Psi^{-1}(y)$ such that $\Psi^* \p_\epsilon^\alpha g|_{\epsilon =0}$
is smooth in $V'$. 
Now (\ref{detect_smoothness_W}) together with the smoothness of $\Psi^* \p_\epsilon^\alpha g|_{\epsilon =0}$ implies that 
the function (\ref{R_nonsmooth}) is smooth at $t=0$. But this is a contradiction.
\end{proof}

In the proof of Theorem \ref{th_detection} we will apply the above lemma in a setting where $\Psi$ is given as the composition of the Fermi coordinates $\Phi_g$ and certain normal coordinates. The next lemma captures the essential features of the setting. We will apply it in the case 
$\Phi(\epsilon, y) = \Phi_{g(\epsilon)}(y)$, $y \in V$.

\begin{lemma}
\label{lem_detect_normal_coords}
Let $g \in C^4(\mathcal E; C^3(M;\Sym))$, with $\mathcal E$ a neighbourhood of the origin in $\R^4$,  
let $\Phi \in C^4(\mathcal E; C^2(V; M))$,
and suppose that $\Phi(\epsilon, \cdot)$ gives local coordinates on $M$. Write $g' = \Phi^* g$, that is,
$(\Phi^* g)(\epsilon,\cdot) = \Phi(\epsilon, \cdot)^* g(\epsilon, \cdot)$.
Let $y \in V$, let $X_j$, $j=0,1,2,3$, be a basis 
of $T_y V$, and consider the normal coordinates 
    \begin{equation}\label{norm_coord_Xi}
\Xi(\epsilon, z) = \exp_{y}^{g'(\epsilon)}
\left(\sum_{j=0}^3 z^j X_j\right).
    \end{equation}
Define $\Psi(\epsilon, z) = \Phi(\epsilon, \Xi(\epsilon, z))$. 
Then there is a neighbourhood $V'$ of $y$ such that $\Psi \in C^4(\mathcal E; C^2(V'; U))$. Furthermore, if $g$ satisfies (\ref{detect_smoothness_beta}) then $\Psi^* g$ satisfies (\ref{detect_smoothness_W}) with $V=V'$.
\end{lemma}
\begin{proof}
Naturality of the exponential map implies that 
    \begin{equation*}
\Psi(\epsilon, z) = \exp_{\Phi(\epsilon, y)}^{g(\epsilon)}
\left(\sum_{j=0}^3 z^j (\p_{y^j} \Phi^k)(\epsilon, y)X_k\right).
    \end{equation*}
For small enough neighbourhood $V'$ of $y$, the smoothness $\Psi \in C^4(\mathcal E; C^2(V'; U))$ follows from 
$g \in C^4(\mathcal E; C^3(M;\Sym))$.
As the exponential map is obtained by a solving system of ordinary differential equations, while varying the initial conditions, the assumption (\ref{detect_smoothness_beta})
implies that $\p_\epsilon^\beta \Psi|_{\epsilon = 0} \in C^\infty(V'; \Sym)$ for any $\beta < \alpha$ where $\alpha = (1,1,1,1)$. Thus (\ref{detect_smoothness_W}) holds with $V' =V$.
\end{proof}

\subsection{Concepts from global Lorentzian geometry}

Up to this point the global geometry of $(M,\hat g)$
has played a minor role. When we combine sending, interaction and receiving of singularities in the next section,
we need to take into account certain global aspects.

Let $p, q\in M$. 
We write $p< q$ when $p$ and $q$ can be joined by a future pointing causal curve in $(M,\hat g)$.
When this is the case, 
we define  the time separation function $\tau(p,q)\in [0,\infty)$
to be the supremum of the lengths 
$$
L(\alpha)=\int_0^1 \sqrt{-\hat g(\dot\alpha(s),\dot\alpha(s))}\,ds
$$
of the piecewise smooth
causal curves $\alpha:[0,1]\to M$ from $p$ to $q$. If  the condition $p< q$ does not
hold, we define $\tau(p,q)=0$. 

When $(x,\xi) \in \Sigma$, we define  $\T(x,\xi)$ to
be the length of
the maximal interval on which $\gamma_{x,\xi}:[0,\T(x,\xi))\to M$ is defined.
Recall that $\gamma_{x,\xi}$ is the geodesic of $(M,\hat g)$
with the initial data $(x,\xi_*)$ where $\xi_*$ is the 
vector $\xi_*^j = \hat g^{jk} \xi_k$.
%Below,  to simplify notations, we sometimes use the notation $\gamma_{x,\xi}([0,\infty))$ for the geodesic $\gamma_{x,\xi}([0,\T(x,\xi))$.
When $\xi$ is future pointing,
we define the modified cut locus functions, c.f.\ \cite[Def.\ 9.32]{Beem},
\beq\label{eq: max time}
& &\rho(x,\xi)=\sup\{s\in [0,\T(x,\xi));\ \tau(x,\gamma_{x,\xi}(s))=0\},
%\\\nonumber
%& &\rho_g(x,\xi_-)=\sup\{s\in [0,\T(x,\xi_-));\ \tau(\gamma_{x,\xi_-}(s),x)=0\}.
\eeq
%The point $\gamma_{x,\xi}(\rho(x,\xi))$ is called the cut point on the geodesic  $\gamma_{x,\xi}$.
Using \cite[Th. 9.33]{Beem}, we see that the function $\rho(x,\xi)$ is lower semi-continuous on the globally hyperbolic
Lorentzian manifold $(M,\hat g)$.

\subsection{Probing the background geometry with propagating singularities}

\begin{figure}%[tbhp]
\centering
\includegraphics[width=0.5\textwidth]{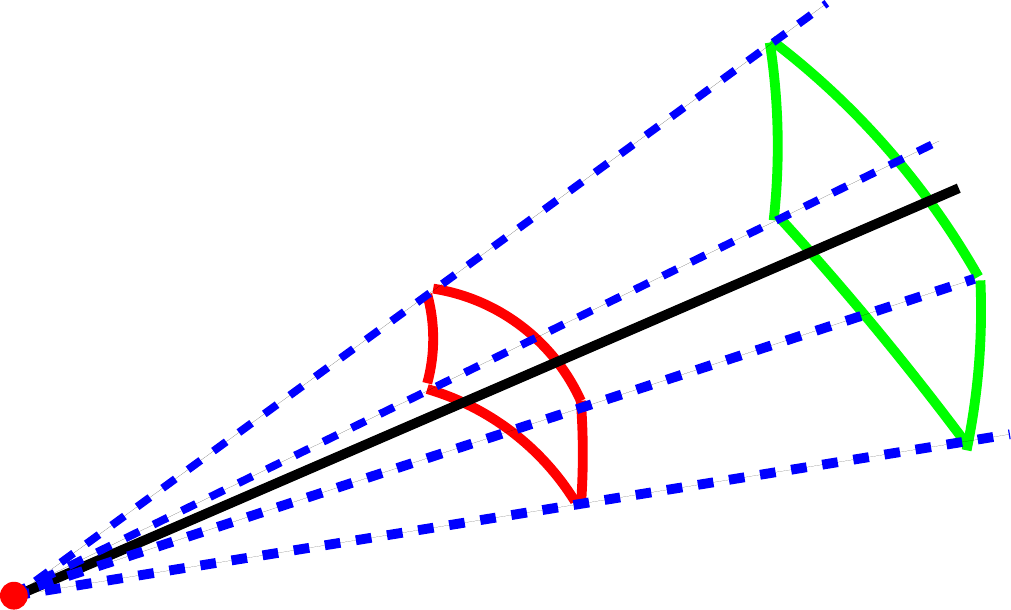}
\caption{
\label{fig_plane_wave}
Schematic of $K(x,\xi; \delta)$ in the Minkowski space.
The projection of $K(x,\xi; \delta)$  at a fixed time is a piece of a sphere. The spherical piece propagates near the lightlike geodesic $\gamma_{x,\xi}$ and is singular on a subset of the light cone emanating from $x$. The black line segment is the
projection of $\gamma_{x,\xi}$ on $\R^3$. 
Two projections of $K(x,\xi; \delta)$ 
on $\R^3$, corresponding to two different times, are drawn 
in red and green.
}
\end{figure}

We will next combine all the steps studied above: 
sending, interaction and receiving of singularities. 
This is done in Theorem \ref{th_detection}, 
that relates analytical information on propagating singularities to purely geometrical information on lightlike geodesics and the causal structure of $(M,\hat g)$.
%
%
%
%Theorem \ref{th_detection} can be roughly summarized by saying that 
%it gives a way to tell apart, using the data (\ref{data}), two geometric situations related to four lightlike geodesics leaving near the world line $\mu$ and a point $y$ near $\mu$.
%
%In the first case, the four geodesics intersect at a single point and $y$ can be joined to the point of intersection by a lightlike geodesic.
%In the second case, the four geodesic do not intersect at a single point, or the point $y$ is not in the causal future of any point where all the four geodesics intersect. 
%The precise formulation in Theorem \ref{th_detection} takes 
%into account several corner cases related to 
%cut points, geodesics returning near $\mu$, 
%and ``spurious'' singularities, in particular, those caused by interactions of singularities propagating on three out the four geodesics. 
%
After establishing the relation between analytical and geometrical information, % in Theorem \ref{th_detection}, 
the final step of the proof of Theorem \ref{th_main} is purely geometrical and coincides with the geometric part of the proof in \cite{KLU-august}.

Consider the Fermi coordinates $\Phi = \Phi_{\hat g}$ defined by (\ref{def_Fermi}). 
Let $x \in V$ and $\xi \in L^+_x V$, 
where the light cone is defined with respect to $\hat g' = \Phi^* \hat g$.
Let $\delta > 0$, denote by $B_x(\xi, \delta) \subset T_{x}^* V$ the ball of radius $\delta$ with respect to 
a fixed auxiliary Riemannian metric on $V$,
and define
\beq \label{notations numbered 3}
K(x,\xi; \delta) = \pi(\Lambda(B_x(\xi, \delta))) \setminus J_{(M, \hat g)}^-(x).
\eeq
Here $\pi : T^* M \to M$ is  
the canonical projection to the base manifold,
$\Lambda$ is the flowout defined by (\ref{def_flowout_S}),
and we have extended the definition as follows
$$
\Lambda(S) = \Lambda(\{(\Phi(x), \Phi_* \xi);\ (x,\xi) \in S\}),
\quad S \subset T^* V.
$$
Observe that $K(x,\xi; \delta)$ is a 3-dimensional $\Lambda$-invariant manifold away from the cut locus of $x$. Moreover, using the notation,
$$
\gamma_{x,\xi} = \gamma_{\Phi(x), \Phi_* \xi},
\quad
\rho(x,\xi) = \rho(\Phi(x), \Phi_* \xi),
\quad
(x,\xi) \in T^* V,
$$
it contains the lightlike geodesic $\gamma_{x,\xi}((0,\rho(x,\xi)))$.

For $z_j \in V$, $\zeta_j \in L_{z_j}^+ V$,
$j=1,2,3,4$, we write $\vec z=(z_j)_{j=1}^4$, $\vec \zeta=(\zeta_j)_{j=1}^4$ and define
    \begin{equation*}
C^+(\vec z,\vec \zeta)= \bigcup_{j=1}^4
J_{(M, \hat g)}^+(\gamma_{z_j,\zeta_j}(\rho(z_j,\zeta_j))),
    \end{equation*}
the causal future of the first cut points on the geodesics $\gamma_{z_j,\zeta_j}$.
Below, it is enough to consider the complement of $C^+(\vec z,\vec \zeta)$, since the cut points are treated in the purely geometric step. 

Denote by $\mathcal K(\vec z,\vec \zeta,\delta)$
the set $\mathcal K$ defined by formula (\ref{bad_set_K})
with 
    \begin{equation*}
K_j = K(z_j,\zeta_j;\delta),\quad \delta > 0,\ j=1,2,3,4.
    \end{equation*}
Recall that $\mathcal K(\vec z,\vec \zeta,\delta)$ is the union of the sets $K_j$ together with a set giving an upper bound for the singular supports of the waves produced by the threefold interactions of waves with singular supports on the sets $K_j$.
When  $\delta\to 0$, the set $\mathcal K(\vec z,\vec \zeta,\delta)$ tends to the set
    \begin{equation*}
\mathcal K_0(\vec z,\vec \zeta)=\bigcap_{\delta>0}\mathcal K(\vec z,\vec \zeta,\delta),
    \end{equation*}
that is of Hausdorff dimension at most two. Below, we consider $\mathcal K_0(\vec z,\vec \zeta)$ as a small exceptional set.

We define
\begin{align*}
\mathbb X(\vec z,\vec \zeta) = \bigcap_{j=1}^4\gamma_{z_j,\zeta_j}((0,\rho(z_j,\zeta_j))).
\end{align*}
Note that $\mathbb X(\vec z,\vec \zeta) \ne \emptyset$
if and only if all the four geodesics $\gamma_{z_j,\zeta_j}$
intersect before their first cut points. 
A piece of information that we want to recover given the data $\mathcal D_{r,k}(\hat g, \hat \phi)$ is 
whether $\mathbb X(\vec z,\vec \zeta) \ne \emptyset$ or not.
Moreover, in the case $\mathbb X(\vec z,\vec \zeta) \ne \emptyset$,
the main piece of information that we want to recover is if $x \in \mathbb X(\vec z,\vec \zeta)$ can be joined to a point $y \in \Phi(V)$ by a lightlike geodesic segment without cut points.
More precisely, we want to determine if $y \in \mathcal E_V(x)$ where 
    \begin{equation}\label{earliest light observation set}
\mathcal E_V(x) = V \cap \mathcal E(x), \quad
\mathcal E(x) = \{\gamma_{x,\xi}(t);\  t \in [0,\rho(x,\xi)],\ \xi \in L^+_x M \}.
    \end{equation}
    
    %{lem_detect}

We denote by $B((\vec z, \vec \zeta), \delta)$,
the ball of radius $\delta > 0$ with respect to a fixed auxiliary Riemannian metric on $(T^* V)^4$,
and use also the shorthand notation $$E = \Sym \oplus \R^L.$$ 
The following two definitions form the basis of a test that allows us to extract geometric information given the data (\ref{data}).

\def\i{\text{in}}
\def\o{\text{out}}
\begin{definition}
\label{def_meas_sing}
Let $p,q \in V$ and suppose that there is open $B \subset \R^4$
such that $J_{(V,\hat g')}(p,q) \subset B$
and $\overline B \subset V$.
Write $V^\i = J_{(V,\hat g')}(p,q)$ and define
    \begin{equation*}
L_4 V^\i = \{ (z_1,\dots, z_4, \zeta_1,\dots,\zeta_4) \in (T^* V^\i)^4;\ 
\zeta_j \in L_{z_j}^+ V,\ z_j \notin J_{(V,\hat g')}^+(z_k),\ j\ne k \}.
    \end{equation*}
Let $\mathcal E \subset \R^4$ be a neighbourhood of the origin.
We say that a family of measurements 
    \begin{equation*}%\label{meas_sing_family}
u(\epsilon) = (g'(\epsilon),\phi'(\epsilon),\mathcal F'(\epsilon)) \in \mathcal D_{r,k}(\hat g, \hat \phi), \quad \epsilon \in \mathcal E,
    \end{equation*}
sends singularities of width $\delta > 0$ from $(\tilde{\vec z}, \tilde{\vec \zeta}) \in L_4 V^\i$ to $y \in V$ if $\mathcal F'(0) = 0$ and the following conditions hold:
\begin{itemize}
\item[(In 1)] 
Writing $(g^\i(\epsilon), \phi^\i(\epsilon), \mathcal F^\i(\epsilon)) = u^\i(\epsilon) = \Psi^\epsilon_* u(\epsilon)$
where $\Psi^\epsilon$ is the solution of (\ref{wave_map_V}) with $g' = g'(\epsilon)$, $\epsilon \in \mathcal E$,
it holds that 
$\supp(\Psi^\epsilon_* \mathcal F'(\epsilon)) \subset V^\i$.
\item[(In 2)]
The family $u^\i$ is $C^4$-smooth with respect to $\epsilon$,
(\ref{no_source_interactions}) holds with $F = \mathcal F^\i$, and writing $f_j=\p_{\epsilon_{j}} \mathcal F^\i(\epsilon)|_{\epsilon = 0}$, it holds that 
    \begin{equation*}
\supp(f_j)\cap J_{(V,\hat g')}^+(\supp(f_k))=\emptyset,\quad j\ne k.
    \end{equation*}
\item[(In 3)]
There are submanifolds 
$Y_j' \subset V^\i$
of codimension one such that the geodesic $\gamma_{z_j,\zeta_j}$ intersects $\Phi(Y_j')$ transversally, and that
    \begin{equation*}
f_j
\in 
I(Y_j' \cap K_{j}'; E),\quad
j = 1,2,3,4,
    \end{equation*}
where $K_j' = \Phi^{-1}(K(\tilde z_j, \tilde \zeta_j;\delta))$, $\tilde{\vec z}=(\tilde z_j)_{j=1}^4$ and $\tilde{\vec \zeta}=(\tilde \zeta_j)_{j=1}^4$.
\item[(Out)] 
Writing $g^\o = \Xi^* g'$ where $\Xi$ are the normal coordinates defined by (\ref{norm_coord_Xi}) and $\alpha = (1,1,1,1)$, it holds that $\p_\epsilon^\alpha g^\o$ is not smooth at the origin.
\end{itemize}
\end{definition}

\begin{definition}
\label{def_detect_fun}
Let $r > 0$ and $k \in \N$, and let $L_4 V^\i$ be as in Definition \ref{def_meas_sing}.
We define a singularity detection function 
$$D_{r,k} : V \times L_4 V^\i \to \{0,1\}$$ by 
$D_{r,k}(y, \vec z, \vec \zeta) = 1$
if and only if
there are $k_0 \in \N$ and $r_0, \delta_0 > 0$
such that $k_0 \ge k$, $r_0 \le r$ and that 
for all $0 < \delta < \delta_0$ 
there are 
$(\tilde{\vec z}, \tilde{\vec \zeta}) \in L_4 V^\i \cap  B((\vec z, \vec \zeta), \delta)$,
a neighbourhood $\mathcal E \subset \R^4$ of the origin 
and a family of measurements 
    \begin{equation*}
u(\epsilon) = (g'(\epsilon),\phi'(\epsilon),\mathcal F'(\epsilon)) \in \mathcal D_{r_0,k_0}(\hat g, \hat \phi), \quad \epsilon \in \mathcal E,
    \end{equation*}
that sends singularities of width $\delta$ from $(\tilde{\vec z}, \tilde{\vec \zeta})$
to $y$. {\mmtext We say that singularities are detected in a stable way at the point $y$ if $D_{r,k}(y,\vec z,\vec \zeta)=1$.}
%\begin{itemize}
%\item[(S1)] $\Phi^{-1} \circ \gamma_j$ intersects $Y_j'$ transversally where $\gamma_j = \gamma_{z_j,\zeta_j}$,
%\item[(S2)] $z_j < z$ with respect to $\hat g'$ for all $z \in Y_j'$,
%\item[(S3)] $Y_j' \cap J_{(V,\hat g')}^+(z) = \emptyset$ 
%for all $z \in Y_k'$ and
%$j \ne k$.
%\end{itemize}
\end{definition}

\begin{theorem}
\label{th_detection}
Suppose that (ND) holds with $U=\overline{\Phi(V)}$.
Let $r > 0$, $k \in \N$ and let $D_{r,k}$ be the singularity detection function. Then the closure $S(\vec z, \vec \zeta)$ of the set $\{\Phi(y);\ y \in V,\ D_{r,k}(y,\vec z, \vec\zeta) = 1\}$
has the following two properties:
\begin{itemize}
\item[(i)] If $\mathbb X(\vec z, \vec \zeta) \ne \emptyset$
then 
$\mathcal E_V(x) \subset S(\vec z, \vec \zeta) \subset J_{(M,\hat g)}^+(x)$.
\item[(ii)] 
If $\mathbb X(\vec z, \vec \zeta) = \emptyset$
then 
$S(\vec z, \vec \zeta) \subset C^+(\vec z,\vec \zeta) 
\cup \mathcal K_0(\vec z,\vec \zeta)$.
\end{itemize}
Furthermore the function $D_{r,k}(y,\vec z,\vec \zeta)$ is determined by the data $\mathcal D_{r,k}(\hat g, \hat \phi)$.

\end{theorem}

Before entering in the rather long proof, let us briefly describe its outline. We send four waves propagating near the four geodesics $\gamma_{z_j,\zeta_j}$. If those geodesics intersect at $x$, we can make a small perturbation in $z_j,\zeta_j$, so that the four geodesics intersect still at $x$, and choose the linearized sources $f_j$, as in Definition \ref{def_meas_sing}, so that the fourfold interaction of the corresponding linearized solutions at $x$ has a principal symbol, in the sense of Proposition \ref{prop_nonvanish_symb}, that does not vanish in the direction towards the point $y$. Moreover, by varying the symbols of $f_j$, we can produce a six dimensional space of singularities at $y$, and at least some of those singularities are observed in the normal coordinates $\Xi$. 

\begin{proof}
It follows from Lemma \ref{lem_wave_coord} 
that the data $\mathcal D_{r,k}(\hat g, \hat \phi)$
determines whether a family 
a family of measurements 
    \begin{equation}\label{detect_pf_family}
u(\epsilon) = (g'(\epsilon),\phi'(\epsilon),\mathcal F'(\epsilon)) \in \mathcal D_{r_0,k_0}(\hat g, \hat \phi), \quad \epsilon \in \mathcal E,
    \end{equation}
with $r_0 \le r$ and $k_0 \ge $,
satisfies (In 1)--(In 3).
Clearly, also the validity of (Out) can be verified by using the data. Thus the function $D_{r,k}(y,\vec z,\vec \zeta)$ is determined by the data $\mathcal D_{r,k}(\hat g, \hat \phi)$.

{\mmtext We note that since $\mathcal E_{0}(x) = \{\gamma_{x,\xi}(t);\  t \in (0,\rho(x,\xi)),\ \xi \in L^+_x M \}$
is a smooth 3-dimensional manifold and the Hausdorff dimension of  $\mathcal K_0(\vec z,\vec \zeta)$ is at most 2,
the set $\mathcal E_{0}(x)\setminus \mathcal K_0(\vec z,\vec \zeta)$ is dense in $\mathcal E_V(x)$.}

%    \begin{equation*}
%\mathcal E_V(x) = V \cap \mathcal E(x), \quad
%\mathcal E(x) = \{\gamma_{x,\xi}(t);\  t \in [0,\rho(x,\xi)],\ \xi \in L^+_x M \}.
%    \end{equation*}
%    
%    %{lem_detect}

The properties (i) and (ii) follow, by passing to the closure,
after we show that for any $y \in V$ 
satisfying 
    \begin{equation}\label{y_assum}
 \Phi(y) \notin C^+(\vec z,\vec \zeta) 
\cup \mathcal K_0(\vec z,\vec \zeta)   
    \end{equation}
the following two implications hold:
(1) if $\mathbb X(\vec z, \vec \zeta) = \emptyset$ or (a) then $D(y,\vec z,\vec \zeta) = 0$; (2) if (b)
then $D(y,\vec z,\vec \zeta) = 1$, where
\begin{itemize}
\item[(a)]
For any point $x \in \mathbb X(\vec z, \vec \zeta)$ it holds that $\Phi(y) \notin J_{(M,\hat g)}^+(x)$.
\item[(b)]
There are $x \in \mathbb X(\vec z, \vec \zeta)$, $\xi \in L^+_x M$ and $t\in (0, \rho(x,\xi))$
 such that ${\Phi(y)=\gamma_{x,\xi}(t)}$.
\end{itemize}

%{\mmtext We note that since $\mathcal E_{in}(x) = \{\gamma_{x,\xi}(t);\  t \in (0,\rho(x,\xi)),\ \xi \in L^+_x M \}$
%is a smooth 3-dimensional manifold and $\mathcal K_0(\vec z,\vec \zeta)$ has Hausdorff dimension 2,
%the set $\mathcal E_{in}(x)\setminus \mathcal K_0(\vec z,\vec \zeta)$ is dense in $\mathcal E_V(x)$.
%Hence, analogously to the proof of Lemma \ref{lem_detect}, we see by using the naturality of the exponential map 
% that $ S(\vec z, \vec \zeta) \subset J_{(M,\hat g)}^+(x)$.}

%    \begin{equation*}
%\mathcal E_V(x) = V \cap \mathcal E(x), \quad
%\mathcal E(x) = \{\gamma_{x,\xi}(t);\  t \in [0,\rho(x,\xi)],\ \xi \in L^+_x M \}.
%    \end{equation*}
%    
%    %{lem_detect}

To establish the implication (1), we show for large $k_0 \in \N$ and small $\delta_0,r_0 > 0$, that if the family (\ref{detect_pf_family})
satisfies
(In 1)--(In 3), with $(\tilde{\vec z}, \tilde{\vec \zeta}) \in L_4 V^\i \cap  B((\vec z, \vec \zeta), \delta)$ and $0 < \delta < \delta_0$, and $\mathbb X(\vec z, \vec \zeta) = \emptyset$ or (a) holds,
then (Out) does not hold. 

We use the shorthand notation $\gamma_j = \gamma_{z_j,\zeta_j}$ and
$$
\mathcal N = 
\hat M \setminus \left(
C^+(\vec z,\vec \zeta) \cup 
\bigcup_{j=1}^4 J_{(M, \hat g)}^-(\Phi(z_j)) \right).
$$
Note that the geodesics $\gamma_j, \gamma_k$ with $j \ne k$
can intersect only once in $\mathcal N$.
Indeed, 
if $\gamma_j$ and $\gamma_k$ intersect at points $x_1 < x_2$,
then the existence of the broken causal geodesic between $z_j$ and $x_2$, 
obtained by switching from $\gamma_j$ to $\gamma_k$
at $x_1$, implies that $\tau(z_j, x_2) > 0$.
This is a contradiction with the fact that $x_2 = \gamma_j(s)$
for some $s < \rho(z_j,\zeta_j)$.

An analogous argument shows that $\gamma_j$ can not have any self-intersections in $\mathcal N$.
In particular, if $\tilde{\mathcal N} \subset \mathcal N$
is open and if its closure is contained in $\mathcal N$,
then for small $\delta$, writing $K_j = K(z_j,\zeta_j;\delta)$, the intersection $K_j \cap \tilde{\mathcal N}$ is a smooth manifold.
The same is true for $\tilde K_j = K(\tilde z_j,\tilde \zeta_j;\delta)$ when $(\tilde z_j, \tilde \zeta_j) \in B((z_j, \zeta_j), \delta)$ for small $\delta$.
Hence the results in Section \ref{sec_interaction} can be applied on $\tilde{\mathcal N}$.
Observe also that when considering a fourfold interaction at $x \in \tilde{\mathcal N}$ the condition (\ref{x_cond_sym}) holds, since there can not exist an earlier threefold interaction.

For large $k_0$ and small $r_0$, Lemma \ref{lem_wave_coord} implies
that the family (\ref{detect_pf_family}) can be also written as 
$$
u(\epsilon) = \Phi_{g(\epsilon)}^* (g(\epsilon),\phi(\epsilon), \mathcal F(\epsilon)), \quad
\Psi^\epsilon_* \mathcal F'(\epsilon) = \Phi^* \mathcal F(\epsilon),
$$
where $(g(\epsilon),\phi(\epsilon),I_{g(\epsilon)} \mathcal F(\epsilon))$
solves (\ref{einred_full}).
We define $v(\epsilon)$ by (\ref{def_v}),
and use the notation $Y' = \bigcup_{j=1}^4 Y_j'$ and $Y = \Phi(Y')$ 
where $Y_j'$ is as in (In 3).
We write also $\alpha = (1,1,1,1)$ and
$$
\mathbb K_\beta = \singsupp(\p_\epsilon^\beta v|_{\epsilon = 0}).
$$
Now (In 3) implies that for small $\delta$,
$$
v_j = \p_{\epsilon_j} v|_{\epsilon = 0}
\in I(\tilde K_j;E) \quad \text{in $\tilde{\mathcal N} \setminus Y$}.
$$
Moreover, 
Lemmas \ref{lem_2_waves} and \ref{lem_3_waves}
imply that 
$\mathbb K_\beta \subset \mathcal K(\tilde{\vec z},\tilde{\vec \zeta},\delta)$ for $\beta < \alpha$.

Let $x \in \mathcal N$ be a point where 
a pair of geodesics $\gamma_j$, $\gamma_k$, $j \ne k$, intersect.
There is at most one such a point,
since any pair of geodesics can intersect only once in $\mathcal N$.
If $\mathbb X(\vec z, \vec \zeta) = \emptyset$, then for small $\delta$,
one of the distributions $v_j$, $j=1,2,3,4$,
is smooth near $x$.
Thus 
Lemma \ref{lem_F4_WF_bound}
implies that $\mathbb K_\alpha \subset \mathcal K(\tilde{\vec z},\tilde{\vec \zeta},\delta)$.
As $y$ satisfies (\ref{y_assum})
it holds for small enough $\delta$ that $\Phi(y) \notin \mathbb K_\alpha$.
On the other hand if (a) holds, then for small $\delta$
it holds again that $\Phi(y) \notin \mathbb K_\alpha$ since $\mathbb K_\alpha$ lies near 
$$
J_{(M, \hat g)}^+(x) \cup \bigcup_{j=1}^4 \gamma_{z_j,\zeta_j}(\R),
$$
and $\bigcup_{j=1}^4 \gamma_{z_j,\zeta_j}((0,\infty)) \subset \mathcal K_0(\vec z, \vec \zeta)$.

In particular, in both the cases $\mathbb X(\vec z, \vec \zeta) = \emptyset$ and (a), it holds for small $\delta$ that 
$$
\Phi(y) \notin \singsupp(\p_\epsilon^\beta g|_{\epsilon = 0}), \quad \beta \le \alpha.
$$

%Using  
%analogously to the proof of Lemma \ref{lem_detect},} % that $ S(\vec z, \vec \zeta) \subset J_{(M,\hat g)}^+(x)$
%Analogously with the proof of Lemma \ref{lem_detect_normal_coords}, 
Using the naturality of the exponential map, {\mmtext analogously with the proof of Lemma \ref{lem_detect_normal_coords},}
we see that $\p_\epsilon^\alpha g^\o$ is smooth near the origin.
This concludes the proof of the implication (1).

Let us now turn to the implication (2) and suppose that (b) holds. 
We choose spacelike submanifolds 
    \begin{equation*}
Y_j' = Y'(z_j,\zeta_j) \subset V^\i, \quad j=1,2,3,4,
    \end{equation*}
of codimension one satisfying the following:
\begin{itemize}
\item[(S1)] $\gamma_j$ intersects $Y_j$ transversally where $Y_j = \Phi(Y_j')$,
\item[(S2)] $z_j < z$ with respect to $\hat g'$ for all $z \in Y_j'$,
\item[(S3)] $Y_j' \cap J_{(V,\hat g')}^+(z) = \emptyset$ 
for all $z \in Y_k'$ and
$j \ne k$.
\end{itemize}
Let 
\begin{align}\label{sigmaj_choice}
\sigma_j \in S^p(Y_j \cap K_j; \SSym \oplus \R^L),
\quad j=1,2,3,4,
\end{align}
and define the map
\begin{align}\label{actual_source}
F(\epsilon, g,\phi) = I_g A_{(g,\phi)}(\mathbb F \sum_{j=1}^4 \epsilon_j \sigma_j), \quad \epsilon = (\epsilon_1,\epsilon_2,\epsilon_3,\epsilon_4),
\end{align}
where $\mathbb F$ and $A_{(g,\phi)}$ are as in Proposition \ref{prop_muloc_linstab}.
We choose the neighbourhood $U$ in Proposition \ref{prop_muloc_linstab} so that $U = U_1 \cup \dots \cup U_4$
where $\overline Y_j \subset U_j$,
$U \subset \Phi(V^\i)$,
and that the 
following perturbed version of (S3) holds:
\begin{itemize}
\item[(S3')] $U_j \cap J_{(V,\hat g')}^+(z) = \emptyset$ 
for all $z \in U_k$ and
$j \ne k$.
\end{itemize}
For any $k_0 \in \N$, the quasilinear wave equation (\ref{einred_full})
with the right-hand side $(F^1,F^2) = F(\epsilon, g,\phi)$ 
has a unique solution 
$$
(g(\epsilon), \phi(\epsilon)) \in C^\infty(\mathcal E; C^{k_0}(\hat M; \Sym \otimes \R^L))
$$ 
when $p \in \R$ is negative enough
and $\mathcal E \subset \R^4$ is a small neighbourhood of the origin.
Note that the neighbourhood $\mathcal E$ depends on the symbols $\sigma_j$, $j=1,2,3,4$.

We define 
$$
\mathcal F(\epsilon) = I_{g(\epsilon)} F(\epsilon, g(\epsilon),\phi(\epsilon)).
$$
By (A1) of Proposition \ref{prop_muloc_linstab},
the family 
$(g(\epsilon), \phi(\epsilon), \mathcal F(\epsilon))$
satisfies the compatibility condition (\ref{compcond}), and therefore it solves (\ref{eq1})--(\ref{init_cond}).
Moreover, the support condition in (A1) implies that 
$\supp(\mathcal F(\epsilon)) \subset U \subset \Phi(V^\i)$.
By (A3) the source $\mathcal F(\epsilon)$ depends smoothly on $\epsilon \in \mathcal E$.
For any $r_0 > 0$ there is a neighbourhood $\mathcal E_{r_0} \subset \mathcal E$ of the origin so that $u(\epsilon) = \Phi_{g(\epsilon)}^* (g(\epsilon), \phi(\epsilon),\mathcal F(\epsilon)) \in \mathcal D_{r_0,k_0}(\hat g, \hat \phi)$
when $\epsilon \in \mathcal E_{r_0}$.

Lemma \ref{lem_wave_coord} implies that 
    \begin{equation}\label{family_uin_by_Phi}
u^\i(\epsilon) = \Psi^\epsilon_* u(\epsilon) = \Phi^* (g(\epsilon), \phi(\epsilon),\mathcal F(\epsilon)).
    \end{equation}
Using again the support condition in (A1) we see that $u^\i(\epsilon)$ satisfies (In 1).
The equation (\ref{family_uin_by_Phi}) implies also that $u^\i$ is smooth with respect to $\epsilon$.
Linearity of $A_{(g,\phi)}$ implies that 
\begin{align}
\label{F_pj}
\p_{\epsilon_j} \mathcal F(\epsilon) |_{\epsilon = 0} = A_{(\hat g, \hat \phi)}(E\sigma_j) \in I(Y_j \cap K_j; E).
\end{align}
To see that (In 2) is satisfied, we use (S3') together with arguments analogous to Lemma \ref{lem_adaptive_deriv} in the below appendix.
The condition (In 3) follows from (\ref{family_uin_by_Phi}), (\ref{F_pj}) and (S1).

We will finish the proof of the implication (2) by showing that for all small $\delta> 0$
there are $(\tilde z_j, \tilde \zeta_j) \in B((z_j,\zeta_j),\delta)$
and symbols $\sigma_j \in S^p(Y_j \cap K_j; \SSym \oplus \R^L)$ such that (Out) holds. 
By (b) there is $x \in \mathbb X(\vec z, \vec \zeta)$, and
by perturbing the geodesics $\gamma_j$ slightly we get the following. 
For any small $\delta > 0$,
there are $(\tilde z_j^n, \tilde \zeta_j^n) \in B((z_j,\zeta_j),\delta)$ such that 
the geodesics $\gamma_j^n = \gamma_{\tilde z_j^n, \tilde \zeta_j^n}$, $j=1,2,3,4$,
intersect at the point $x$ and 
the corresponding manifolds $K_j^n$
satisfy $(\vec b^n, \vec s^n)_{n=1}^6 \in \mathcal U$
for some $\vec s^n$
where 
$\vec b^n = (b_1^n,b_2^n,b_3^n,b_4^n)$, 
$N_x^* K_j^n = \R b_j^n$ and $\mathcal U$ is the set in 
Proposition \ref{prop_nonvanish_symb}.
For small $\delta$, the geodesic $\gamma_j^n$
intersects $Y_j$ and we denote by $z_j^n$ the point of intersection. 
Moreover, we write $\xi = \xi_1^n +\dots +\xi_4^n$
where $\xi_j^n \in \R b_j^n \setminus 0$.
Then there is a unique covector $\zeta_j^n \in N_{z_j^n}^* (Y_j \cap K_j^n)$ 
that lies on the bicharacteristic of $\Box_{\hat g, \hat \phi}$
through $(x, \xi_j^n)$.
As discussed in Section \ref{sec_prop_sing_ein},
taking $\sigma_j = \sigma_j^n$ in (\ref{sigmaj_choice}), the map 
$$
\sigma_j^n(z_j^n,\zeta_j^n) \mapsto \sigma[v_j](x,\xi_j^n)
$$
is an isomorphism 
between the fibres $(\SSym_{Y_j \cap K_j^n} \oplus \R^L)_{(z_j^n,\zeta_j^n)}$
and $(\SSym_{K_j^n} \oplus \R^L)_{(x,\xi_j^n)}$.
Thus we can choose the symbols $\sigma_j^n$ so that 
$\vec s^n = (\sigma[v_j](x,\xi_j^n))_{j=1}^4$, $n=1,\dots,6$.

We write 
$$
\tilde E_x = \mathop{span}\{ \sigmaJ(\vec b^n, \vec s^n);\ n=1,\dots,6\},
$$
and $K_5 = K(x,\xi; \delta')$, $\delta' > 0$.
Observe that $K_5$ is a smooth manifold near $\Phi(y)$ since 
$t < \rho(x,\xi)$.
There is a unique covector $\eta \in N_y^* K_5$ that lies on the bicharacteristic of $\Box_{\hat g, \hat \phi}$ through $(x, \xi)$.
Recalling (\ref{interaction_prinsymb}),
the discussion in Section \ref{sec_prop_sing} implies that the map 
$$
\sigma[\I_4(v_1, v_2, v_3, v_4)](x,\xi)
\mapsto \sigma[v_{1234}](y,\eta)
$$
is an isomorphism between the fibres 
$E_x$ and $E_y$.
Hence
it maps the subspace $\tilde E_x$ onto a 6-dimensional 
subspace $\tilde E_y$ of $E_y$.
Moreover, by Remark \ref{rem_restriction_S},
the subspace 
$\tilde E_y$ is contained in $\Sym_y$
since $\tilde E_x$ is contained in $\Sym_x$.
The intersection of $\tilde E_y$
and the $6$-dimensional subspace $\Sym(K_5)_y \subset \Sym_y$
is at least of dimension $2$.
Now Lemmas \ref{lem_detect} and \ref{lem_detect_normal_coords} imply that there is $n=1,\dots,6$ such that (Out) holds when $\sigma_j = \sigma_j^n$ in (\ref{sigmaj_choice}).
\end{proof}

{\mmtext

Let us summarize the above considerations. When the small (2-dimensional)  sets  $\mathcal K_0(\vec z,\vec \zeta)$ are not considered, Theorem \ref{th_detection} essentially states the following:
First, if the geodesics
$\gamma_{z_j,\zeta_j}$, $j=1,2,3,4$, intersect  before their first cut points 
and the intersection point is $x$, then the set $S(\vec z, \vec \zeta)\subset \Phi(V)$ where 
singularities can be detected in a stable way 
satisfies $\mathcal E_V(x) \subset S(\vec z, \vec \zeta) \subset J^+(x)$.
%contains the set $\mathcal E_V(x)$.
The set $\mathcal E_V(x)$,  defined in (\ref{earliest light observation set}), is called the earliest light observation set corresponding to
the source point $x$. %Note that $\mathcal E_V(x)$ does not intersect the past $J^+(p_j)$ of the points $p_j$. 
Second, if the geodesics  $\gamma_{z_j,\zeta_j}$ do not intersect  before their first cut points, % $p_j$,
then singularities are not detected in a stable way outside the sets $J^+(p_j)$ where $p_j=\gamma_{z_j,\zeta_j}(\rho(z_j,\zeta_j))$, $j=1,2,3,4$, are the first cut points.

Roughly speaking, the above means that when we send singular waves
along the geodesics $\gamma_{z_j,\zeta_j}$, and the geodesics intersect at a point $x$ before their cut points, the non-linear interaction creates an artificial point source 
at the point $x$. This point source sends singularities along the light cone emanating from the point $x$ and the earliest arriving singularities are detected in 
$\Phi(V)$. Thus the inverse problem in Theorem \ref{th_main} is reduced to the geometric problem of determining the conformal type of $ I(\mu(0), \mu(1))$ when we are given the sets $S(\vec z, \vec \zeta)$ for all initial directions $(\vec z, \vec \zeta) \in L_4 V^\i$.  
 This geometric problem is solved \cite{KLU-august} where the inverse
problem for the passive observations of the light cones emanating from the point sources is studied.
Thus} the rest of the proof of Theorem \ref{th_main} is purely geometrical
and coincides with the geometric part of the proof in \cite{KLU-august}. As our notations differ slightly from those in \cite{KLU-august}, we spell out the relation to \cite{KLU-august} in the short proof below.

\begin{proof}[Proof of Theorem \ref{th_main}]
Similarly to Definition 4.2 of \cite{KLU-august} we define  the set 
${\mathcal S}_e(\vec z,\vec \zeta)\subset  {S}(\vec z,\vec \zeta)$ 
of the earliest points where singularities are detected in a stable way. 
Lemma 4.4 of \cite{KLU-august} shows that if 
$x \in \mathbb X(\vec z,\vec \zeta)$ then ${\mathcal S}_e(\vec z,\vec \zeta)$ coincides with the earliest observation set $\mathcal E_{V}(x)$ corresponding to the point $x$ and the observation set $\Phi(V)$.
Theorem \ref{th_detection} above, together with  Theorem 4.5 and Lemma 4.4 of \cite{KLU-august}, implies that the data set (\ref{data}) determines uniquely the collection of the earliest light observations sets 
$$
\Phi^{-1}(\mathcal E_{V}(W))=\{\Phi^{-1}(\mathcal E_{V}(x));\ x\in W\}
$$ 
with source points in the chronological diamond $W = I(\mu(0), \mu(1))$.
By \cite[Th. 1.2]{KLU-august} the collection $\mathcal E_V(W)$, or equivalently, the collection $\Phi^{-1}(\mathcal E_{V}(W))$, determines the conformal type of the set $(W,\hat g)$.
\end{proof}

Corollary \ref{cor_to_main} follows in the same way as 
Corollary 1.3 of \cite{KLU-august}.

\section*{Appendix A: Proofs for linearization stability}

The next lemma is well-known but we give a short proof for the convenience of the reader.

\begin{lemma}
\label{lem_from_div_to_wgauge}
Let $(M, \hat g)$ be a background spacetime,
and suppose that a Lorentzian metric $g \in C^3(M; \Sym)$ satisfies the initial condition (\ref{init_cond}), $(M,g)$ is globally hyperbolic, and that $\Sigma_0 = \{0\} \times N$ is a Cauchy surface with respect to $g$.
Define the reduced Einstein tensor by 
$\Ein_{\hat g}(g) = I_g \Ric_{\hat g}(g)$.
Then $\div_g \Ein_{\hat g}(g) = 0$ implies that $H_{\hat g}(g)=0$.
\end{lemma}

\begin{proof}
We use the shorthand notation $H = H_{\hat g}(g)$, $\nabla = \nabla_g$, $\hat \nabla = \nabla_{\hat g}$,
$\nabla_j = \nabla_{\p_j}$ and $\hat \nabla_j = \hat \nabla_{\p_j}$, and define also 
the operator $\flat_g^s$ that
lowers indices and symmetrizes two tensors $S$ as follows,
$$
(\flat_g^s S)_{jk} 
= \frac 12 (g_{jp} S_k^p + g_{kp} S_j^p).
$$
In local coordinates, $\tr_g \flat_g^s \hat \nabla H = \hat \nabla_p H^p$,
whence
$$
\Ein(g)_{jk}-\Ein_{\hat g}(g)_{jk}
=
\frac12(g_{jp}\hat \nabla_k H^p + g_{kp}\hat \nabla_j H^p
- g_{jk} \hat \nabla_p H^p).
$$
As also $\Ein(g)$ has vanishing divergence, 
\begin{align*}
0 
&= 
g^{qj} \nabla_q (g_{jp} \hat \nabla_k H^p + g_{kp}\hat \nabla_j H^p
- g_{jk} \hat \nabla_p H^p)
%\\&= 
%\nabla_p \hat \nabla_k H^p + 
%g_{kp} g^{qj} \nabla_q \hat \nabla_j H^p
%- \nabla_k \hat \nabla_p H^p
\\&= 
 g_{kp} g^{qj} \nabla_q \hat \nabla_j H^p + L_{kp} H^p,
\end{align*}
where $L_{kp} H^p = \nabla_p \hat \nabla_k H^p - \nabla_k \hat \nabla_p H^p$.
Note that $L_{kp}$
is a linear first order differential operator.
We see that $H$ satisfies the linear second order hyperbolic system 
\begin{align}
\label{lin_syst_for_H}
\begin{cases}
g^{qj} \nabla_q \hat \nabla_j H^r + g^{rk} L_{kp} H^p = 0, & \text{in $M$},
\\H^r = 0, & \text{in $(-\infty, 0) \times N$},
\end{cases}
\quad r = 1,2,3,4,
\end{align}
whose principal part coincides with the principal part of 
$\square_g$. 
As $(M,g)$ is globally hyperbolic and $\Sigma_0$ is a Cauchy surface 
for $g$, $H = 0$ on $M$, see e.g. \cite[Cor. 12.14]{Ringstrom}.
\end{proof}

\begin{corollary}
\label{cor_reduced}
Let $(M, \hat g)$ and $g \in C^3(M; \Sym)$ be as in Lemma 
\ref{lem_from_div_to_wgauge}.
Let 
$$
\phi \in C^2(M; \R^L), 
\quad F^1 \in C^1(M; \Sym),
\quad F^2 \in C^0(M; \R^L).
$$
Suppose that $S = (g,\phi,F^1,F^2)$ 
satisfies (\ref{einred_full}).
Then $S$
satisfies (\ref{eq1})--(\ref{eq2})
if and only if the compatibility condition (\ref{compcond})  holds with $\F^1 = I_g F^1$
and $\F^2 = F^2$.
\end{corollary}

\begin{proof}
If $S$ satisfies (\ref{eq1})--(\ref{eq2}) then 
a direct computation shows that 
(\ref{compcond}) holds.
Suppose now that (\ref{compcond}) holds.
We apply $I_g$ on the first equation in (\ref{einred_full}),
and get 
$$
\Ein_{\hat g}(g) = \mathbb T(g,\phi) + \F^1.
$$ 
Now a computation shows that (\ref{einred_full}) and 
(\ref{compcond}) imply 
$$
\div_g \Ein_{\hat g}(g) 
= (\square_g - \mathcal V'(\phi))\phi \cdot d \phi + \div_g \F^1
=  \div_g \F^1 + \F^2 \cdot d \phi = 0.
$$
Lemma \ref{lem_from_div_to_wgauge} implies $H_{\hat g}(g) = 0$,
and therefore
$$
\Ein(g) = \Ein_{\hat g}(g) = T(g,\phi) + \F^1. 
$$
\end{proof}

\begin{proof}[Proof of Lemma \ref{lem_adaptive}]
We omit writing $\phi$ as a subscript in the proof.
We denote by $e_\ell$, $\ell = 1,\dots,L$, the constant frame 
corresponding to the standard basis of $\R^L$,
and define a vector bundle homomorphism 
$A^1 \in C^\kappa(U; \Hom(T^* U, \R^L))$ by
\begin{align*}
A^1(d \phi_{j}) = e_{j}, \quad j=1,2,3,4.
\end{align*}
Then 
$
A^1(d \phi_{j}) \cdot d \phi
= e_{j} \cdot d \phi = d \phi_{j},
$
and $A^1$ satisfies (i).

We define
\begin{alignat*}{2}
A^2(e_{\ell}) &= 0, &&\ell=1,2,3,4,
\\
A^2(e_{\ell})  &= e_{\ell} - \sum_{j=1}^4 c_{\ell j} e_{j}, \quad &&\ell=5,\dots,L,
\end{alignat*}
where $c_{\ell j} \in C^{\kappa}(U;\R)$ 
are the coefficients of $d \phi_{\ell}$
in the frame $d \phi_j$, $j=1,2,3,4$.
Note that for $\ell=5,\dots,L$, 
$$
A^2 (e_{\ell}) \cdot d \phi
= d \phi_{\ell} - \sum_{j=1}^4 c_{\ell j} d \phi_{j}
= 0.
$$
Hence $A^2$ satisfies (ii).

Let us now turn to the property (iii). 
Suppose that $w = \sum_{\ell = 1}^L a_\ell e_{\ell}$.
Then 
\begin{align*}
A^1(w \cdot d \phi)
&= \sum_{\ell = 1}^L a_\ell A^1 (d \phi_{\ell})
= \sum_{j=1}^4 a_j e_{j} + \sum_{\ell = 5}^L \sum_{j=1}^4 a_\ell c_{\ell j} e_{j}, 
\\
A^2(w) &= \sum_{\ell = 5}^L a_\ell e_{\ell} - \sum_{\ell = 5}^L a_\ell \sum_{j=1}^4 c_{\ell j} e_{j},
\end{align*}
and thus $A^1(w \cdot d \phi) + A^2(w) =w$, that is, (iii) holds.
\end{proof}

Note that $A_\phi^1 (v) = 0$ and 
$A_\phi^2 (w) = 0$
outside the supports of $v$ and $w$, respectively,
simply because these maps are vector bundle homomorphisms. 
We need a similar property also for the derivatives of $A_\phi^j$, $j=1,2$, with respect to $\phi$.
We will give a proof only in the case of $j=1$, the other case being similar. 

\begin{lemma}
\label{lem_adaptive_deriv}
Let $\hat \phi$ be as in Lemma \ref{lem_adaptive}.
Let $\mathcal E \subset \R$ be a neighbourhood of the origin, and
let $\phi(\epsilon) \in C^\infty(\mathcal E; C^{\kappa+1}(U;\R^L))$ satisfy $\phi(0) = \hat \phi$.
Then $\p_\epsilon A_\phi^1|_{\epsilon = 0}(v) = 0$
if the supports of $\p_\epsilon \phi|_{\epsilon = 0}$
and $v$ do not intersect.
\end{lemma}
\begin{proof}
Observe that in local coordinates, 
$$
A_{\phi(\epsilon)}^1 
= \begin{pmatrix}
\Phi(\epsilon)^{-1}
\\
0
\end{pmatrix}
$$ 
where $\Phi(\epsilon)$ is the matrix that has the columns $d \phi_j$, $j=1,2,3,4$.
Thus, writing $\hat \Phi = \Phi(0)$, we have
$$
\p_\epsilon A_\phi^1|_{\epsilon = 0}
= \begin{pmatrix}
-\hat \Phi^{-1}\, \p_\epsilon \Phi|_{\epsilon = 0}\, \hat \Phi^{-1}
\\
0
\end{pmatrix},
$$
and the claim follows.
\end{proof}

\section*{Appendix B: Local existence of solutions}\label{subsec: Direct problem}

In this section we outline briefly the results on quasilinear wave equations that we need. The reduced Einstein's equations coupled with scalar fields (\ref{einred_full}) can be written in the abstract form (\ref{ein_abs}), that is,
\begin{align}
\label{ein_abs_repeated}
g^{-1}(D,D) v + A(x,v,Dv) = F.
\end{align} 
Here $v = (g - \hat g, \phi - \hat \phi)$, and when considering sources of the form (\ref{actual_source}) we allow also $F$ to be a function of $v$ and $Dv$, that is, $F = F(x,v,Dv)$.
Recall that a globally hyperbolic Lorentzian manifold is isometric to a product manifold $\R \times N$ with metric of the form 
\begin{equation*}
%\label{gl_hyp_form}
\hat g = -\beta(t,y) dt^2 + \kappa(t,y).
\end{equation*} 
Then (\ref{ein_abs_repeated}) can be written in the form 
\begin{align}\label{ein_abs_glhyp}
\beta(t,y; v) \p_t^2 v = \sum_{j,k=1}^3 \kappa^{jk}(t,y; v) \p_{y^j} \p_{y^k} v + A(t,y; v,Dv) - F(t,y; v,Dv),
\end{align}
where $g(t,y) = -\beta(t,y; v) dt^2 + \kappa(t,y; v)$.
Well-posedness of equations of this form was established in \cite{HKM}.
 The results there are not directly applicable in our case since \cite{HKM} considers perturbations of the initial condition whereas we consider perturbations via a small compactly supported source $F$. 
However, the fixed point argument there can be adapted to our case, and we refer to Appendix B in \cite{preprint} for the details.

Let $\hat T > 0$ and consider the spaces
$$
E^k=\bigcap_{j=0}^k C^j(0,\hat T;H^{k-j}(N; \Sym \oplus \R^L)),\quad s\in \N.
$$
Let $K \subset (0,\hat T) \times N$ be compact.
Suppose that for all large $k \in \N$, 
there is a neighbourhood $U_k$ of the origin in $E^k$
such that for $v \in U_k$ the function
$f(t,y) = F(t,y; v, Dv)$ satisfies 
$\supp(f) \subset K$.
Then for all large $k \in \N$ there is $r > 0$ such that if 
$\norm{f}_{E^{k}} < r$ for all $v \in U_k$ then there is a unique solution $v \in E^k$ of (\ref{ein_abs_glhyp}) satisfying the initial condition (\ref{init_cond}).

Observe that the above result is local in the sense that only a fixed time interval is considered and the source is assumed to be supported in a fixed compact set. 

 {\nntext 
 %On local solvability of the Einstein and the Eistein-scalar field systems in
%less smooth functions spaces, see \cite{Dafermos2,Dafermos1,Ringstrom,Ringstrom2}. 
In the case of Cauchy problem for Einstein's equations, there are also global results available, see e.g. 
\cite{Ch-K,Hintz-Vasy,Li1,Ringstrom-futurestab,Svedberg-futurestab}.
Above we paid also no attention to the size of the smoothness index $k$. 
For results on the Cauchy problem for
% \HOX{References Choquet-Bruhat2006,
%Choquet-Bruhat2006,
%Choquet-Bruhat2007a are replaced by other references.}
the Einstein and the Eistein-scalar field systems
in low regularity settings see \cite{Dafermos2,Dafermos1,Ringstrom,Ringstrom2}.}

%In the case of Cauchy problem for Einstein's equations, there are also global results available, see e.g. \cite{Chrusciel2005}. Above we paid also no attention to the size of the smoothness index $k$. 
%For results on the Cauchy problem for Einstein's equations in low regularity settings see \cite {Choquet-Bruhat2006, Choquet-Bruhat2007a}.

\section*{Appendix C: On active measurements}

By an active measurement
we mean a model where we can control some of the physical fields  and
the other  physical  fields adapt to the changes in all the fields so that the  conservation law holds. Roughly speaking,
we can consider measurement devices as a complicated
process that changes one energy form to other forms of energy, like a system
of explosives  that transform some potential energy to kinetic energy.
This process creates a perturbation of the metric, a gravitational wave, and it also perturbs the matter fields.
We can then observe these perturbations in a subset of the spacetime.

There are several proposals on generation of large enough gravitational waves so that they can in principle be detected \cite{Chapline1974,Fuezfa2016,Portilla2001,Weber1960}.
Also, gravitational waves from two black holes, in close orbit around each other, were recently observed \cite{Abbott2016}, however, no artificially generated gravitational waves have been measured at the present time. 
In this paper,  our aim has  been to consider
a mathematical model that can be rigorously analyzed. 

Let us take $\mathcal V_\ell(\phi)= \frac 12 m^2 \phi_\ell^2$
for simplicity, and let $g$ and $\phi$ satisfy Einstein's equations coupled with scalar fields
\begin{align}\label{eq: adaptive model}
&\Ein_{\hat g}(g) =P+{\mathbb T}(g,\phi),\hspace{-5mm}\\ \nonumber
&\square_g\phi_\ell -m^2\phi_\ell=S_\ell 
\quad
\hbox{in }\hat M,\quad \ell=1,2,3,\dots,L,
\\ \nonumber
&S_\ell=Q_\ell+{\mathcal S}_\ell^{2nd}(g,\phi, Q, P),\quad
\hbox{in } \hat M,
 \\ \nonumber
& g=\hat g,\quad \phi_\ell=\hat \phi_\ell,\quad
\hbox{in } (-\infty,0) \times N.
\end{align}
Here, in comparison to (\ref{eq1})--(\ref{eq2}), we have simply written $P = \mathcal F^1$ and ${\mathcal S}_\ell^{2nd} = \mathcal F^2_\ell - Q_\ell$, where $Q_\ell$ is arbitrary.
This is to emphasize that $P$ and $Q_\ell$ are considered as primary sources produced by the measurement devices, and the functions ${\mathcal S}^{2nd}_\ell(g,\phi, Q, P)$ are secondary sources that transform energy and momentum from the $\phi$ fields to the sources $Q_\ell$ and $P$. Mathematically speaking the secondary sources cause an effect that guarantees the conservation law (\ref{conservation_law}) is satisfied. 

We allow ${\mathcal S}_\ell^{2nd}$ to depend also on the first derivatives of $(g,\phi, Q, P)$ but hide this in the notation.
The precise construction of ${\mathcal S}^{2nd}_\ell$ is quite technical, see (\ref{def_E_and_A}), but these functions can be viewed as the instructions on  how to build a device that can be used 
to measure the structure of the spacetime far away. 
Put differently they can be seen as recipe to transform energy between the gravitational field and the scalar fields.
We note that if the primary sources satisfy $P=\epsilon f_1$ and $Q=\epsilon f_2$ and the pair $(f_1,f_2)$ satisfies the linearized conservation law (\ref{compcond_lin}), then the first $\epsilon$ derivative of the secondary source is zero.

Above, the source term $P$ can be written  in the form   
\begin{equation}  
\label{P_fluids}  P_{jk}(x)=\sum_{\kappa=1}^K  \mu_\kappa(x) v_j^\kappa(x) v_k^\kappa(x),
\end{equation}  
where $v^\kappa$ are timelike covector fields,  $\mu_\kappa$ real-valued functions and $K\geq 10$.  When $\mu_\kappa\geq 0,$ this decomposition corresponds to non-interacting fluids   with densities $\mu_k$, 4-velocity vectors $v^\kappa$,  and zero pressures (i.e.\ particle or dust flows).   Recall that above   we have considered the case when the source field  $P$ in (\ref{eq: adaptive model})  vanishes for the background solution $(\widehat g,\widehat \phi)$.   This condition can be relaxed, for instance, we can assume that   $(\widehat g,\widehat \phi)$ satisfy (\ref{eq: adaptive model}) with background source $\widehat P$ for which the conservation law $\div_{\hat g}(\widehat P)=0$ holds.  When we assume that $\widehat P$ is written using suitable timelike fields $\hat v^\kappa$, $\kappa=1,2,\dots,K$, $K\geq 10$ and densities $\hat \mu_\kappa$  that are positive in $V$,  any small symmetric perturbation $P$  of $\hat P$ in the set $V$ can be obtained using the fluids $(\mu_\kappa, \hat v_\kappa)$  such that the perturbed densities $ \mu_\kappa$ are positive.   The inverse  problem studied in   Theorem \ref{alternative main thm Einstein} can be  modified to the setting where the background source $\widehat P$  in given in the form (\ref{P_fluids})  with fluids  $(\hat \mu_\kappa, \hat v_\kappa)$ that are non-vanishing  near the geodesic  $\mu$ and the measurements are implemented by perturbing these fluids, but the detailed considerations of this model are outside the scope of this paper.

\bigskip
\paragraph{\bf Acknowledgements.} The authors express their gratitude to MSRI, the Newton  
Institute, the Fields Institute and  
the Mittag-Leffler Institute, where parts of this work have been done.
YK and LO were partly supported by EPSRC.
ML was partly supported by the Finnish Centre of
Excellence in Inverse Problems Research 2012-2017.
GU was partly supported  by NSF, a Clay Senior Award at MSRI,  a Chancellor Professorship at UC Berkeley,  
a Rothschild Distinguished Visiting Fellowship at the Newton Institute, the Fondation de Sciences Math\'ematiques de Paris,
FiDiPro professorship,
and a Simons Fellowship.

\end{document}